\documentclass[final,3p,times]{elsarticle}
\usepackage[T1]{fontenc}

\usepackage{amssymb,color,booktabs}
\usepackage{amsmath}
 \usepackage{amsthm}
\newtheorem{theorem}{Theorem}
\newtheorem{corollary}{Corollary}
\newtheorem{lemma}{Lemma}
\newtheorem{definition}{Definition}
\newtheorem{remark}{Remark}
\newcommand{\diag}{\operatorname*{diag}}

\usepackage[normalem]{ulem}
\def\pep#1#2{\ifx\detokenize{#1}{}\relax\else{\color{red}\sout{#1}}\fi\ifx\detokenize{#2}{}\relax\else{\color{blue}#2}\fi}

\journal{Journal of Computational Physics}
\begin{document}

\bibliographystyle{elsarticle-num}

\numberwithin{equation}{section}
\begin{frontmatter}



\title{Invariant-region-preserving WENO schemes for one-dimensional multispecies kinematic flow models}


\author[1]{Juan Barajas-Calonge}
\ead{juan.barajas2001@alumnos.ubiobio.cl}
\author[2,5]{Raimund Bürger}
\ead{rburger@ing-mat.udec.cl}
\author[3]{Pep Mulet}
\ead{mulet@uv.es}
\author[1,5]{Luis Miguel Villada\corref{cor}}
\ead{lvillada@ubiobio.cl}

\cortext[cor]{Corresponding author}


\affiliation[1]{organization={GIMNAP-Departamento de Matem{\'{a}}tica, Universidad del Bío-Bío},
	addressline={Avenida Collao 1202}, 
	city={Concepción},
	country={Chile}}
	
\affiliation[2]{organization={Departamento de Ingeniería
		Matem\'{a}tica, Universidad
		de Concepción},
	addressline={Casilla 160-C}, 
	city={Concepción},
	country={Chile}}
	
\affiliation[3]{organization={Departament de Matemàtiques, Universitat de València},
		addressline={Av. Vicent Andrés Estellés s/n}, 
		city={Burjassot},
		country={Spain}}
	
\affiliation[5]{organization={CI$^2$MA, Universidad de Concepción},
		addressline={Casilla 160-C}, 
		city={Concepción},
		country={Chile}}

\begin{abstract}
Multispecies kinematic flow models are defined by systems of $N$ strongly coupled, nonlinear first-order conservation laws, where the solution
is a vector of $N$~partial volume fractions or densities. These models 
 arise in various applications including multiclass vehicular traffic and sedimentation of polydisperse suspensions. The solution vector   should  
 take values in a  set   of physically relevant values 
  (i.e., the components are nonnegative and sum up at most to a given maximum value).  It is demonstrated that this 
   set, the so-called invariant region,  is preserved by numerical solutions produced by a new family of high-order finite volume numerical schemes adapted to 
   this class of models. To achieve this property,  and motivated by [X.\ Zhang, C.-W.\ Shu, On maximum-principle-satisfying high order schemes for scalar conservation laws, 
    J.\    Comput.\  Phys.\ 
229  (2010) 3091--3120], a   pair of linear scaling limiters is applied to a high-order  weighted essentially non-oscillatory (WENO)  polynomial reconstruction to obtain  
  invariant-region-preserving (IRP)  high-order polynomial reconstructions. These reconstructions are  combined
    with a local Lax-Friedrichs (LLF) or Harten-Lax-van Leer (HLL)   numerical flux   to obtain a  high-order numerical scheme 
     for the system of conservation laws. It is proved  that this scheme satisfies an  
       IRP property under a suitable Courant-Friedrichs-Lewy (CFL)  condition. 
        The theoretical analysis is corroborated with  
         numerical simulations   for models of multiclass traffic flow and  polydisperse sedimentation.
\end{abstract}




\end{frontmatter}



\section{Introduction}
\subsection{Scope} \label{subsec:scope} 
This work concerns high-order numerical schemes for spatially one-dimensional systems of $N$ first-order nonlinear conservation laws
\begin{align}\label{eq:main-eq} \begin{split} 
	& \partial_t \Phi+\partial_x \boldsymbol{f}(\Phi)=\boldsymbol{0},\quad 
	\Phi = ( \phi_1, \dots, \phi_N)^{\mathrm{T}}, \quad 
	\boldsymbol{f}(\Phi)\coloneqq \bigl(f_1(\Phi),\dots,f_N(\Phi)\bigr)^{\mathrm{T}};  \\ 
 &  f_i(\Phi)\coloneqq \phi_i v_i(\Phi), \quad i=1, \dots, N; \quad x\in I\coloneqq [0,L]\subset \mathbb{R}, \quad t >0,  \end{split} 
\end{align}
 where $t$~is time, $x$~is the spatial coordinate, and the solution  $\Phi=\Phi(x,t)$ usually denotes a vector of partial concentrations $\phi_1, \dots, \phi_N$  (volume fractions or densities) of a number~$N$ of species. The components of~$\Phi$ 
   should be nonnegative and sum up at most to some maximum value~$\phi_{\max}$ that depends on the physical system under consideration. 
   Consequently, $\Phi$~is   assumed to take  values in 
   the set 
\begin{align*} 
\mathcal{D}_{\phi_{\mathrm{max}}}\coloneqq \bigl\{\Phi = ( \phi_1, \dots, \phi_N)^{\mathrm{T}} \in \mathbb{R}^N :  \phi_1\geq0,\dots, \phi_N\geq0, \; \phi\coloneqq \phi_1+\cdots+\phi_N \leq\phi_{\mathrm{max}} \bigr\}. 
\end{align*} 
The model  \eqref{eq:main-eq} is called ``kinematic'' since it is assumed that the velocities~$v_1, \dots, v_N$ are not determined by additional  
 balance equations but are explicitly given functions of~$\Phi$.   
The system     \eqref{eq:main-eq} is equipped with the initial condition
\begin{equation}\label{eq:ini-cond}
	\Phi(x,0) = \Phi_0(x)\in\mathcal{D}_{\phi_{\mathrm{max}}},\quad x\in I, 
\end{equation}
and either periodic boundary conditions
\begin{align}\label{eq:periodic-boundary}
	\Phi(0,t) = \Phi(L,t), \quad t>0
\end{align}
or zero-flux boundary conditions 
\begin{align} \label{bc:zeroflux} 
	\boldsymbol{f}|_{x=0}= \boldsymbol{0}, \quad \boldsymbol{f}|_{x=L}= \boldsymbol{0},\quad t >0.  
\end{align}
  We assume that there exists a piecewise differentiable function $w=w(\phi)$ such that
  \begin{align}  \label{w-phi-ass} 
    \sum_{i=1}^N  \phi_i v_i ( \Phi) = w( \phi) \boldsymbol{\kappa}^{\mathrm{T}} \Phi, 
   \quad \text{where $\phi= \phi_1 + \dots + \phi_N$ and 
  $\boldsymbol{\kappa}\coloneqq (\kappa_1, \dots, \kappa_N)^{\mathrm{T}}$,}
\end{align} 
where $\kappa_1 \geq \kappa_2 \geq \dots \geq \kappa_N >0$ are parameters and the function~$w$ is assumed to satisfy 
\begin{align} 
& w(\phi) \geq 0 \quad \text{for all $\phi \in [0, \phi_{\max}]$,} \label{wpos}  \\
& w(\phi_{\max} ) =0,  \label{wphimax} \\ 
& \text{$w' (\phi)$ is nondecreasing, i.e., $w'(\phi) \leq w'(\tilde{\phi})$ if $0 \leq \phi \leq \tilde{\phi} \leq \phi_{\max}$.} 
 \label{wconvex} 
\end{align} 
Furthermore, we assume that the system \eqref{eq:main-eq} is hyperbolic on $\mathcal{D}_{\phi_{\max}}$, 
 and if $\lambda_1 ( \Phi) \geq \lambda_2 ( \Phi) \geq \dots \geq \lambda_N ( \Phi)$ are the eigenvalues of the 
  Jacobian matrix   $\boldsymbol{\mathcal{J}}_{\boldsymbol{f}} = \boldsymbol{\mathcal{J}}_{\boldsymbol{f}}(\Phi)$
   at $\Phi \in \mathcal{D}_{\phi_{\max}}$, then it is assumed that there exists a piecewise continuous  function $\psi( \phi)$ such that 
   \begin{align}  \label{psiass} 
    \psi ( \phi) \boldsymbol{\kappa}^{\mathrm{T}} \Phi \leq \lambda_N( \Phi) \quad \text{and} \quad w'( \phi) \geq \psi( \phi)
     \quad \text{for all $\Phi \in \mathcal{D}_{\phi_{\max}}$.}  
   \end{align} The properties  \eqref{w-phi-ass}--\eqref{psiass} related to the function~$w=w(\phi)$ are necessary to prove the 
    IRP property for first-order numerical schemes.

System of PDEs like  \eqref{eq:main-eq} arise in various applications including multiclass vehicular traffic \cite{wong2002multi,benzoni2003populations},  the sedimentation of polydisperse suspensions \cite{burger2000numerical,buerger2002model,bbt05,burger2010hyperbolicity}, or the separation kinetics of a dispersion of two immiscible liquids
under the action of gravity \cite{rosso2001gravity}. 
 We herein focus on  two alternative applications. The first is  multiclass 
     vehicular traffic, where we distinguish between $N$~classes of vehicles 
      differing in preferential (freeway) velocities, and consider  
      a closed road of length~$L$ so the final model is   the initial-boundary value problem 
      \eqref{eq:main-eq}--\eqref{eq:periodic-boundary}.  
     The second application is  sedimentation of a 
 polydisperse suspension of small solid particles  dispersed  in a viscous fluid, where we distinguish between $N$~size classes 
  (that settle at different velocities and hence segregate), 
 the $x$-coordinate  usually is a vertical one aligned with gravity,  and 
  the zero-flux boundary conditions \eqref{bc:zeroflux}  are appropriate. Thus, within   that  application, 
   the initial-boundary value problem \eqref{eq:main-eq}--\eqref{eq:ini-cond}, \eqref{bc:zeroflux}  describes settling of 
    a suspension of initial composition~$\Phi_0$ in a closed column. In both applications we
       seek solutions that satisfy an  {invariant region preservation} (IRP)  property, namely,  we expect that 
        if $\Phi_0(x) \in \mathcal{D}_{\max}$ for all~$x$, then $\Phi (x,t) \in \mathcal{D}_{\max}$ for all~$x$ and~$t>0$.
         (A function $w= w(\phi)$ satisfying   \eqref{w-phi-ass}--\eqref{psiass} 
          can easily found for each of the  multiclass traffic and polydisperse sedimentation models studied herein.)

     In the scalar case~$N=1$, and where the unknown is $\phi= \phi(x,t)$, the 
      governing equation turns into the scalar conservation law $\partial_t \phi + \partial_x f(\phi)=0$, 
     where $f(\phi) = \phi v(\phi)$ is assumed to be a non-negative Lipschitz continuous function with support in $(0, \phi_{\max})$.  
      The periodic and zero-flux  initial-boundary value problems then possess unique entropy solutions $\phi(x,t)$  taking values  in 
       $[0, \phi_{\max}]$ for all $t>0$ provided that $\phi_0(x)\in [0,\phi_{\max}]$ for all~$x\in I$. (This is the IRP property in the 
        scalar case.) 
       However, existing high-order numerical schemes for solving these problems do not necessarily yield numerical values 
        in that interval but produce oscillations with slight under- and overshoots. To  handle this  shortcoming, Zhang and Shu  \cite{zhang2010maximum} proposed a technique 
          to create high-order schemes that do satisfy the~IRP property. For one space dimension,  the idea  
          is based on utilizing a standard two-point monotone numerical flux, which if directly applied to the 
            cell averages of the numerical solution  adjacent to a cell interface would generate an IRP 
            scheme   of first-order accuracy only. The key idea to ensure high-order accuracy consists in 
             evaluating the numerical flux on  high-order reconstructions of the unknown to both sides of the 
              cell interface. These reconstructions are obtained by linear scaling around the cell average of the solution 
               in those neighboring cells in conjunction with a limiter function that controls the extrema of the reconstruction polynomial 
                to ensure that  the  reconstructed solution values lie within $[0, \phi_{\max}]$ (in our setting). This procedure 
                 ensures high-order accuracy in space, and is combined with a strong stability preserving (SSP) Runge-Kutta (RK) time 
                 discretization \cite{gottlieb2001strong,gottlieb2009high} to ensure high-order accuracy in time. The  computation 
                  of the limiter involves  determining the extrema of the reconstruction polynomial on the whole cell; this task is 
                  simplified if these extrema are replaced by evaluations of the polynomial on the 
                   finite set of nodes of the Legendre-Gauss-Lobatto quadrature formula on each interval.

       It is the  purpose of the present work to introduce   high-order finite volume schemes with the IRP property for the general case ($N\geq1$) of the 
        initial-boundary value problems \eqref{eq:main-eq}--\eqref{eq:periodic-boundary}   and 
      \eqref{eq:main-eq}--\eqref{eq:ini-cond},  \eqref{bc:zeroflux}. To put the novelty  into the proper perspective, we 
       mention that  
        first-order finite volume schemes with the IRP property  for at least one of these problems 
         have been studied  in \cite{frid1995invariant,burger2008family,burger2013antidiffusive}  while   high-order WENO 
          approximations for the same  models have been advanced in 
          \cite{burger2007adaptive,burger2011implementation,donat2010secular,burger2013spectral,burger2016polynomial,burger2019entropy}. 
           However, the IRP property is not ensured a priori by any of these available  WENO-based 
           treatments, in which  
            spurious oscillations and  negative solution values can be observed in the numerical results. 
Based on the ideas of \cite{zhang2010maximum,zhang2010positivity,zhang2011maximum},   we herein construct  high-order finite volume numerical schemes 
 for the 
        initial-boundary value problems  \eqref{eq:main-eq}--\eqref{eq:periodic-boundary} and 
      \eqref{eq:main-eq}--\eqref{eq:ini-cond}, \eqref{bc:zeroflux}   
       that do 
      satisfy the IRP property. In fact,   
        the linear scaling limiter advanced in \cite{zhang2010maximum} cannot be applied directly
          to this class of problems since if we denote  by $\smash{\phi_{i,j+1/2}^{\mathrm{L},\mathrm{R}}}$, $i = 1, \dots, N$ the 
           reconstructed values at the cell boundaries, and  supposedly applied  the linear scaling limiter of \cite{zhang2010maximum} to ensure that 
            $\smash{\phi_{i,j+ 1/2}^{\mathrm{L},\mathrm{R}}\in[0,\phi_{\max}]}$ for $i=1,\dots,N$, then it would    not be ensured that 
             $\smash{\phi_{j+1/2}^{\mathrm{L},\mathrm{R}} \coloneqq \phi_{1,j+1/2}^{\mathrm{L},\mathrm{R}}+ \dots +  \phi_{N,j+1/2}^{\mathrm{L},\mathrm{R}} \leq \phi_{\max}}$. Consequently, 
               the corresponding reconstructed solution vector~$\Phi_{j+1/2}^{\mathrm{L},\mathrm{R}}$  would not  belong to~$\mathcal{D}_{\phi_{\max}}$. To solve this issue, 
               that is,   to guarantee that 
               \begin{align} \label{nice-prop} 
               \Phi_{j+1/2}^{\mathrm{L},\mathrm{R}} \in \mathcal{D}_{\phi_{\max}} \quad \text{for all~$j$},
               \end{align}     
                and to point out the novelty, 
               we herein propose a two-step approach: in the first step, a high-order polynomial reconstruction of each 
                component~$\phi_i$ on each cell is applied to ensure that the reconstructed cell interface values for each species are  nonnegative, which is achieved 
                 by a modified version of the linear scaling limiter originally proposed in  \cite{zhang2010maximum}; and in the 
                  second step, another linear scaling limiter is applied to the sum of the reconstruction polynomials 
                   (those of the first step) to ensure that the final  reconstructed concentration values at the cell interfaces 
                    sum up to at most to $\phi_{\max}$.  It is proven that this procedure ensures that \eqref{nice-prop} holds. 
                  
 The cell- and species-wise high-order polynomial reconstructions (to which the above-mentioned two-step procedure is applied in each time step)                
     are chosen as  central weighted essentially  non-oscillatory reconstructions 
       proposed by Levy, Puppo, and Russo in
       \cite{levy1999central}. The advantage of   this 
        choice  is that the underlying reconstructions do not only provide single point values but a complete spatial reconstruction   {with adaptive spatial order} in every time step, which is beneficial for evaluating the 
         extrema of the reconstruction polynomials
         in each cell. However, it is possible to use other alternatives such as  WENO-JS \cite{jiang1996efficient}, MR-WENO \cite{zhu2018new}, 
          WENO-Z \cite{borges2008improved} or   {WENO-AO \cite{arbogast18}} reconstructions (with some modifications)  in the implementation of the limiters. 
          The numerical flux is chosen as the   local Lax-Friedrichs (LLF) or    Harten-Lax-van Leer (HLL) numerical 
             flux \cite{harten1983upstream}. For both   we prove 
            that,  under a  suitable CFL condition, the first-order method applied to either initial-boundary value  problem 
            has the 
            IRP property.  That said, it is possible to utilize alternative numerical fluxes such as the 
              Hilliges-Weidlich (HW)  flux \cite{burger2008family,hilliges1995phenomenological} for the MCLWR model. 
              As stated above,  the fully discrete scheme employs an SSP total variation diminishing  (TVD)  Runge-Kutta time  discretization.

\subsection{Related work}  \label{subsec:related} 
Numerical methods preserving physical properties of multiclass models have  been explored in several works. For instance, Jaouen and 
   Lagouti{\`e}re \cite{jaouen2007numerical}   propose a second-order numerical algorithm for the transport of an arbitrary number of materials that is conservative for the mass of each component, i.e., each mass fraction stays in $[0,1]$ and the sum of all mass fractions does not exceed one. 
    Ancellin et al.\  \cite{ancellin2023extension} advance a volume-of-fluid (VOF) method that guarantees natural properties of the volume fractions of a multi-phase flow. 
     In the same direction, Baumgart  and Blanquart \cite{baumgart2024ensuring} propose a simple method for preserving the sum of mass fractions in transport without penalizing the inert species. Recently, Huang and Johnsen   \cite{huang2024consistent} proposed a general numerical approach, consisting of a consistent limiter and the multiphase reduction-consistent formulation,   to solve the multiphase Euler/phase-field
model for compressible $N$-phase $(N \geq 1)$ flows in a consistent and conservative form. Finally, B\"{u}rger et al.\  \cite{burger2013antidiffusive} 
  proposed a first-order antidiffusive and Lagrangian-remap scheme for the MCLWR traffic model which  has the  IRP property  under certain CFL condition 
   depending on the number of species~$N$. 
    The study of high-order IRP numerical methods  is of interest also for other kinds of problems. For instance, one can find schemes preserving density and pressure for  the Euler equations of compressible gas dynamics \cite{zhang2010positivity,zhang2011positivity,zhang2021high} or numerical methods preserving positivity for the shallow water equations \cite{xing2010positivity,xing2013positivity}. An exhaustive   survey of property-preserving numerical schemes for conservation laws and further  references   
  are provided  in  the recent monograph by  Kuzmin and Hajduk \cite{kuzmin2024property}.

\subsection{Outline of the paper} \label{subsec:outline}   
The remainder of the paper is organized as follows. The multispecies  kinematic flow models  studied 
  are described in Section \ref{sec:models}, including  models of 
  multiclass vehicular traffic  (Section~\ref{subsec:mclwr}) and 
  polydisperse sedimentation (Section~\ref{subsec:mlb})  
   along  with their respective bounds of the eigenvalues  of the Jacobian matrix~$\boldsymbol{\mathcal{J}}_{\boldsymbol{f}} ( \Phi)$  
    of the flux vector~$\boldsymbol{f} ( \Phi)$. These are obtained with the secular equation (see Appendix~A).
     In what follows, we refer to the particular  extension of the LWR  traffic model \cite{LW,R} to the multiclass case 
        summarized in Section~\ref{subsec:mclwr} simply as  ``MCLWR model''
      and 
     to the model of polydisperse sedimentation in  conjunction with the Masliyah-Lockett-Bassoon velocity functions 
      (see Section~\ref{subsec:mlb})  as ``MLB model.''  
     In Section \ref{sec:lfinv}
we first introduce the basic time and space finite volume discretizations (Section \ref{subsec:discretizations}). Then, in Section~\ref{subsec3.2} 
 we introduce  the   first-order LLF method and prove in that it satisfies the invariance of the region $\mathcal{D}_{\phi_{\mathrm{max}}}$ for both MCLWR and MLB   models, under a certain CFL condition,  while in Section \ref{subsec3.4} we do the same for the HLL method. In Section~\ref{sec:invcweno} we 
  briefly describe  WENO  reconstructions for  
  scalar conservation laws (Section~\ref{subsec4.1}) and then, in Section~\ref{subsec4.2},   explain   component-wise WENO reconstructions  systems of conservation laws. 
  Then, in Section~\ref{subsec:zs-limiters}, which is at the core of the present work, we 
  advance  a  modification to Zhang and Shu's limiters \cite{zhang2010maximum} for  multiclass models like \eqref{eq:main-eq}  to ensure that 
   the   reconstruction polynomials have the IRP  property. Then we consider a $G$-point Legendre Gauss-Lobatto quadrature
    formula (for ease of computation of the limiter) and  prove that the resulting 
     high-order WENO finite volume scheme is IRP under a CFL condition that depends on one of the quadrature weights. 
       {Moreover, in Section~\ref{subsec:aol} 
      we adapt the ideas from \cite{zhang2010positivity,zhang2017positivity} to rigorously discuss 
       the accuracy of the limiter, with the conclusion   that the
         modification introduced in Section~\ref{subsec:zs-limiters}  does not affect the order of the reconstruction.} 
      To obtain the correct order of convergence of the method we use an SSP Runge-Kutta time discretization scheme, see  Section~\ref{subsec:rk}. 
      In Section \ref{sec:numexa}  
     we present five examples  in which   the numerical methods introduced  in Section~\ref{sec:lfinv},  
       together with the  invariant region preserving  WENO reconstructions of Section~\ref{sec:invcweno}, 
        are used to  numerically solve \eqref{eq:main-eq} for the MCLWR  and  MLB   models. 
          Finally, in Section \ref{sec:conc} we draw some conclusions and discuss open issues.

\section{Multispecies kinematic flow models}\label{sec:models}

\subsection{Multiclass LWR traffic model}\label{subsec:mclwr} The well-known Lighthill-Whitham-Richards (LWR)  kinematic traffic model 
 \cite{LW,R} 
   describes the evolution of the vehicle density $\phi(x,t)$ on a single-lane road by a scalar conservation law $\partial_t \phi +\partial_x(\phi v(\phi))=0,$ where  the velocity function $v=v(\phi)$ is nonnegative and non-increasing  ($v'\leq 0$). In \cite{wong2002multi,benzoni2003populations} this model is generalized to $N$ classes of vehicles with individual densities $\phi_i(x,t)$, $i=1,\dots,N$. The governing equations of the 
resulting multiclass LWR traffic (MCLWR) model  are \eqref{eq:main-eq}, where the main assumption specific for the MCLWR  model is that
\begin{align*} 
	v_i(\phi) = \beta_i v(\phi), \quad i=1,\dots,N,  \quad \beta_1>\cdots>\beta_N>0,  \quad \phi \coloneqq \phi_1+\cdots+\phi_N, 
\end{align*}
that is, drivers of different classes adjust their speed to the total traffic  $\phi$ through the same
function $v(\phi)$, and $\beta_i$ is the free-flowing speed of vehicles  of class~$i$ on an empty highway. The
behavioral law $\phi \mapsto v(\phi)$ may be taken from standard speed-density relations like the Greenshields (GS) model $v(\phi) = 1-\phi/\phi_{\mathrm{max}}$ \cite{greenshields1935study}, where $\phi_{\mathrm{max}}$ represents a maximal car density, or the  Dick–Greenberg (DG) model \cite{dick1966speed,greenberg1959analysis}  
\begin{equation}\label{eq:dgmodel}
	v(\phi) = \min\{1,-C \ln (\phi/\phi_{\mathrm{max}})\} \quad \text{with a constant $C >0$.} 
\end{equation}

The MCLWR model is  strictly hyperbolic whenever $\phi_i>0$ and $\phi<\phi_{\mathrm{max}}$ 
  \cite{donat2008characteristic,zhang2006hyperbolicity}. The eigenvalues $\lambda_i(\Phi)$ of the  
  Jacobian  matrix of the flux $\boldsymbol{\mathcal{J}}_{\boldsymbol{f}}(\Phi) \coloneqq 
   (\partial f_i ( \Phi) / \partial\phi_j)_{1 \leq i,j \leq N}$ satisfy the interlacing property \eqref{eq:int-prop} with the bounds  
\begin{align}\label{eq:bounds-MCLWR}
	M_1(\Phi) = \beta_N v (\phi)+v'(\phi) \boldsymbol{\beta}^{\mathrm{T}} \Phi, \quad 
	M_2(\Phi)  = \beta_1 v(\phi), 
\end{align}
where $\boldsymbol{\beta} \coloneqq (\beta_1, \dots, \beta_N)^{\mathrm{T}}$.  
We assume that the function~$v(\phi)$ has the following properties: 
	\begin{align} \label{mclwrprop}  
	\begin{split} 
		&v(\phi) 
		>0  \quad  \text{for $0 \leq \phi < \phi_{\max}$,} \quad   v( \phi_{\max} ) =0, \quad  
		v'(\phi)\leq 0  \quad   \text{for $0 \leq \phi \leq  \phi_{\max}$; } \\ 
		& \text{$v' (\phi)$ is nondecreasing, i.e., $v'(\phi) \leq v'(\tilde{\phi})$ if $0 \leq \phi \leq \tilde{\phi} \leq \phi_{\max}$.}   
	\end{split} 
\end{align} 
This occurs, for instance, if $v(\phi)$ is chosen according to the  DG or GS models. The  properties \eqref{mclwrprop} ensure that  for the MCLWR model 
\begin{align}  \label{M1M2ineq} 
	M_1 ( \Phi) <  M_2 ( \Phi) \quad \text{for all $\Phi \in \mathcal{D}_{\phi_{\max}}$},   
\end{align} 
 {and it can easily be shown that for the MCLWR model, we may choose $w(\phi) = v(\phi)$ and $\boldsymbol{\kappa} \coloneqq \boldsymbol{\beta}$  along with $\psi ( \phi) = v'(\phi)$ such that  \eqref{w-phi-ass}--\eqref{psiass} 
  are satisfied.}

\subsection{Polydisperse sedimentation}\label{subsec:mlb}
Polydisperse suspensions consist of small solid spherical particles that belong to a number~$N$ of species 
 that differ in size or density, and which are dispersed in a viscous fluid. Here we assume that 
 all solid particles have the same density~$\rho_{\mathrm{s}}$ and that 
  $D_i$~is the diameter of  particle class (species)~$i$, where   $D_1>D_2>\cdots>D_N$. 
   The sedimentation of such a mixture of given initial concentration $\Phi_0(x)$ in a column of depth~$L$ can be then described by 
    the initial-boundary value problem  \eqref{eq:main-eq}--\eqref{eq:ini-cond}, \eqref{bc:zeroflux}, 
     where $\phi_i$ denotes the local
volume fraction of particle species~$i$. A widely used 
   choice of the velocity functions~$v_i$, which is also supported by experimental evidence,   is due to 
  Masliyah  \cite{masliyah1979hindered} and Lockett and Bassoon  
 \cite{lockett1979sedimentation} (``MLB model''). This model arises from the continuity and linear momentum balance equations for the solid species and the fluid through   constitutive assumptions and simplifications \cite{buerger2002model}. For   equal-density particles, the MLB velocities $v_1(\Phi),\dots, v_N(\Phi)$ are given by 
	\begin{align}  \label{eq:MLB} \begin{split} 
		& v_i(\Phi) \coloneqq C(1-\phi)V(\phi)(\delta_i- \boldsymbol{\delta}^{\mathrm{T}} \Phi),\quad  \text{where}  \quad 
		 \delta_i \coloneqq D_i^2/D_1^2, \quad i=1,\dots,N, \\ &   \boldsymbol{\delta} \coloneqq  (\delta_1=1, \delta_2, \dots,\delta_N)^{\mathrm{T}}, \quad     
		C \coloneqq \dfrac{\rho_{\mathrm{s}}-\rho_{\mathrm{f}}}{18 \mu_{\mathrm{f}}} gD_1^2, \end{split} 
	\end{align}
where $\rho_{\mathrm{f}}$ is fluid density, $g$ is the acceleration of gravity, $\mu_{\mathrm{f}}$ is the fluid viscosity, $\phi \coloneqq  \phi_1 +\cdots+ \phi_N$ is the total solids volume fraction, and $V(\phi)$ is a hindered settling factor that is assumed to satisfy~$V(0)=1$, $V(\phi_{\mathrm{max}})=0$, and $V'(\phi)\leq 0$ for $\phi\in [0,\phi_{\mathrm{max}}]$. A standard choice is the Richardson-Zaki equation \cite{rz54} 
\begin{equation}\label{eq:hind-set-fun}
	V(\phi)= \tilde{V} ( \phi) \coloneqq \begin{cases} (1-\phi)^{n_{\mathrm{RZ}}-2} & \text{if $\Phi\in \mathcal{D}_{\phi_{\mathrm{max}}}$}, \\
	0 & \text{otherwise}, \end{cases}  \quad n_{\mathrm{RZ}}>3.
\end{equation}
Following \cite{burger2016polynomial} and if $\phi_{\max} <1$, we use a ``soft cutoff'' version of \eqref{eq:hind-set-fun} to avoid the discontinuity at  $\phi = \phi_{\mathrm{max}}$, namely 
\begin{equation}\label{eq:hind-set-fun-2}
	V(\phi) \coloneqq  \begin{cases}
		(1-\phi)^{n_{\mathrm{RZ}}-2}& \text{ for }0<\phi<\phi_{*},\\
		\tilde{V} (\phi_{*})+ \tilde{V}'(\phi_{*})(\phi-\phi_{*})&\text{ for }\phi_{*}\leq \phi \leq \phi_{\mathrm{max}}, \\
		0&\text{otherwise,}
	\end{cases} \quad n_{\mathrm{RZ}}>3.
\end{equation}
Here $\tau(\phi)  \coloneqq  \tilde{V}(\phi_{*})+ \tilde{V}'(\phi_{*})(\phi-\phi_{*})$ is the tangent to $\tilde{V}(\phi)$ at $(\phi_{*},\tilde{V}(\phi_{*}))$, where 
$\phi_{*}$ is chosen such that $\tau(\phi_{\mathrm{max}})=0$, i.e., 
\begin{equation*}
	\phi_{*} = \dfrac{(n_{\mathrm{RZ}}-2)\phi_{\mathrm{max}}-1}{n_{\mathrm{RZ}}-3}. 
\end{equation*}

According to \cite{burger2010hyperbolicity}, when $\phi_i>0$ and $\phi<\phi_{\mathrm{max}}$, the MLB model is strictly hyperbolic and the eigenvalues $\lambda_i(\Phi)$ of 
 the Jacobian matrix of the flux $\boldsymbol{\mathcal{J}}_{\boldsymbol{f}}  (\Phi)$ interlace with the velocities $v_i ( \Phi)$, 
  as is expressed by  \eqref{eq:int-prop} and where  
\begin{align}\label{eq:bounds-MLB}
		M_1(\Phi)  = C\bigl(\delta_N V(\phi) + \bigl((1-\phi)V'(\phi)-2V(\phi)\bigr) 
		 \boldsymbol{\delta}^{\mathrm{T}} \Phi \bigr) , \quad 
		M_2(\Phi)  = v_1(\Phi).
	\end{align}
	In what follows we always assume  that 
	\begin{align} \label{mlbprop}  
	 & V(\phi) 
	 >0  \quad  \text{for $0 \leq \phi < \phi_{\max}$,} \quad   V( \phi_{\max} ) =0, \quad  \text{and} \quad 
	V'(\phi)\leq 0. 
	\end{align} 
	These properties hold,  for instance, if $V(\phi)$ is chosen according to \eqref{eq:hind-set-fun}, \eqref{eq:hind-set-fun-2}
 with $0 < \phi_{\max} < 1$. The  properties \eqref{mlbprop} ensure  that   \eqref{M1M2ineq}
  holds  for the MLB model, where $M_1 ( \Phi)$ and $M_2 ( \Phi)$ are  defined by  \eqref{eq:bounds-MLB}.
	 {We may then  choose
	$w(\phi) \coloneqq  C(1-\phi)^2 V(\phi)$ and $\boldsymbol{\kappa} \coloneqq \boldsymbol{\delta}$ 
 along with $\psi( \phi)  \coloneqq C((1- \phi) V'(\phi) - 2 V(\phi))$ 
	such that  \eqref{w-phi-ass}--\eqref{psiass} 
  are satisfied.}

\section{First-order invariant-region-preserving schemes}\label{sec:lfinv}
\subsection{Discretizations}\label{subsec:discretizations}  \label{subsec3.1} 
We first discretize the domain $[0,L]\times[0,T]$. 
 For the spatial interval $[0,L]$, we choose $M \in \mathbb{N}$,  a meshwidth $\Delta x \coloneqq L/ M$,  and define the cell centers $x_{j}\coloneqq (j+1/2)\Delta x$ for $j\in \{0,\dots,M-1\}$ and the cell interfaces $x_{j+1/2} =  (j+1) \Delta x$ for $j\in \{0,\dots,M \} \eqqcolon \mathbb{Z}_M$. With this setup $x_{-1/2}=0$ and $x_{M-1/2}=L$. In this way, we subdivide the interval $[0,L]$ into cells $I_j\coloneqq  [x_{j-1/2},x_{j+1/2})$, $j\in \{0,\dots,M-1\}$.  Similarly, for the time interval $[0,T]$ we select  $N_{T} \in \mathbb{N}$ and a sequence of temporal mesh widths $\Delta t_n$, and defining $t_0\coloneqq 0$ and $t_{n+1}\coloneqq  t_n+\Delta t_n$ for $n\in \{0,\dots,N_T\}$ subject to the condition $\Delta t_0+\cdots+\Delta t_{N_T-1}
   = T$. This leads to time strips $I^n\coloneqq  [t_n,t_{n+1})$, $n\in \{0,\dots,N_T-1\}$. The ratio $\lambda_n\coloneqq \Delta t_n/\Delta x$ is assumed to 
    satisfy  a CFL condition that will be specified later. The numerical schemes  produce an approximation $\Phi_{j}^n\approx \Phi(x_j,t_n)$ defined at the mesh points $(x_j,t_n)$ for $j\in \mathbb{Z}_M$ and $n\in \{0,\dots, N_T\}$.

We then define a marching formula for the solution as
\begin{align} \label{eq:march} 
	\Phi_j^{n+1} = \Phi_j^n -\lambda_n \bigl(\boldsymbol{\mathcal{F}}_{j+1/2}^n-\boldsymbol{\mathcal{F}}_{j-1/2}^n \bigr), \quad j= 0, \dots, M-1, \quad 
	  n = 0,\dots, N_T-1. 
\end{align}
For the  periodic boundary conditions \eqref{eq:periodic-boundary},  we set 
\begin{align} \label{bcperiodicdiscr} 
  \boldsymbol{\mathcal{F}}_{-1/2}^n =  \boldsymbol{\mathcal{F}}_{M-1/2}^n,   
 \end{align} 
 and for the zero-flux boundary conditions \eqref{bc:zeroflux}, 
 \begin{align} \label{bczerofluxdiscr} 
  \boldsymbol{\mathcal{F}}_{-1/2}^n = \boldsymbol{0}, \quad   \boldsymbol{\mathcal{F}}_{M-1/2}^n = \boldsymbol{0}.  
 \end{align} 
The computation of the 
  numerical flux vector $\smash{\boldsymbol{\mathcal{F}}_{j+1/2}^n}$ in all other cases is  described in what follows.  

\subsection{IRP property of the LLF scheme} \label{subsec3.2} 
We now describe  a first-order method that satisfies the  IRP property with respect to  $\mathcal{D}_{\phi_{\mathrm{max}}}$.  We utilize the local  Lax-Friedrichs (LLF) 
 numerical flux  
\begin{equation}\label{eq:llf}
	\boldsymbol{\mathcal{F}}_{j+1/2}^{\mathrm{LLF}}\coloneqq \dfrac{1}{2}\big(\boldsymbol{f}(\Phi_j^n)+\boldsymbol{f}(\Phi_{j+1}^n)-\alpha_{j+1/2}^n(\Phi_{j+1}^n-\Phi_j^n)\big),
	 \quad \text{where} \quad  \alpha_{j+1/2}^n \coloneqq  \max \bigl\{|S_{\mathrm{L},j+1/2}^n|,|S_{\mathrm{R},j+1/2}^n| \bigr\},  
\end{equation}
where $\smash{S_{\mathrm{L},j+1/2}^n}$ and $\smash{S_{\mathrm{R},j+1/2}^n}$ denote  
 lower and upper bounds for the eigenvalues $\lambda_1(\Phi), \dots, \lambda_N( \Phi)$ of 
$\boldsymbol{\mathcal{J}}_{\boldsymbol{f}}(\Phi)$ at the interface $x_{j+1/2}$. These values can be estimated by setting 
\begin{align}
		S_{\mathrm{L},j+1/2}^n &=\min_{0\leq s\leq 1}M_1\bigl(s\Phi_{j+1}^n+(1-s)\Phi_{j}^n\bigr), \label{eq:disc-boundsuuy}  \\  
		S_{\mathrm{R},j+1/2}^n &=\max_{0\leq s\leq 1}M_2\bigl(s\Phi_{j+1}^n+(1-s)\Phi_{j}^n\bigr), \label{eq:disc-bounds} 
	\end{align}
where $M_1$ and $M_2$ are given by \eqref{eq:bounds-MCLWR} or \eqref{eq:bounds-MLB} (depending on the model under study). 

Now we can prove an IRP property  for the first-order scheme \eqref{eq:march}, equipped with the LLF numerical flux \eqref{eq:llf}.  For the proof, we use a slightly smaller bound  for the smallest eigenvalue than the 
one stipulated by \eqref{eq:bounds-MCLWR} for the MCLWR model and by \eqref{eq:bounds-MLB} for the MLB model; namely, we employ 
$\smash{\tilde{M}_1 ( \Phi)\coloneqq\psi(\phi)\boldsymbol{\kappa}^{\mathrm{T}} \Phi}$, where we recall  {from Sections~\ref{subsec:mclwr} and~\ref{subsec:mlb} that} 
\begin{equation}\label{eq:kappafun} 
	\psi(\phi) \coloneqq 
	\begin{cases}  
		v' (\phi)  & \text{for the  MCLWR model},\\
		C[(1-\phi) V'(\phi)-2V(\phi)] &\text{for the MLB model,} 
	\end{cases} 
	\qquad
	\boldsymbol{\kappa} \coloneqq 
	\begin{cases}  
		\boldsymbol{\beta} & \text{for the  MCLWR model},\\
		\boldsymbol{\delta} &\text{for the MLB model.} 
	\end{cases} 
      \end{equation}
Notice that $\psi(\phi) \leq 0$ for $0 \leq \phi \leq \phi_{\max}$ for both models   due to the 
explicit respective assumptions \eqref{mclwrprop} and \eqref{mlbprop}.   For both models, 
  $\smash{\tilde{M}_1 ( \Phi) \leq M_1(\Phi)}$, and  to prove the IRP property for both models and 
 the LLF and HHL schemes we employ the slightly smaller lower  estimate (instead of  \eqref{eq:disc-boundsuuy})
\begin{align}  \label{sl34} 
	S_{\mathrm{L},j+1/2}^n =\min_{0\leq s\leq 1}\tilde{M}_1\bigl(\Phi_{j+1/2}^n (s) \bigr), \quad \text{where} \quad  
	\Phi_{j+1/2}^n (s) \coloneqq s\Phi_{j+1}^n+(1-s)\Phi_{j}^n. 
\end{align}

\begin{theorem}\label{thm:lfinv}
		Consider the LLF scheme defined by the marching formula 
	\begin{align}  \label{llf-march} 
		\Phi_j^{n+1} = \Phi_j^n - \lambda_n\bigl(    \boldsymbol{\mathcal{F}}^{\mathrm{LLF}}_{j+1/2} - \boldsymbol{\mathcal{F}}^{\mathrm{LLF}}_{j-1/2}
		\bigr), 
	\end{align} 
	where the LLF numerical flux is given by \eqref{eq:llf} along with the definitions \eqref{eq:disc-bounds} and \eqref{sl34} of 
	$S_{\mathrm{R},j+1/2}$ and $S_{\mathrm{L},j+1/2}$, respectively. If the CFL condition 
	\begin{align} \label{llfcfl}  
	\alpha	\lambda_n  \leq 1 , \quad \text{\em where} \quad \alpha \coloneqq \max_j \bigl\{  |S_{\mathrm{L},j+1/2}|,  |S_{\mathrm{R},j+1/2}|  \bigr\} 
	\end{align}
	is in effect, then the LLF scheme for  \eqref{eq:main-eq},   with the velocity functions~$v_i(\Phi)$ 
	  {chosen such that there exists a function $w=w(\phi)$ satisfying  \eqref{w-phi-ass}--\eqref{psiass}}, satisfies the invariant region preservation property
	 \begin{align}  \label{eq3.7} 
	  \text{\em for all $n=0, \dots, N_T-1$:}  \quad 
	  \Phi_j^{n}\in \mathcal{D}_{\phi_{\mathrm{max}}} \quad \text{\em for all $j\in\mathbb{Z}_M$} 
	   \Rightarrow   \Phi_j^{n+1}\in \mathcal{D}_{\phi_{\mathrm{max}}} \quad \text{\em for all $j\in\mathbb{Z}_M$}.  
	   \end{align}
\end{theorem}

\begin{proof}
	For simplicity,  in the proof we omit the  index~$n$  in the right-hand side of the marching formula, i.e., we set 
	 $\smash{\Phi_j\coloneqq \Phi_j^n}$, $\smash{\alpha_{j+1/2} \coloneqq  \alpha_{j+1/2}^n}$, and $\lambda \coloneqq  \lambda_n$. Let us  assume that both~$\smash{\boldsymbol{\mathcal{F}}_{j+1/2}^{\mathrm{LLF}}}$ 
	 and~$\smash{\boldsymbol{\mathcal{F}}_{j-1/2}^{\mathrm{LLF}}}$ are given by appropriate versions of~\eqref{eq:llf} (see Remark~\ref{llfrem} below for 
	  the boundary flux vectors defined by \eqref{bcperiodicdiscr} or \eqref{bczerofluxdiscr}). Then   
	\begin{align*}
		\boldsymbol{\mathcal{F}}_{j+1/2}^{\mathrm{LLF}}-\boldsymbol{\mathcal{F}}_{j-1/2}^{\mathrm{LLF}} &=  \dfrac{\alpha_{j+1/2}+\alpha_{j-1/2}}{2} \Phi_j-\dfrac{ \alpha_{j+1/2}}{2}\biggl(\Phi_{j+1}-\dfrac{1}{\alpha_{j+1/2}}\boldsymbol{f}(\Phi_{j+1}) \biggr) -\dfrac{ \alpha_{j-1/2}}{2}\biggl(\Phi_{j-1}+\dfrac{1}{\alpha_{j-1/2}}\boldsymbol{f}(\Phi_{j-1}) \biggr), 
	\end{align*}
	hence by the marching formula \eqref{eq:march} 
	and defining 
		\begin{align*}
 	\boldsymbol{\mathcal{G}}_1 (\Phi_{j+1}) \coloneqq \Phi_{j+1}-\dfrac{1}{\alpha_{j+1/2}}\boldsymbol{f}(\Phi_{j+1}) \quad \text{and} \quad 
			\boldsymbol{\mathcal{G}}_2(\Phi_{j-1})&\coloneqq \Phi_{j-1}+\dfrac{1}{\alpha_{j-1/2}}\boldsymbol{f}(\Phi_{j-1}), 
		\end{align*}
	 we can write 
	\begin{equation}\label{eq:convex}
		\Phi_j^{n+1}=\left(1-\dfrac{\lambda(\alpha_{j+1/2}+\alpha_{j-1/2})}{2}\right) \Phi_j+\dfrac{ \lambda\alpha_{j+1/2}}{2}\boldsymbol{\mathcal{G}}_1(\Phi_{j+1})+\dfrac{ \lambda \alpha_{j-1/2}}{2}\boldsymbol{\mathcal{G}}_2(\Phi_{j-1}). 
	\end{equation}
	Consequently, 
	 the  CFL condition  \eqref{llfcfl} implies that $\smash{\Phi_j^{n+1}}$ can be represented as a convex combination of~$\Phi_j$, $\smash{\mathcal{G}_1(\Phi_{j+1})}$, 
	  and $\mathcal{G}_2(\Phi_{j-1})$.
	 %
	%
    %
	Noting  that  
	\begin{align*}
		\mathcal{G}_{1,i} (\Phi_{j+1}) &  = \phi_{i,j+1}-\dfrac{ f_i(\Phi_{j+1}) }{\alpha_{j+1/2}}
		 =\phi_{i,j+1}\biggl(1-\dfrac{v_i(\Phi_{j+1})}{\alpha_{j+1/2}}\biggr) 
		 \geq \phi_{i,j+1}\biggl(1-\dfrac{|v_1(\Phi_{j+1})| }{\alpha_{j+1/2}}\biggr) 
		\geq \phi_{i,j+1} \biggl(1-\dfrac{|S_{\mathrm{R},j+1/2}| }{\alpha_{j+1/2}}\biggr) \quad \text{and}   \\ 
		\mathcal{G}_{2,i}(\Phi_{j-1}) &  = \phi_{i,j-1}+\dfrac{f_i(\Phi_{j-1})}{\alpha_{j-1/2}} =\phi_{i,j-1}\biggl(1+\dfrac{v_i(\Phi_{j-1})}{\alpha_{j-1/2}}\biggr) 
		 \geq \phi_{i,j-1}\biggl(1+\dfrac{M_1(\Phi_{j-1})}{\alpha_{j-1/2}}  \biggr) 
		 \geq \phi_{i,j-1}\biggl(1+\dfrac{S_{\mathrm{L},j-1/2}}{\alpha_{j-1/2}}  \biggr) 
		\\ &  \geq \phi_{i,j-1}\biggl(1-\dfrac{|S_{\mathrm{L},j-1/2}| }{\alpha_{j-1/2}}  \biggr)
	\end{align*}
	 and recalling  the definition of~$\alpha_{j+1/2}$ in \eqref{eq:llf},  
	  we see that $\mathcal{G}_{1,i} (\Phi_{j+1}) \geq 0$ and~$\mathcal{G}_{2,i}(\Phi_{j-1}) \geq 0$.
	  Since  the coefficients of $\phi_{i,j}$, $\phi_{i,j+1}$, and $\phi_{i,j-1}$ are all nonnegative (this follows in case 
	  of the coefficient of~$\phi_{i,j}$ from the CFL condition), we deduce that if $\phi_{i,j}\geq 0$, $\phi_{i,j+1}\geq 0$, and $\phi_{i,j-1}\geq 0$,  
	  then $\phi_{i,j}^{n+1}  \geq 0$   
	    for all $i\in \{1,\dots,N\}$.

 It remains to prove that if $\smash{\phi_j \coloneqq \phi_j^n \coloneqq   \phi_{1,j}^n   + \dots  + \phi_{N,j}^n  \leq \phi_{\max}}$ for all~$j$, then 
$\smash{ \phi_j^{n+1} \leq \phi_{\max}}$ for all~$j$. Component~$i$, $i \in \{ 1, \dots, N \}$ of \eqref{eq:convex} can be rewritten as
\begin{align}\label{eq11}
	\begin{split}
			\phi_{i,j}^{n+1} &= \left(1- \lambda \dfrac{\alpha_{j+1/2}+\alpha_{j-1/2}}{2}\right) \phi_{i,j}+\dfrac{\lambda \alpha_{j+1/2}}{2}\left(\phi_{i,j+1}-\dfrac{\phi_{i,j+1}v_i(\Phi_{j+1})}{\alpha_{j+1/2}}\right)+\dfrac{\lambda \alpha_{j-1/2}}{2}\left(\phi_{i,j-1}+\dfrac{\phi_{i,j-1}v_i(\Phi_{j-1})}{\alpha_{j-1/2}}\right)\\
			 &=\left(1-\lambda \dfrac{\alpha_{j+1/2}+\alpha_{j-1/2}}{2}\right)\phi_{i,j}+ \dfrac{\lambda \alpha_{j+1/2} }{2}  \phi_{i,j+1}+ \dfrac{\lambda \alpha_{j-1/2} }{2} \phi_{i,j-1}-\dfrac{\lambda}{2} \phi_{i,j+1} v_i(\Phi_{j+1})+\dfrac{\lambda}{2} \phi_{i,j-1} v_i(\Phi_{j-1})
	\end{split}
\end{align}
Summing \eqref{eq11} over $i=1, \dots, N$, we get 
\begin{align} \label{eq12}  
		\phi_j^{n+1} & =\left(1-\lambda \dfrac{\alpha_{j+1/2}+\alpha_{j-1/2}}{2}\right) \phi_{j}+ \dfrac{\lambda \alpha_{j+1/2} }{2} \phi_{i,j+1}+\dfrac{\lambda \alpha_{j-1/2} }{2} \phi_{i,j-1} 
		-   \frac{\lambda}{ 2} 
		\sum_{i=1}^N v_i ( \Phi_{j+1}) \phi_{j+1} 
		+  \frac{\lambda}{ 2} 
		\sum_{i=1}^N v_i ( \Phi_{j-1}) \phi_{j-1}. 
\end{align} 
 {According to \eqref{w-phi-ass}, we can write \eqref{eq12} as  
\begin{align} \label{eq12a}  
		\phi_j^{n+1} & =\left(1-\lambda \dfrac{\alpha_{j+1/2}+\alpha_{j-1/2}}{2}\right) \phi_{j}+ \dfrac{\lambda \alpha_{j+1/2} }{2} \phi_{j+1}+\dfrac{\lambda \alpha_{j-1/2} }{2} \phi_{j-1} 
		-   \frac{\lambda}{ 2} w(\phi_{j+1}) \boldsymbol{\kappa}^{\mathrm{T}} 
		  \Phi_{j+1} 
		+  \frac{\lambda}{ 2}  w(\phi_{j-1})  \boldsymbol{\kappa}^{\mathrm{T}}  
		\Phi_{j-1} . 
\end{align}
In light of \eqref{wpos} we obtain  
\begin{align}  \nonumber 
		\phi_j^{n+1} &  \leq \left(1-\lambda \dfrac{\alpha_{j+1/2}+\alpha_{j-1/2}}{2}\right) \phi_{j}+ \dfrac{\lambda \alpha_{j+1/2} }{2} \phi_{j+1}+\dfrac{\lambda \alpha_{j-1/2} }{2} \phi_{j-1} 
		+  \frac{\lambda}{ 2}  w(\phi_{j-1})  \boldsymbol{\kappa}^{\mathrm{T}}  
		\Phi_{j-1}  \\  \label{eq12a} 
		& = \left(1-\lambda \dfrac{\alpha_{j+1/2}+\alpha_{j-1/2}}{2}\right) \phi_{j}+ \dfrac{\lambda \alpha_{j+1/2} }{2} \phi_{j+1}+\dfrac{\lambda \alpha_{j-1/2} }{2} \phi_{\max} + \mathcal{Y}, 
\end{align}
where we define 
\begin{align} \label{calY1def-llf-new} 
	\mathcal{Y}  & \coloneqq \dfrac{\lambda}{2} \bigl(  \alpha_{j-1/2}  ( \phi_{j-1} - \phi_{\max} ) 
	+    w(\phi_{j-1})  
	\boldsymbol{\kappa}^{\mathrm{T}}\Phi_{j-1} \bigr).
\end{align}
We assume that $\phi_j,\phi_{j+1} \in [0, \phi_{\max}]$. If $\phi_{j-1} =
 \phi_{\max}$ then $\mathcal{Y}=0$, so we  deduce from \eqref{wphimax}, by  a convex combination argument, that $\phi_j^{n+1} \leq \phi_{\max}$. To handle the remaining cases, i.e. $0\leq \phi_{j-1}<\phi_{\max}$, we divide \eqref{calY1def-llf-new} by $ \phi_{j-1} - \phi_{\max}$ to obtain
  \begin{align}  \label{huzzut-llf} 
 	\frac{\mathcal{Y}}{  \phi_{j-1} - \phi_{\max} } & = \dfrac{\lambda \alpha_{j-1/2} }{2} +  \frac{\lambda }{ 2}  \frac{w(\phi_{j-1})}{\phi_{j-1}- \phi_{\max}}  \boldsymbol{\kappa}^{\mathrm{T}} \Phi_{j-1} . 
 \end{align}  
 Notice that there exists a number $\xi_{j-1} \in [ \phi_{j-1}, \phi_{\max}]$ such that 
 \begin{align*} 
 	\frac{w(\phi_{j-1})}{\phi_{j-1}- \phi_{\max}} =  \frac{w(\phi_{j-1}) - w( \phi_{\max}) }{\phi_{j-1}- \phi_{\max}}
 	= w' ( \xi_{j-1})  \geq w' ( \phi_{j-1}), 
 \end{align*} 
where the last inequality holds due to \eqref{wconvex}.  Consequently, from \eqref{huzzut-llf} and 
 $\alpha_{j-1/2} \geq  | S_{\mathrm{L},j-1/2} | \geq   - S_{\mathrm{L},j-1/2}^-$ 
we get
  \begin{align}  \label{kjjy-llf} 
 	\frac{\mathcal{Y}}{  \phi_{j-1} - \phi_{\max} } \geq 
 	\frac{  \lambda }{  2} \bigl( - S_{\mathrm{L},j-1/2}^- +  w'( \phi_{j-1} )  \boldsymbol{\kappa}^{\mathrm{T}} \Phi_{j-1} \bigr). 
 \end{align} 
 Since 
 \begin{align}  \label{SLineq1} 
  S_{\mathrm{L}, j-1/2}^- \leq  S_{\mathrm{L}, j-1/2} = \min_{ 0 \leq s \leq 1} 
    \psi \bigl( \phi_{j-1/2} (s) \bigr) \boldsymbol{\kappa}^{\mathrm{T}} \Phi_{j-1/2} (s) 
     \leq \psi ( \phi_{j-1}) \boldsymbol{\kappa}^{\mathrm{T}} \Phi_{j-1} \leq w'( \phi_{j-1} ) 
      \boldsymbol{\kappa}^{\mathrm{T}} \Phi_{j-1}, 
 \end{align} 
 we get that $\smash{-  S_{\mathrm{L}, j-1/2}^-  \geq - w'( \phi_{j-1} ) 
      \boldsymbol{\kappa}^{\mathrm{T}} \Phi_{j-1}}$. Thus, after  multiplying \eqref{kjjy-llf} with $\phi_{j-1} - \phi_{\max} <0$, 
       we get $\mathcal{Y} \leq 0$ and therefore, from  \eqref{eq12a}, 
       \begin{align}  \nonumber 
		\phi_j^{n+1} &  
	 \leq  \left(1-\lambda \dfrac{\alpha_{j+1/2}+\alpha_{j-1/2}}{2}\right) \phi_{j}+ \dfrac{\lambda \alpha_{j+1/2} }{2} \phi_{j+1}+\dfrac{\lambda \alpha_{j-1/2} }{2} \phi_{\max}.  
\end{align}
The right-hand side is a convex combination of $\phi_j$, $ \phi_{j+1} $, and $\phi_{\max}$, so if 
 $\phi_j \in [0, \phi_{\max}]$ and $\phi_{j+1} \in [0, \phi_{\max}]$, then $\smash{\phi_j^{n+1} \leq \phi_{\max}}$.   This concludes the 
 proof of Theorem~\ref{thm:lfinv}. 
 }
\end{proof}

\begin{remark} \label{llfrem}  {Clearly,} the  previous proof also handles the case of periodic boundary conditions \eqref{bcperiodicdiscr}  if we assume that the 
 marching  formula \eqref{llf-march}  verbatim for all $j=0, \dots, M-1$ but understand all indices~$j$ ``modulo~$M$'', with the effect that 
   \eqref{bcperiodicdiscr}  is indeed enforced. As for the zero-flux boundary conditions  \eqref{bczerofluxdiscr}, consider for example the case 
    $\smash{\boldsymbol{\mathcal{F}}^n_{-1/2} = \boldsymbol{0}}$. The corresponding boundary scheme then becomes 
    \begin{align}  \label{bcscheme}
    \Phi_0^{n+1} = \Phi_0^n - \lambda \boldsymbol{\mathcal{F}}^{\mathrm{LLF}}_{1/2}= \Phi_0^n - \frac{\lambda}{2} \bigl( 
     \boldsymbol{f} ( \Phi_0^n) +   \boldsymbol{f} ( \Phi_1^n)  - \alpha_{1/2}^n  ( \Phi_1^n - \Phi_0^n)  \bigr). 
    \end{align} 
   Simpler versions of the arguments  used in the proof of Theorem~\ref{thm:lfinv} 
   {now}   suffice to show that also the boundary scheme \eqref{bcscheme} satisfies the IRP property. For instance, 
     one may formally set $\boldsymbol{\mathcal{G}}_2 (\Phi_{-1}) \coloneqq  \boldsymbol{0}$  
      and $\alpha_{-1/2} \coloneqq 0$ to demonstrate that $\smash{\phi_{i,0}^{n+1} \geq 0}$  for all $i = 1, \dots, N$, 
       and $\phi_{-1} \coloneqq \phi_{\max} $ and $\alpha_{-1/2} =0$ in \eqref{eq12a}  to deduce that $\phi_0^{n+1} \leq \phi_{\max}$ 
        (under  the conditions stated in Theorem~\ref{thm:lfinv}). The boundary condition $\smash{\boldsymbol{\mathcal{F}}_{M-1/2}^n = \boldsymbol{0}}$ is treated by similar arguments. 
        \end{remark} 

\subsection{IRP property of the HLL scheme} \label{subsec3.4} 

Let us now consider the HLL scheme   \cite{harten1983upstream} defined by  the numerical flux 
	\begin{equation}\label{eq:hll}
		\boldsymbol{\mathcal{F}}^{\mathrm{HLL}}_{j+1/2}\coloneqq \dfrac{S_{\mathrm{R},j+1/2}^{+}\boldsymbol{f}(\Phi_{j})-S_{\mathrm{L},j+1/2}^{-}\boldsymbol{f}(\Phi_{j+1})+S_{\mathrm{L},j+1/2}^{-}S_{\mathrm{R},j+1/2}^{+}(\Phi_{j+1}-\Phi_j)}{S_{\mathrm{R},j+1/2}^{+}-S_{\mathrm{L},j+1/2}^{-}},
	\end{equation}
	where $S_{\mathrm{L},j+1/2}$ and $S_{\mathrm{R},j+1/2}$ are defined as in \eqref{eq:disc-bounds}, and we have used the notation $a^{-}=\min(a,0)$ and $a^{+}=\max(a,0)$. First of all, observe that the HLL numerical flux \eqref{eq:hll} is well defined for both the MCLWR  and MLB  models, 
	  since \eqref{M1M2ineq} ensures that always 
	  \begin{align} \label{walways} 
	  S_{\mathrm{L},j+1/2}^{-} < S_{\mathrm{R},j+1/2}^{+} \quad \text{for all~$j$}.
	  \end{align}  
For the proof of the invariant region principle of the HLL scheme, we also use the bound for the smallest eigenvalue stipulated by
  \eqref{eq:kappafun} and \eqref{sl34}.

\begin{theorem} \label{th:hll} Consider the HLL scheme defined by the marching formula 
\begin{align}  \label{hll-march} 
  \Phi_j^{n+1} = \Phi_j^n - \lambda_n\bigl(    \boldsymbol{\mathcal{F}}^{\mathrm{HLL}}_{j+1/2} - \boldsymbol{\mathcal{F}}^{\mathrm{HLL}}_{j-1/2}
   \bigr), 
   \end{align} 
   where the HLL numerical flux is given by \eqref{eq:hll} along with 
   the definitions \eqref{eq:disc-bounds} and \eqref{sl34} of 
	$S_{\mathrm{R},j+1/2}$ and $S_{\mathrm{L},j+1/2}$, respectively.  
    If the CFL condition 
    \begin{align} \label{hllcfl}  
       \alpha \lambda_n \leq \frac{1}{2} , \quad \text{\em where} \quad \alpha \coloneqq \max_j \bigl\{  |S_{\mathrm{L},j+1/2}|,  |S_{\mathrm{R},j+1/2}|  \bigr\} 
     \end{align} is in effect, then the HLL scheme satisfies the invariant region preservation property \eqref{eq3.7} for the multispecies kinematic flow model \eqref{eq:main-eq}  with the velocity functions~$v_i(\Phi)$   {chosen such that there exists a function $w=w(\phi)$ satisfying  \eqref{w-phi-ass}--\eqref{psiass}}.  
 \end{theorem} 

\begin{proof} For simplicity,  let us keep the same notation as in the proof of Theorem \ref{thm:lfinv}, i.e., we set 
	$\smash{\Phi_j\coloneqq \Phi_j^n}$ and $\lambda \coloneqq  \lambda_n$. 
A straightforward computation
 reveals that the marching formula of the HLL scheme \eqref{hll-march}  
  can be written as follows, where we assume that $\smash{S_{\mathrm{R}, j+1/2}^+ >0}$, $\smash{S_{\mathrm{L}, j+1/2}^-<0}$,
   $\smash{S_{\mathrm{R}, j-1/2}^+ >0}$, and $\smash{S_{\mathrm{L}, j-1/2}^-<0}$:  
\begin{align*}   
\Phi_j^{n+1} &=  \Phi_j - \lambda \left( 
\dfrac{S_{\mathrm{R},j+1/2}^{+}\boldsymbol{f}(\Phi_{j})-S_{\mathrm{L},j+1/2}^{-}\boldsymbol{f}(\Phi_{j+1})+S_{\mathrm{L},j+1/2}^{-}S_{\mathrm{R},j+1/2}^{+}(\Phi_{j+1}-\Phi_j)}{S_{\mathrm{R},j+1/2}^{+}-S_{\mathrm{L},j+1/2}^{-}} \right.  \\
& \quad \qquad \qquad - \left. \dfrac{S_{\mathrm{R},j-1/2}^{+}\boldsymbol{f}(\Phi_{j-1})-S_{\mathrm{L},j-1/2}^{-}\boldsymbol{f}(\Phi_{j})+S_{\mathrm{L},j-1/2}^{-}S_{\mathrm{R},j-1/2}^{+}(\Phi_{j}-\Phi_{j-1})}{S_{\mathrm{R},j-1/2}^{+}-S_{\mathrm{L},j-1/2}^{-}}  \right)  \\ 
& = \left(1+\dfrac{2\lambda S_{\mathrm{L},j+1/2}^{-} S_{\mathrm{R},j+1/2}^{+}}{S_{\mathrm{R},j+1/2}^{+}-S_{\mathrm{L},j+1/2}^{-}}+\dfrac{2\lambda S_{\mathrm{L},j-1/2}^{-} S_{\mathrm{R},j-1/2}^{+}}{S_{\mathrm{R},j-1/2}^{+}-S_{\mathrm{L},j-1/2}^{-}}\right) \Phi_j + \dfrac{-\lambda  S_{\mathrm{L},j+1/2}^{-} S_{\mathrm{R},j+1/2}^{+}}{S_{\mathrm{R},j+1/2}^{+}-S_{\mathrm{L},j+1/2}^{-}} 
 \left(\Phi_{j+1}-\dfrac{1}{S_{\mathrm{R},j+1/2}^{+}} \boldsymbol{f}(\Phi_{j+1})\right)\\
&\quad  + \dfrac{-\lambda  S_{\mathrm{L},j+1/2}^{-} S_{\mathrm{R},j+1/2}^{+}}{S_{\mathrm{R},j+1/2}^{+}-S_{\mathrm{L},j+1/2}^{-}} 
 \left(\Phi_{j}+\dfrac{1}{(-S_{\mathrm{L},j+1/2}^{-})} \boldsymbol{f} (\Phi_{j})\right) 
   + \dfrac{-\lambda  S_{\mathrm{L},j-1/2}^{-} S_{\mathrm{R},j-1/2}^{+}}{S_{\mathrm{R},j-1/2}^{+}-S_{\mathrm{L},j-1/2}^{-}}
    \left(\Phi_{j-1}+\dfrac{1}{(-S_{\mathrm{L},j-1/2}^{-})} \boldsymbol{f} (\Phi_{j-1})\right)\\
&\quad  + \dfrac{-\lambda  S_{\mathrm{L},j-1/2}^{-} S_{\mathrm{R},j-1/2}^{+}}{S_{\mathrm{R},j-1/2}^{+}-S_{\mathrm{L},j-1/2}^{-}} \left(\Phi_{j}-\dfrac{1}{S_{\mathrm{R},j-1/2}^{+}}\boldsymbol{f} (\Phi_{j})\right).
  \end{align*} 
Consequently, if we define the coefficients 
\begin{align} \label{gammajphdef} 
	\gamma_{j+1/2}\coloneqq\dfrac{-S_{\mathrm{L},j+1/2}^{-} S_{\mathrm{R},j+1/2}^{+}}{S_{\mathrm{R},j+1/2}^{+}-S_{\mathrm{L},j+1/2}^{-}}>0  \quad \text{for all~$j$} 
\end{align}  
and  the functions
\begin{align*} 
		\mathcal{G}_1(\Phi_{j+1})& \coloneqq \Phi_{j+1}-\dfrac{1}{S_{\mathrm{R},j+1/2}^{+}} \boldsymbol{f} (\Phi_{j+1}),\quad &\mathcal{G}_2(\Phi_{j}) \coloneqq \Phi_j+\dfrac{1}{(-S_{\mathrm{L},j+1/2}^{-})} \boldsymbol{f}(\Phi_j),\\
		\mathcal{G}_3(\Phi_{j-1})& \coloneqq \Phi_{j-1}+\dfrac{1}{(-S_{\mathrm{L},j-1/2}^{-})} \boldsymbol{f}(\Phi_{j-1})	,\quad &\mathcal{G}_4(\Phi_{j}) \coloneqq \Phi_{j-1}-\dfrac{1}{S_{\mathrm{R},j-1/2}^{+}} \boldsymbol{f}(\Phi_j), 
	\end{align*}
 then we obtain 
\begin{equation}\label{eq:convex-comb}
	\phi_{i,j}^{n+1} = \bigl(1-2(\gamma_{j+1/2}+\gamma_{j-1/2}) \bigr) \phi_{i,j}+\gamma_{j+1/2} \mathcal{G}_{1,i}(\Phi_{j+1})+\gamma_{j+1/2}\mathcal{G}_{2,i}(\Phi_j)+\gamma_{j-1/2}\mathcal{G}_{3,i}(\Phi_{j-1})+\gamma_{j-1/2}\mathcal{G}_{4,i}(\Phi_j). 
\end{equation}
  Since $\gamma_{j+1/2} > 0$ and  $\gamma_{j-1/2} > 0$, under the CFL condition \eqref{hllcfl} this identity represents a convex combination of~$\phi_{i,j}$,   $\mathcal{G}_{1,i}  ( \Phi_{j+1})$, 
    $\mathcal{G}_{2,i}  ( \Phi_{j})$,  
  $\mathcal{G}_{3,i}  ( \Phi_{j-1})$, and~$\mathcal{G}_{4,i}  ( \Phi_{j})$. If  we assume that $\phi_{i,j} \geq 0$ for all~$i=1, \dots, N$
    and $j \in \mathbb{Z}_M$ and take into account  that $f_i ( \Phi) = \phi_i v_i (\Phi)$,   then we get 
 \begin{align*}
 	\mathcal{G}_{1,i} (\Phi_{j+1}) & = \phi_{i,j+1}\left(1-\dfrac{v_i(\Phi_{j+1})}{S_{\mathrm{R},j+1/2}^{+}}\right)\geq \phi_{i,j+1}\left(1-\dfrac{M_2(\Phi_{j+1})}{S_{\mathrm{R},j+1/2}^{+}}\right)\geq 0,\\
 	\mathcal{G}_{2,i} (\Phi_{j}) &=\phi_{i,j}\left(1+\dfrac{v_i(\Phi_{j})}{(-S_{\mathrm{L},j+1/2}^{-})}\right)\geq \phi_{i,j}\left(1+\dfrac{M_1(\Phi_{j})}{(-S_{\mathrm{L},j+1/2}^{-})}\right)\geq 0,\\
 	\mathcal{G}_{3,i} (\Phi_{j-1}) &=\phi_{i,j-1}\left(1+\dfrac{v_i(\Phi_{j-1})}{(-S_{\mathrm{L},j-1/2}^{-})}\right)\geq \phi_{i,j}\left(1+\dfrac{M_1(\Phi_{j-1})}{(-S_{\mathrm{L},j-1/2}^{-})}\right)\geq 0,\\
 	\mathcal{G}_{4,i} (\Phi_{j}) & = \phi_{i,j}\left(1-\dfrac{v_i(\Phi_{j})}{S_{\mathrm{R},j-1/2}^{+}}\right)\geq \phi_{i,j}\left(1-\dfrac{M_2(\Phi_{j})}{S_{\mathrm{R},j-1/2}^{+}}\right)\geq 0.
 \end{align*}
 Thus, $\mathcal{G}_{k,i} (\Phi_j) \geq 0$ for all~$i=1, \dots, N$,   $k=1, \dots,4$,  and $j \in \mathbb{Z}_M$.
 Since the coefficient of $\phi_{i,j}$ is non-negative (this follows from the CFL condition),  we deduce that if $\phi_{i,j}\geq 0$, $\phi_{i,j+1}\geq 0$, and $\phi_{i,j-1}\geq 0$, 
   then $\phi_{i,j}^{n+1}  \geq 0$.  
   
   It remains to prove that if $\phi_j  \leq \phi_{\max}$ for all~$j$, then 
 $\smash{ \phi_j^{n+1} \leq \phi_{\max}}$ for all~$j$. To this end, we rewrite \eqref{eq:convex-comb} as
 	\begin{align}\label{eq144} \begin{split} 
 	\phi_{i,j}^{n+1} 
 	& =    \Biggl( 1 - \lambda  \dfrac{S_{\mathrm{R},j+1/2}^{+} }{S_{\mathrm{R},j+1/2}^{+}-S_{\mathrm{L},j+1/2}^{-}}
 	\bigl( v_i ( \Phi_j) - S_{\mathrm{L},j+1/2}^- \bigr)  
 	+  \dfrac{ \lambda  S_{\mathrm{L},j-1/2}^{-}}{S_{\mathrm{R},j-1/2}^{+}-S_{\mathrm{L},j-1/2}^{-}} 
 	\bigl( S_{\mathrm{R},j-1/2}^+ - v_i ( \Phi_j) \bigr) \Biggr) \, \phi_{i,j}  \\ & \quad 
 	+ \dfrac{  - \lambda  S_{\mathrm{L},j+1/2}^{-}}{S_{\mathrm{R},j+1/2}^{+}-S_{\mathrm{L},j+1/2}^{-}}
 	\bigl( S_{\mathrm{R}, j+1/2}^+- v_i ( \Phi_{j+1}) \bigr)  \phi_{i,j+1} 
 	+\dfrac{ \lambda S_{\mathrm{R},j-1/2}^{+} }{S_{\mathrm{R},j-1/2}^{+}-S_{\mathrm{L},j-1/2}^{-}}
 	\bigl( v_i( \Phi_{j-1} ) - S_{\mathrm{L},j-1/2}^-\bigr) \phi_{i,j-1}. 
 	\end{split} 
 \end{align} 
  Summing \eqref{eq144} over $i=1, \dots, N$, we get 
   \begin{align} \label{eq145}  \begin{split} 
 \phi_{j}^{n+1} 
 & =  
\sum_{i=1}^N   \Biggl( 1 - \lambda  \dfrac{S_{\mathrm{R},j+1/2}^{+} }{S_{\mathrm{R},j+1/2}^{+}-S_{\mathrm{L},j+1/2}^{-}}
  \bigl( v_i ( \Phi_j) - S_{\mathrm{L},j+1/2}^- \bigr)  
   +  \dfrac{ \lambda  S_{\mathrm{L},j-1/2}^{-}}{S_{\mathrm{R},j-1/2}^{+}-S_{\mathrm{L},j-1/2}^{-}} 
    \bigl( S_{\mathrm{R},j-1/2}^+ - v_i ( \Phi_j) \bigr) \Biggr) \, \phi_{i,j}  \\ & \quad 
        +\sum_{i=1}^N   \dfrac{  - \lambda  S_{\mathrm{L},j+1/2}^{-}}{S_{\mathrm{R},j+1/2}^{+}-S_{\mathrm{L},j+1/2}^{-}}
     \bigl( S_{\mathrm{R}, j+1/2}^+- v_i ( \Phi_{j+1}) \bigr)  \phi_{i,j+1} 
      + \sum_{i=1}^N  \dfrac{ \lambda S_{\mathrm{R},j-1/2}^{+} }{S_{\mathrm{R},j-1/2}^{+}-S_{\mathrm{L},j-1/2}^{-}}
       \bigl( v_i( \Phi_{j-1} ) - S_{\mathrm{L},j-1/2}^-\bigr) \phi_{i,j-1}.
       \end{split} 
       \end{align} 
   From \eqref{eq145}, we get 
\begin{align} \label{eq146}  \begin{split} 
\phi_j^{n+1} & = ( 1- \lambda \gamma_{j+1/2}-  \lambda \gamma_{j-1/2}   ) \phi_j + \lambda \gamma_{j+1/2} \phi_{j+1} + 
 \lambda \gamma_{j-1/2} \phi_{j-1}  \\ & \quad
-   \frac{\lambda S_{\mathrm{L,j-1/2}}^-}{  S_{\mathrm{R,j-1/2}}^+ -  S_{\mathrm{L,j-1/2}}^-} 
 \sum_{i=1}^N v_i ( \Phi_{j}) \phi_{i,j} 
  -  \frac{\lambda S_{\mathrm{R,j+1/2}}^+}{  S_{\mathrm{R,j+1/2}}^+ -  S_{\mathrm{L,j+1/2}}^-} 
 \sum_{i=1}^N v_i ( \Phi_j) \phi_{i,j} \\ & \quad +   \frac{\lambda S_{\mathrm{L,j+1/2}}^-}{  S_{\mathrm{R,j+1/2}}^+ -  S_{\mathrm{L,j+1/2}}^-} 
 \sum_{i=1}^N v_i ( \Phi_{j+1}) \phi_{i,j+1}  +   \frac{\lambda S_{\mathrm{R,j-1/2}}^+}{  S_{\mathrm{R,j-1/2}}^+ -  S_{\mathrm{L,j-1/2}}^-} 
 \sum_{i=1}^N v_i ( \Phi_{j-1}) \phi_{i,{j-1}}. 
\end{split} 
\end{align} 
 {Utilizing \eqref{w-phi-ass} we may rewrite this identity as 
\begin{align} \label{eq146}  \begin{split} 
\phi_j^{n+1} & = ( 1- \lambda \gamma_{j+1/2}-  \lambda \gamma_{j-1/2}   ) \phi_j + \lambda \gamma_{j+1/2} \phi_{j+1} + 
 \lambda \gamma_{j-1/2} \phi_{j-1}  \\ & \quad
-   \frac{\lambda S_{\mathrm{L,j-1/2}}^- w( \phi_j) \boldsymbol{\kappa}^{\mathrm{T}} \Phi_j }{  S_{\mathrm{R,j-1/2}}^+ -  S_{\mathrm{L,j-1/2}}^-}  
  -  \frac{\lambda S_{\mathrm{R,j+1/2}}^+ w( \phi_j) \boldsymbol{\kappa}^{\mathrm{T}} \Phi_j}{  S_{\mathrm{R,j+1/2}}^+ -  S_{\mathrm{L,j+1/2}}^-} 
 +   \frac{\lambda S_{\mathrm{L,j+1/2}}^- w( \phi_{j+1}) \boldsymbol{\kappa}^{\mathrm{T}} \Phi_{j+1} }{  S_{\mathrm{R,j+1/2}}^+ -  S_{\mathrm{L,j+1/2}}^-} 
  +   \frac{\lambda S_{\mathrm{R,j-1/2}}^+ w( \phi_{j-1}) \boldsymbol{\kappa}^{\mathrm{T}} \Phi_{j-1} }{  S_{\mathrm{R,j-1/2}}^+ -  S_{\mathrm{L,j-1/2}}^-} . 
\end{split} 
\end{align} 
Since by \eqref{wpos}, the second and third term in the second line of the right-hand side of \eqref{eq146} are nonpositive, we get 
\begin{align} \label{eq146a}  \begin{split} 
\phi_j^{n+1} & \leq  ( 1- \lambda \gamma_{j+1/2}-  \lambda \gamma_{j-1/2}   ) \phi_j + \lambda \gamma_{j+1/2} \phi_{j+1} + 
 \lambda \gamma_{j-1/2} \phi_{j-1} 
-   \frac{\lambda S_{\mathrm{L,j-1/2}}^- w( \phi_j) \boldsymbol{\kappa}^{\mathrm{T}} \Phi_j }{  S_{\mathrm{R,j-1/2}}^+ -  S_{\mathrm{L,j-1/2}}^-}  
  +   \frac{\lambda S_{\mathrm{R,j-1/2}}^+ w( \phi_{j-1}) \boldsymbol{\kappa}^{\mathrm{T}} \Phi_{j-1} }{  S_{\mathrm{R,j-1/2}}^+ -  S_{\mathrm{L,j-1/2}}^-} . 
\end{split} 
\end{align} 
 If we assume that $\phi_j = \phi_{\max}$, $\phi_{j-1} = \phi_{\max}$, and $\phi_{j+1} \in [0, \phi_{\max}]$, then we   deduce from 
   \eqref{wphimax}, by  a convex combination argument, that $\phi_j^{n+1} \leq \phi_{\max}$. To handle the remaining cases, we rewrite \eqref{eq146a} as 
 \begin{align*} 
\phi_j^{n+1} & \leq    \biggl( \frac12 - \lambda \gamma_{j+1/2}\biggr) \phi_j + \biggl( \frac12 -  \lambda \gamma_{j-1/2}   \biggr) \phi_{\max}  + \lambda \gamma_{j+1/2} \phi_{j+1} + 
 \lambda \gamma_{j-1/2} \phi_{\max}   
 + \mathcal{Y}_1 + \mathcal{Y}_2 , 
 \end{align*} 
 where we define 
 \begin{align} \label{calY1defnew} 
 \mathcal{Y}_1  & \coloneqq \lambda \gamma_{j-1/2} ( \phi_{j-1} - \phi_{\max} ) 
 +    \frac{\lambda S_{\mathrm{R,j-1/2}}^+ w(\phi_{j-1}) \boldsymbol{\kappa}^{\mathrm{T}} \Phi_{j-1} }{  S_{\mathrm{R,j-1/2}}^+ -  
 S_{\mathrm{L,j-1/2}}^-}  , \\
    \label{calY2defnew} 
 \mathcal{Y}_2 & \coloneqq   \biggl( \frac12 - \lambda \gamma_{j-1/2} \biggr) ( \phi_{j} - \phi_{\max} )  
    -  \frac{\lambda S_{\mathrm{L,j-1/2}}^- w(\phi_j) \boldsymbol{\kappa}^{\mathrm{T}} \Phi_j }{  S_{\mathrm{R,j-1/2}}^+ -  S_{\mathrm{L,j-1/2}}^-}.  
 \end{align}
  If $\phi_{j-1} = \phi_{\max}$, then $\mathcal{Y}_1 =0$; otherwise we may divide \eqref{calY1defnew} by $ \phi_{j-1} - \phi_{\max}$ to obtain 
 \begin{align}  \label{huzzutnew} 
 \frac{\mathcal{Y}_1}{  \phi_{j-1} - \phi_{\max} } & = \lambda \gamma_{j-1/2} + \lambda 
   \frac{ S_{\mathrm{R,j-1/2}}^+ }{  S_{\mathrm{R,j-1/2}}^+ -  S_{\mathrm{L,j-1/2}}^-} \frac{w(\phi_{j-1})}{\phi_{j-1}- \phi_{\max}}  \boldsymbol{\kappa}^{\mathrm{T}} \Phi_{j-1} . 
 \end{align}  
 Notice that there exists a number $\xi_{j-1} \in [ \phi_{j-1}, \phi_{\max}]$ such that 
 \begin{align*} 
   \frac{w(\phi_{j-1})}{\phi_{j-1}- \phi_{\max}} =  \frac{w(\phi_{j-1}) - w( \phi_{\max}) }{\phi_{j-1}- \phi_{\max}}
   = w' ( \xi_{j-1})  \geq w' ( \phi_{j-1}), 
 \end{align*} 
 where we have used \eqref{wconvex}  to establish the last inequality.  
 Thus, from \eqref{huzzutnew} we get
 \begin{align}  \label{kjjy-new} 
  \frac{\mathcal{Y}_1}{  \phi_{j-1} - \phi_{\max} } \geq  \lambda 
   \frac{ S_{\mathrm{R,j-1/2}}^+ }{  S_{\mathrm{R,j-1/2}}^+ -  S_{\mathrm{L,j-1/2}}^-} \bigl( - S_{\mathrm{L},j-1/2}^- +  w'( \phi_{j-1} )  \boldsymbol{\kappa}^{\mathrm{T}} \Phi_{j-1} \bigr). 
 \end{align} 
 Arguing exactly as in   the discussion of inequality \eqref{kjjy-llf} in the proof of Theorem~\ref{thm:lfinv}, 
  we now deduce from \eqref{kjjy-new}  that   
  $\mathcal{Y}_1/(  \phi_{j-1} - \phi_{\max} )\geq  0$,
         and therefore 
          $\mathcal{Y}_1 \leq  0$.  Furthermore, if $\phi_j = \phi_{\max}$, then $\mathcal{Y}_2 =0$, otherwise we may divide \eqref{calY2defnew} 
by $\phi_j- \phi_{\max}$ to obtain 
\begin{align*} 
\frac{\mathcal{Y}_2}{\phi_j- \phi_{\max}} = \frac12 - \lambda \gamma_{j-1/2} - \lambda \frac{S_{\mathrm{L,j-1/2}}^-}{  S_{\mathrm{R,j-1/2}}^+ -  S_{\mathrm{L,j-1/2}}^-}
 \frac{w( \phi_j)}{\phi_j - \phi_{\max}} \boldsymbol{\kappa}^{\mathrm{T}} \Phi_j.  
\end{align*} 
In light of \eqref{hllcfl},  noticing that there exists $\xi_j \in [\phi_j, \phi_{\max}]$ such that 
\begin{align*} 
 \frac{w( \phi_j)}{\phi_j - \phi_{\max}}  =  \frac{w( \phi_j) - w ( \phi_{\max})} {\phi_j - \phi_{\max}} = w' ( \xi_j)\geq w'(\phi_j) 
\end{align*} 
 (where we once again utilize \eqref{wconvex}),   taking into account that  analogously to 
  \eqref{SLineq1},  
  \begin{align}  \label{SLineq2} 
  S_{\mathrm{L}, j-1/2}^- \leq  S_{\mathrm{L}, j-1/2} = \min_{ 0 \leq s \leq 1} 
    \psi \bigl( \phi_{j-1/2} (s) \bigr) \boldsymbol{\kappa}^{\mathrm{T}} \Phi_{j-1/2} (s) 
     \leq \psi ( \phi_{j}) \boldsymbol{\kappa}^{\mathrm{T}} \Phi_{j} \leq w'( \phi_{j} ) 
      \boldsymbol{\kappa}^{\mathrm{T}} \Phi_{j}, 
 \end{align}   
 hence $-S_{\mathrm{L}, j-1/2}^-  \geq   - w' ( \phi_{j} )  \boldsymbol{\kappa}^{\mathrm{T}} \Phi_{j}$, and  appealing to the CFL condition, 
 we obtain 
 \begin{align*} 
\frac{\mathcal{Y}_2}{\phi_j- \phi_{\max}} =  \frac12  - \lambda \frac{-S_{\mathrm{L,j-1/2}}^-}{  S_{\mathrm{R,j-1/2}}^+ -  S_{\mathrm{L,j-1/2}}^-} \left(S_{\mathrm{R},j-1/2}^+-w'(\xi_j)  \boldsymbol{\kappa}^{\mathrm{T}} \Phi_j \right)\geq \frac12+ \lambda S_{\mathrm{L,j-1/2}}^-\geq 0,
\end{align*} 
and therefore 
$\mathcal{Y}_2 \leq 0$. Consequently, we deduce that $\mathcal{Y}_1 + \mathcal{Y}_2 \leq 0$ 
which in the present case means that 
\begin{align*} 
\phi_j^{n+1} & \leq  \biggl( \frac12 - \lambda \gamma_{j+1/2}\biggr) \phi_j + \biggl( \frac12 -  \lambda \gamma_{j-1/2}   \biggr) \phi_{\max}  + \lambda \gamma_{j+1/2} \phi_{j+1} + 
 \lambda \gamma_{j-1/2} \phi_{\max}. 
 \end{align*} 
The right-hand side is a convex combination of $\phi_j$, $ \phi_{j+1} $, and $\phi_{\max}$, so  if 
 $\phi_j \in [0, \phi_{\max}]$ and $\phi_{j+1} \in [0, \phi_{\max}]$, then $\smash{\phi_j^{n+1} \leq \phi_{\max}}$.  Thus, we  always get  $\smash{\phi_j^{n+1} \leq \phi_{\max}}$. This concludes the 
    proof of Theorem~\ref{th:hll}.  
} 
      \end{proof} 
    
    \begin{remark} \label{hllrem} A remark similar to Remark~\ref{llfrem} is in place here. First of all, note that a simpler proof
     of $\smash{\phi_{i,j}^{n+1} \geq 0}$ applies if one or two of the bounds  
    $\smash{S_{\mathrm{R}, j+1/2}^+}$, $\smash{S_{\mathrm{L}, j+1/2}^-}$,
   $\smash{S_{\mathrm{R}, j-1/2}^+}$, and $\smash{S_{\mathrm{L}, j-1/2}^-}$
    are zero (we recall that \eqref{walways} is always in effect).  In the latter case we have, 
     of course, $\gamma_{j-1/2} = 0$ or  $\gamma_{j+1/2} = 0$ (cf.\ \eqref{gammajphdef}). 
     On the other hand, as for  the LLF scheme, it suffices to take all $j$-indices ``modulo $M$'' 
      to demonstrate that the proof of of Theorem~\ref{th:hll}  also handles periodic boundary conditions. In the case 
       of zero-flux boundary conditions  \eqref{bczerofluxdiscr}, the bound $\phi_{i,j}^{n+1} \geq 0$ for $j=0$ or $j=M-1$ follows  by taking into account that 
        certain terms in the discussion leading to \eqref{eq:convex-comb} are zero.  With respect to the upper bound
         of $\smash{\phi_0^{n+1}}$, we note that if $\smash{\boldsymbol{\mathcal{F}}_{-1/2}^n = \boldsymbol{0}}$, then $\gamma_{-1/2}=0$ and the last two terms  on the 
          right-hand side of \eqref{eq146a}  are zero, and we deduce that $\smash{\phi_0^{n+1} \leq \phi_{\max}}$ by applying a standard 
           convex combination argument to the remaining inequality. For $j=M-1$, and considering the boundary condition $\smash{\boldsymbol{\mathcal{F}}_{M-1/2}^n = \boldsymbol{0}}$, 
          \eqref{eq:convex-comb}  is valid for $j=M-1$ if we set $\gamma_{M-1/2} =0$. The remainder of the proof for both models remains valid. 
           Consequently,   Theorems~\ref{thm:lfinv} and~\ref{th:hll} remain valid if appropriate numerical fluxes in the respective 
            marching formulas  \eqref{llf-march} and \eqref{hll-march} are replaced by  boundary fluxes coming  
             from~\eqref{bcperiodicdiscr} or~\eqref{bczerofluxdiscr}. 
        \end{remark}

    \begin{remark}
    	Notice that the CFL conditions \eqref{llfcfl} and \eqref{hllcfl} are utilized to calculate $\Delta t$ 
	  adaptively in every time step. However, one may employ a more restrictive but fixed CFL condition 
	   defined by    bounds for $M_1(\Phi)$ and $M_2(\Phi)$. For instance, 
	    for the  MCLWR model one may employ the condition 
		$\lambda \max \{\beta_1\|v\|_{\infty},\beta_1 \phi_{\max}\|\psi\|_{\infty} \}\leq \mu$,
    	while for MLB model one can consider 
    		$\lambda \max  \{\delta_1\|w_1\|_{\infty},\delta_1 \phi_{\max}\|\psi\|_{\infty} \}\leq \mu$,
    	where $w_1(\phi) = C(1-\phi)V(\phi)$, $\psi$~is defined by~\eqref{eq:kappafun},  and 
    	\begin{equation}\label{eq:mu}
    		\mu = 
    		\begin{cases}
    			1& \text{\em  for the  LLF scheme,} \\
    			 1/2 &\text{\em  for  the HLL scheme.} 
    		\end{cases}
    	\end{equation}
    \end{remark}
    
\section{Invariant-region-preserving WENO reconstruction}\label{sec:invcweno}
\subsection{Weighted essentially non-oscillatory (WENO) reconstruction (scalar case)} \label{subsec4.1} 
The WENO reconstruction  \cite{levy1999central}  for the initial value problem of a 
 one-dimensional scalar conservation law  
	\begin{align*} 
		 \partial_t  u + \partial_x f(u) =0, \quad x \in \mathbb{R}, \quad t \in (0,T); \quad 
		u(x,0) =u_{0}(x), \quad x \in \mathbb{R}   
	\end{align*}
 is independent of the time variable, hence it is sufficient to  consider $u=u(x)$ as  a function of the spatial
  coordinate  only. We define the  cell averages~$\bar{u}_j$ of~$u$ over all  cells $I_j\coloneqq [x_{j-1/2},x_{j+1/2}]$: 
\begin{equation*}
	\bar{u}_j \coloneqq \dfrac{1}{\Delta x}\int_{x_{j-1/2}}^{x_{j+1/2}} u(x) \,  \mathrm{d}x, \quad j \in \mathbb{Z}. 
\end{equation*}
In what follows,   $\Pi_R$~denotes  the set of all  polynomials with real coefficients of maximal degree~$R \in \mathbb{N}_0$.

\begin{definition}
	Consider a set of data (point values or cell averages) and a polynomial $P_{\mathrm{opt}} \in \Pi_R$, which interpolates in  the appropriate sense all the given data (optimal polynomial). The 
	 {\em WENO operator}   computes a reconstruction polynomial
	 $P_{\mathrm{rec}}=\mathrm{WENO}(P_{\mathrm{opt}},P_1,\dots,P_m)\in \Pi_{R}$ 
	from $P_{\mathrm{opt}}\in \Pi_{R}$ and a set of $m$ lower-order alternative polynomials $P_1,\dots, P_m\in \Pi_{r}$, where $r<R$ and $m\geq 1$. The definition of $P_{\mathrm{rec}}$ depends on the choice of  a set of  real coefficients 
	\begin{align*} 
	C_0,C_1,\dots, C_m\in [0,1], \quad \text{\em where} \quad C_0 + C_1 + \cdots +  C_m =1, \quad C_0 >0, 
	\end{align*} 
	 the so-called  ``linear weights,''   as follows:
	\begin{enumerate}
		\item First, we define   $P_0 \in \Pi_R $ by  
		\begin{align*} 
			P_0(x) \coloneqq \dfrac{1}{C_0}\left(P_{\mathrm{opt}}(x)-\sum_{k=1}^{m}C_k P_k(x)\right).
		\end{align*}
		
		\item Then the nonlinear weights $\omega_0, \dots, \omega_m $ are computed from the linear ones as
		\begin{equation}\label{eq:nonl-weights}
			\alpha_k \coloneqq \dfrac{C_k}{(IS_k+\varepsilon)^p},\quad \omega_k \coloneqq \dfrac{\alpha_k}{ \alpha_0 + \cdots + \alpha_m }, \quad k=0, \dots, m,  
		\end{equation}
		where $\varepsilon>0$ is a small parameter,  $p\geq 2$, and  $IS_k$ denotes a suitable regularity indicator, 
		 e.g.\ the Jiang-Shu indicator {\em  \cite{jiang1996efficient}}  
		\begin{align*} 
			IS_k \coloneqq \sum_{\ell=1}^{\mathrm{degree}(P_k)} \Delta x
                        ^{2\ell-1} \int_{x_{j-1/2}}^{x_{j+1/2}}  \left(\dfrac{\mathrm{d}^{\ell} P_k(x)}{\mathrm{d} x^{\ell} }\right)^2 \, \mathrm{d}x, \quad k=0, \dots, m. 
		\end{align*}
		\item  Finally, the reconstruction polynomial $P_{\mathrm{rec}} \in \Pi_R$ is defined as 
		\begin{align}\label{eq:rec-poly}
			P_{\mathrm{rec}}(x) \coloneqq  \omega_0 P_0(x) + 	\omega_1 P_1(x) + \cdots + \omega_m P_m(x). 
			\end{align} 
			\end{enumerate}
\end{definition}

The  reconstructed values of~$u$ at the  boundaries~$x_{i-1/2}$ and $x_{j+1/2}$ of cell~$I_j$ are now given by 
\begin{align}\label{eq:recon}
	u_{j-1/2}^{\mathrm{R}} \coloneqq P_{\mathrm{rec}}(x_{j-1/2}), \quad 
		u_{j+1/2}^{\mathrm{L}} \coloneqq P_{\mathrm{rec}}(x_{j+1/2}).	\end{align}
Thus, based on cell averages $\bar{u}_j$ over all $I_j$,   a WENO reconstruction of order $2r+1$ ($r\in \{1,2\}$ within this work), 
 we consider on each interval a polynomial $P_{\mathrm{opt}} \in \Pi_{R=2r}$ along with  $m=r+1$ polynomials $P_1,\dots,P_{r+1} \in \Pi_r$ defined by 
	\begin{align*} 
		P_{\mathrm{opt}}(x) \coloneqq  \sum_{s=0}^{2r} a_{s}^{(0)}(x-x_j)^{s} \quad \text{and} \quad 
		P_{k}(x) \coloneqq \sum_{s=0}^{r} a_{s}^{(k)}(x-x_j)^{s}, \quad  k = 1,\ldots,r+1, 
	\end{align*}
where $P_{\mathrm{opt}}$ and $P_k$ interpolate  the cell averages associated with  
 the respective stencils 
 \begin{align*} 
 S_{0}=\bigcup_{\ell=-r}^{r} I_{j+\ell}, \quad \textrm{and}\quad S_{k}=\bigcup_{\ell=0}^{r} I_{j-r+k+\ell-1}, \quad  k = 1,\dots,r+1.
 \end{align*} 
 For a fixed cell $I_j$ these polynomials satisfy  
\begin{align*} 
	\dfrac{1}{\Delta x} \int_{I_{j+\ell}} P_{\mathrm{opt}}(x) \, \mathrm{d} x & =\bar{u}_{j+\ell},\quad  \ell  = -r,\dots,r, \\ 
 \label{eq:cell-avg-2}
	\dfrac{1}{\Delta x} \int_{I_{j-r+k+\ell-1}} P_k(x) \, \mathrm{d}x & =\bar{u}_{j-r+k+\ell-1},\quad   \ell = 0,\dots,r, \quad k=1, \dots, r+1. 
\end{align*}
The reconstructed  values \eqref{eq:recon} are easily computable since 
	\begin{align*} 
		P_0(x_{j\pm1/2})&=\sum_{s=0}^{2r} a_{s}^{(0)}\biggl(\pm\dfrac{\Delta x}{2}\biggr)^{s}=\sum_{s=0}^{2r} \alpha_{s}^{(0)} \bar{u}_{j+s-r} \quad \text{and} \\
		P_k(x_{j\pm1/2})&=\sum_{\ell=0}^{r} a_{\ell}^{(k)}\biggl(\pm\dfrac{\Delta x}{2}\biggr)^{\ell}=\sum_{\ell=0}^{r} \alpha_{\ell}^{(k)} \bar{u}_{j-r+k+\ell-1}, \quad 
		  k = 1,\dots,r+1. 
	\end{align*}

\subsection{WENO  reconstruction for systems of conservation laws } \label{subsec4.2} 

 Let us consider a first-order IRP numerical flux $\smash{\boldsymbol{\mathcal{F}}_{j+1/2}^{n}}$ for \eqref{eq:main-eq}, for 
  instance the LLF flux \eqref{eq:llf} or the HLL flux \eqref{eq:hll}. To describe WENO reconstructions for~\eqref{eq:main-eq}, we define $\smash{\boldsymbol{\mathcal{F}}(\Phi_{j}^n,\Phi_{j+1}^n)\coloneqq \boldsymbol{\mathcal{F}}_{j+1/2}^{n}}$, so we can write \eqref{eq:march} as
\begin{equation}\label{eq:marching-2}
		\Phi_j^{n+1}=\Phi_j^{n}-\lambda_n\big(\boldsymbol{\mathcal{F}}(\Phi_{j}^n,\Phi_{j+1}^n)-\boldsymbol{\mathcal{F}}(\Phi_{j-1}^n,\Phi_{j}^n)\big).
\end{equation}
In this work, we will use a component-wise WENO reconstruction, i.e., 
\begin{align*} 
& \Phi_{j+1/2}^{\mathrm{L}}\coloneqq \bigl(\phi_{1,j+1/2}^{\mathrm{L}},\dots,\phi_{N,j+1/2}^{\mathrm{L}}\bigr)^{\mathrm{T}}, \quad  
 \Phi_{j+1/2}^{\mathrm{R}}\coloneqq \bigl(\phi_{1,j+1/2}^{\mathrm{R}},\dots,\phi_{N,j+1/2}^{\mathrm{R}}\bigr)^{\mathrm{T}}, \\  
  & \text{where} \quad  
		\phi_{i,j-1/2}^{\mathrm{R}} \coloneqq  P^{(i)}_j(x_{j-1/2}),\quad 
		\phi_{i,j+1/2}^{\mathrm{L}} \coloneqq   P^{(i)}_j(x_{j+1/2}),  \quad i =1, \dots, N, 
	\end{align*}
and $\smash{P_j^{(i)}}$ denotes the reconstruction polynomial given by \eqref{eq:rec-poly}  for $i=1,\dots,N$. Then, we replace \eqref{eq:marching-2} by the  WENO marching formula 
\begin{equation}\label{eq:marching-cweno}
	\Phi_j^{n+1}=\Phi_j^{n}-\lambda_n \bigl( \boldsymbol{\mathcal{F}} (\Phi_{j+1/2}^{\mathrm{L}},\Phi_{j+1/2}^{\mathrm{R}} )-\boldsymbol{\mathcal{F}} (\Phi_{j-1/2}^{\mathrm{L}},\Phi_{j-1/2}^{\mathrm{R}} )\bigr),
\end{equation}
where 
\begin{align} \label{eq:recvalues} 
		\Phi_{j-1/2}^{\mathrm{R}} \coloneqq  \boldsymbol{P}_j(x_{j-1/2}), \quad 
		\Phi_{j+1/2}^{\mathrm{L}} \coloneqq \boldsymbol{P}_j(x_{j+1/2}),
	\end{align}
and we define the vector of polynomials  $\smash{ \boldsymbol{P}_j\coloneqq (P_j^{(1)},\dots,P_j^{(N)})^{\mathrm{T}}}$. In addition, to have the correct order of convergence of the method we need to use high-order time integration schemes;  we focus on this aspect later. 

\subsection{Zhang and Shu limiters for  multispecies kinematic flow models}\label{subsec:zs-limiters}
It is well known that the reconstruction procedure described in Section~\ref{sec:invcweno}  leads to  schemes that resolve discontinuities sharply but   in some cases fail to preserve the invariant region $\mathcal{D}_{\phi_{\mathrm{max}}}$ (see Examples~2 to~5 in Section~\ref{sec:numexa} or numerical examples in \cite{burger2011implementation,donat2010secular,burger2016polynomial}). This shortcoming motivated the development  of the linear scaling limiter 
 by Zhang and Shu \cite{zhang2010maximum}.  We herein   slightly modify this limiter to handle  equations of the form \eqref{eq:main-eq}.

    {The first step is to limit each concentration $\phi_i$. Let us assume that $\smash{\Phi_j^n\in \mathcal{D}_{\phi_{\mathrm{max}}}^{\delta}}$ for some small number $\delta>0$, where 
 	\begin{align*}
 		\mathcal{D}_{\phi_{\mathrm{max}}}^{\delta}\coloneqq \bigl\{\Phi\in \mathbb{R}^N: \phi_1\geq \delta, \dots, \phi_N\geq \delta, \, 
 		\phi\coloneqq \phi_1+\cdots+\phi_N \leq \phi_{\mathrm{max}} \bigr\}.
 	\end{align*}
 	For example, in all our numerical tests we set $\delta = 10^{-12}$. We then replace the polynomials $P_j^{(i)}(x)$ by
 	\begin{equation}\label{eq:ptilde}
 		\tilde{P}_j^{(i)}(x) \coloneqq \theta_i \bigl(P_j^{(i)}(x)-\phi_{i,j}^n \bigr)+\phi_{i,j}^n,\quad \theta_i
 		\coloneqq \min\Biggl\{\dfrac{\phi_{i,j}^n-\delta}{\phi_{i,j}^n-m_j^{(i)}},1\Biggr\},
 		\quad m_j^{(i)} \coloneqq \min_{x\in I_j}P_j^{(i)}(x), \quad i=1,\dots, N.   
 	\end{equation}
 	Then the cell average of $	\tilde{P}_j^{(i)}(x)$ over $I_j$ is still $\phi_{i,j}^n$ and  $\smash{\tilde{P}_{j}^{(i)}(x)\geq 0}$  for all $x\in I_j$ and $i=1,\dots,N$.} Next, we define  the polynomial 
\begin{equation}\label{eq:phat}
	\hat{P}_j(x) \coloneqq \hat{\theta} \left(\sum_{i=1}^{N} \tilde{P}_j^{(i)}(x)-\phi_j^n \right)+\phi_{j}^n, \quad	\hat{\theta}  \coloneqq 
	  \min\Biggl\{\Biggl|\dfrac{\phi_{\mathrm{max}}-\phi_{j}^n}{M_j-\phi_{j}^n}\Biggr|,1\Biggr\}, \quad 
	M_j \coloneqq \max_{x\in I_j}\left(\sum_{i=1}^{N }\tilde{P}_j^{(i)}(x)\right).
\end{equation}
Thus,  $\smash{\hat{P}_j(x)\leq \phi_{\mathrm{max}}}$  for all $x\in I_j$.
Finally,  we define the  modified polynomials 
\begin{equation}\label{eq:mod-poly}
	\bar{P}_j^{(i)}(x): = \hat{\theta} \bigl(\tilde{P}_j^{(i)}(x)-\phi_{i,j}^n \bigr)+\phi_{i,j}^n, \quad i=1, \dots, N, 
\end{equation}
and  replace  the reconstructed values \eqref{eq:recvalues} by 
\begin{equation}\label{eq:newrecvalues}
		\Phi_{j-1/2}^{\mathrm{R}}= \boldsymbol{\bar{P}}_{j}(x_{j-1/2}), \quad \Phi_{j+1/2}^{\mathrm{L}}=\boldsymbol{\bar{P}}_{j}(x_{j+1/2}).
\end{equation}

The quantities $\smash{m_j^{(i)}}$ in \eqref{eq:ptilde} and $M_j$ in \eqref{eq:phat} require evaluating the extrema of polynomials on each cell.  
 This inconvenience can be avoided if we define simplified limiters (as in  \cite{zhang2010maximum}) that reduce the 
  evaluations of each polynomial to a finite number of nodes of a 
 $G$-point Legendre-Gauss-Lobatto quadrature rule  on the interval $I_j = [x_{j-1/2},x_{j+1/2}]$. 
  This formula is exact  for the integral of polynomials of degree up to $2G-3$. We denote the 
quadrature points on~$I_j$ by 
\begin{equation}\label{eq:gauss-legendre}
	S_j \coloneqq  \bigl\{x_{j-1/2}=\hat{x}_j^{1},\hat{x}_j^{2},\dots,\hat{x}_j^{G-1},\hat{x}_j^{G} = x_{j+1/2} \bigr\}.
\end{equation}
Let $\hat{w}_{\alpha}$ be the quadrature weights for the interval $[-1/2, 1/2]$ such that $\hat{w}_{1} + \cdots + \hat{w}_{G} =1$. For instance, we have used a three-point ($G=3$) rule with weights  
\begin{equation*}
 \hat{w}_1=\frac{1}{6}, \quad  \hat{w}_2=\frac{2}{3}, \quad  \hat{w}_3 = \frac{1}{6} 
      \end{equation*}
      and a four-point ($G=4$) rule with weights
      \begin{equation*}
      	\hat{w}_1=\frac{1}{12}, \quad  \hat{w}_2=\frac{5}{12}, \quad  \hat{w}_3 = \frac{5}{12}, \quad 
      	\hat{w}_4 = \frac{1}{12} 
      \end{equation*}
      for the third and fifth order WENO reconstructions, respectively (notice that $\hat{w}_1=\hat{w}_G$). Then we can write 
\begin{equation}\label{eq:media-gl}
		{\Phi}_{j}^n = \dfrac{1}{\Delta x} \int_{I_j} \boldsymbol{\bar{P}}_j (x) \, \mathrm{d}x = \sum_{\alpha=1}^{G} \hat{w}_{\alpha} \boldsymbol{\bar{P}}_j 
		 \bigl(\hat{x}_j^{\alpha} \bigr)= \sum_{\alpha=2}^{G-1} \hat{w}_{\alpha} \boldsymbol{\bar{P}}_j \bigl(\hat{x}_j^{\alpha}\bigr)+\hat{w}_{1}	\Phi_{j-1/2}^{\mathrm{R}}+	\hat{w}_{G}	\Phi_{j+1/2}^{\mathrm{L}}.
	\end{equation}

We first prove  the following lemma related to  the polynomials $\boldsymbol{\bar{P}}_j\coloneqq (\bar{P}_j^{(1)},\dots,\bar{P}_j^{(N)})^{\mathrm{T}}$. 
\begin{lemma}\label{lemma:poly}
	Consider the reconstruction polynomials $\smash{\boldsymbol{\bar{P}}_j(x)}$ defined by \eqref{eq:mod-poly}. If   $\smash{{\Phi}_{j}^{n}\in \mathcal{D}_{\phi_{\mathrm{max}}}}$, then 
	$\smash{\boldsymbol{\bar{P}}_j(x)\in \mathcal{D}_{\phi_{\mathrm{max}}}}$ for all $x\in I_j$.  In particular, this is true for all $x\in S_j$, where $S_j$ is the stencil  \eqref{eq:gauss-legendre} 
	  of Legendre-Gauss-Lobatto quadrature points for~$I_j$. 
        \end{lemma}
       
\begin{proof}
	Let $x\in I_j$. By definition,  we know that 
	\begin{align} \label{eq4.22} 
		\bar{P}_j^{(i)}(x)=\hat{\theta} \bigl(\tilde{P}_j^{(i)}(x)-\phi_{i,j}^n \bigr)+\phi_{i,j}^n = \hat{\theta}\tilde{P}_j^{(i)}(x)+(1- \hat{\theta})\phi_{i,j}^n, 
		 \quad i=1, \dots, N. 
	\end{align}
	Since $0\leq \hat{\theta}\leq 1$ and $\tilde{P}_j^{(i)}(x)\geq 0$, there holds  $\smash{\bar{P}_j^{(i)}(x)\geq 0}$  for $i=1,\dots,N$. In addition, by \eqref{eq:phat}, 
	\begin{align} \label{eq4.23} 
		\sum_{i=1}^{N} \bar{P}_j^{(i)}(x) &= \hat{\theta} \left(\sum_{i=1}^{N} \tilde{P}_j^{(i)}(x)-\sum_{i=1}^{N} \phi_{i,j}^n\right)+\sum_{i=1}^{N} \phi_{i,j}^n 
		 = \hat{\theta} \left(\sum_{i=1}^{N} \tilde{P}_j^{(i)}(x) - {\phi}_{j}^n\right)+{\phi}_{j}^n =\hat{P}_j(x)  \leq \phi_{\mathrm{max}}.
	\end{align}
	 Combining \eqref{eq4.22} and \eqref{eq4.23} we deduce that 
	 $\smash{\boldsymbol{\bar{P}}_j(x)\in \mathcal{D}_{\phi_{\mathrm{max}}}}$  for all~$x\in I_j$.
\end{proof}
Now, we are in position to state the following result.
\begin{theorem}\label{teo:weno-inv}
	Consider the finite volume scheme \eqref{eq:marching-cweno}  associated with the  reconstruction 
	polynomials $\smash{\boldsymbol{\bar{P}}_j(x)}$ defined by \eqref{eq:mod-poly} in the sense that 
	\eqref{eq:newrecvalues}  is  used,  and where the  quantities $\smash{m_j^{(i)}}$ in \eqref{eq:ptilde} and $M_j$ in \eqref{eq:phat}
	 are redefined by  
	 \begin{align*} 
	  m_j^{(i)} \coloneqq \min_{x\in S_j}P_j^{(i)}(x), \quad i=1,\dots, N; \quad  M_j \coloneqq \max_{x\in S_j}\left(\sum_{i=1}^{N }\tilde{P}_j^{(i)}(x)\right), 
	 \end{align*} 
	 where $S_j$ is the stencil \eqref{eq:gauss-legendre} 
	  of Legendre-Gauss-Lobatto quadrature points for~$I_j$. Let $\mu$ be as in \eqref{eq:mu}. If the CFL condition
	\begin{equation}\label{eq:CFLMMP}
		\alpha \lambda_n  \leq \mu \hat{w}_1, \quad \alpha \coloneqq \max_j \bigl\{  |S_{\mathrm{L},j+1/2}|,  |S_{\mathrm{R},j+1/2}|  \bigr\}, 
	\end{equation}
	 is in effect and  $\smash{{\Phi}_{j}^{n}\in \mathcal{D}_{\phi_{\mathrm{max}}}}$, then $\smash{{\Phi}_{j}^{n+1}\in \mathcal{D}_{\phi_{\mathrm{max}}}}$. 
	  \end{theorem}
\begin{proof}
	We use the marching formula~\eqref{eq:marching-cweno} along with~\eqref{eq:media-gl} and  add and subtract 
	 $\smash{\boldsymbol{\mathcal{F}} (\Phi_{j-1/2}^{\mathrm{R}},\Phi_{j+1/2}^{\mathrm{L}})}$ to get
	\begin{align*}
		{\Phi}_{j}^{n+1}&=\sum_{\alpha=2}^{G-1} \hat{w}_{\alpha} \boldsymbol{\bar{P}}_j \bigl(\hat{x}_j^{\alpha}\bigr) +\hat{w}_G
		 \biggl(\Phi_{j+1/2}^{\mathrm{L}}- \dfrac{\lambda_n}{\hat{w}_G} \bigl( \boldsymbol{\mathcal{F}} (\Phi_{j+1/2}^{\mathrm{L}},\Phi_{j+1/2}^{\mathrm{R}} )-\boldsymbol{\mathcal{F}} (\Phi_{j-1/2}^{\mathrm{R}},\Phi_{j+1/2}^{\mathrm{L}} )\bigr) \biggr)\\
		&\quad +\hat{w}_1\biggl(\Phi_{j-1/2}^{\mathrm{R}}- \dfrac{\lambda_n}{\hat{w}_1} \bigl( \boldsymbol{\mathcal{F}} (\Phi_{j-1/2}^{\mathrm{R}},\Phi_{j+1/2}^{\mathrm{L}} )-\boldsymbol{\mathcal{F}} (\Phi_{j-1/2}^{\mathrm{L}},\Phi_{j-1/2}^{\mathrm{R}} )\bigr) \biggr)\\
		&=\sum_{\alpha=2}^{G-1} \hat{w}_{\alpha} \boldsymbol{\bar{P}}_j \bigl(\hat{x}_j^{\alpha} \bigr)+\hat{w}_G \boldsymbol{\mathcal{H}}
		 \bigl(\Phi_{j-1/2}^{\mathrm{R}},\Phi_{j+1/2}^{\mathrm{L}},\Phi_{j+1/2}^{\mathrm{R}}\bigr) +\hat{w}_1 \boldsymbol{\mathcal{H}} 
		  \bigl(\Phi_{j-1/2}^{\mathrm{L}},\Phi_{j-1/2}^{\mathrm{R}},\Phi_{j+1/2}^{\mathrm{L}}\bigr),
	\end{align*}
	where $\boldsymbol{\mathcal{H}}(\cdot,\cdot,\cdot)$ is the three-point operator analyzed in Theorems~\ref{thm:lfinv}-\ref{th:hll}. By Lemma~\ref{lemma:poly} we know that  $\smash{\boldsymbol{\bar{P}}_j (\hat{x}_j^{\alpha})\in \mathcal{D}_{\phi_{\mathrm{max}}}}$  for $j=2,\dots,G-1$, and Theorems~\ref{thm:lfinv}-\ref{th:hll}, in conjunction   with the  CFL condition~\eqref{eq:CFLMMP},  guarantees that   
	\begin{align*}
		\boldsymbol{\mathcal{H}} \bigl(\Phi_{j-1/2}^{\mathrm{R}},\Phi_{j+1/2}^{\mathrm{L}},\Phi_{j+1/2}^{\mathrm{R}}\bigr), \;  \boldsymbol{\mathcal{H}}
		 \bigl(\Phi_{j-1/2}^{\mathrm{L}},\Phi_{j-1/2}^{\mathrm{R}},\Phi_{j+1/2}^{\mathrm{L}} \bigr) \in \mathcal{D}_{\phi_{\mathrm{max}}}. 
	\end{align*}
	Consequently,   $\smash{\Phi_j^{n+1}}$ can be expressed as a convex combination of terms in $\mathcal{D}_{\phi_{\mathrm{max}}}$. This   concludes the proof.
\end{proof}

  {
\subsection{Accuracy of the limiter}\label{subsec:aol} We now adapt the ideas in \cite{zhang2010positivity,zhang2017positivity} to  rigorously discuss  the accuracy of the limiter \eqref{eq:mod-poly}. To this end, let us suppose that $\smash{\boldsymbol{P}_{j}(x)}$  approximates a smooth solution $\smash{\Phi(x)=(\phi_1(x),\dots,\phi_N(x))^{\mathrm{T}}\in \mathcal{D}_{\phi_{\mathrm{max}}}^{\delta}}$, i.e.,  $\smash{\|\boldsymbol{P}_j(x)-\Phi(x)\|=\mathcal{O}(\Delta x^{r+1})}$ for any $x\in I_j$. Our goal is to show that limiting polynomials $\smash{\boldsymbol{\bar{P}}_{j}(x)}$ given by \eqref{eq:mod-poly} also satisfies $\smash{\|\boldsymbol{\bar{P}}_j(x)-\Phi(x)\|=\mathcal{O}(\Delta x^{r+1})}$ for any $x\in I_j$. First, we analyze the accuracy of the limiter \eqref{eq:ptilde}. For a given~$i$, we only discuss the case when $\smash{\tilde{P}_j^{(i)}(x)}$ is not a constant polynomial. If $\theta_i = 1$, then $\smash{\tilde{P}_j^{(i)}(x)=P_j^{(i)}(x)}$ and we have nothing to show. So, let us assume that
 \begin{align} \label{thetaiass} 
 \theta_i=\frac{\phi_{i,j}^n-\delta}{\phi_{i,j}^n-m_j^{(i)}}. 
 \end{align} 
 We observe that
\begin{align*}
	\tilde{P}_j^{(i)}(x)-	P_j^{(i)}(x) &=\theta_i \bigl(P_j^{(i)}(x)-\phi_{i,j}^n \bigr)+\phi_{i,j}^n-P_j^{(i)}(x) = (\theta_i-1) \bigl(P_j^{(i)}(x)-\phi_{i,j}^n \bigr)  \\
	&=  \dfrac{m_j^{(i)}-\delta}{\phi_{i,j}^n-m_j^{(i)}}  \bigl(P_j^{(i)}(x)-\phi_{i,j}^n \bigr)= \bigl(\delta-m_j^{(i)} \bigr) \dfrac{\phi_{i,j}^n-P_j^{(i)}(x)}{\phi_{i,j}^n-m_j^{(i)}},
\end{align*}
so we get 
\begin{align*} 
\bigl|\tilde{P}_j^{(i)}(x)-	P_j^{(i)}(x)\bigr|\leq  \bigl|\delta-m_j^{(i)}\bigr|\left|\frac{\phi_{i,j}^n-P_j^{(i)}(x)}{\phi_{i,j}^n-m_j^{(i)}}\right|.
\end{align*} 
 The assumption that $\theta_i$ given by \eqref{thetaiass} satisfies $\theta_i <1$ implies that $\phi_{i,j}^n-\delta<\phi_{i,j}^n-m_j^{(i)}$, so we get $\smash{m_j^{(i)}<\delta}$. So, if $\smash{x_{*}^{(i)}\in I_j}$ is such that $\smash{\tilde{P}_j^{(i)}(x_{*}^{(i)}) = m_j^{(i)}}$, then   $\smash{\tilde{P}_j^{(i)}(x_{*}^{(i)}) = m_j^{(i)}<\delta\leq \phi_i(x_{*}^{(i)})}$. From this we get the bound 
  \begin{align*} 
  \bigl|\delta-m_j^{(i)} \bigr|\leq \bigl|\phi_i(x_{*}^{(i)})-\tilde{P}_j^{(i)}(x_{*}^{(i)})\bigr|\leq 
   \max_{x\in I_j} \bigl|\phi_i(x)-\tilde{P}_j^{(i)}(x) \bigr|. 
   \end{align*} 
    Consequently, 
   we get $ \delta-m_j^{(i)} = \mathcal{O}(\Delta x^{r+1})$. Thus, we only have to establish upper bounds for the scalar quantities 
   \begin{align*} 
    \left|\frac{\phi_{i,j}^n-P_j^{(i)}(x)}{\phi_{i,j}^n-m_j^{(i)}}\right|, \quad i=1,\dots, N, 
    \end{align*} 
    which can be done by following  \cite[Appendix C]{zhang2017positivity}. We have shown that $\smash{\|\boldsymbol{\tilde{P}}_j(x)-	\boldsymbol{P}_j(x)\| = \mathcal{O}(\Delta x^{r+1})}$ for any $x\in I_j$.}

  {Now we address the limiter \eqref{eq:mod-poly}. For a given $i$, we assume that $\smash{\bar{P}_j^{(i)}(x)}$ is not a constant polynomial. If $\hat{\theta} = 1$, then $\smash{\bar{P}_j^{(i)}(x)=\tilde{P}_j^{(i)}(x)}$ and we have nothing to show. So, let us consider 
\begin{align} \label{thetahatass}
\hat{\theta}=\frac{\phi_{\max}-\phi_{j}^n}{M_j-\phi_{j}^n}. 
\end{align} 
We observe that
\begin{align*}
	\bar{P}_j^{(i)}(x)-\tilde{P}_j^{(i)}(x) &=\hat{\theta} \bigl(\tilde{P}_j^{(i)}(x)-\phi_{i,j}^n \bigr)
	 +\phi_{i,j}^n-\tilde{P}_j^{(i)}(x) = (\hat{\theta}-1) \bigl(\tilde{P}_j^{(i)}(x)-\phi_{i,j}^n \bigr)  \\
	&= \frac{\phi_{\max}-M_j}{M_j-\phi_{j}^n} \bigl(\tilde{P}_j^{(i)}(x)-\phi_{i,j}^n \bigr)
	 = (\phi_{\max}-M_j) \frac{\tilde{P}_j^{(i)}(x)-\phi_{i,j}^n}{M_j-\phi_{j}^n},
\end{align*}
so we obtain 
\begin{align*} 
 \bigl|\bar{P}_j^{(i)}(x)-\tilde{P}_j^{(i)}(x)\bigr|\leq 
 \bigl|\phi_{\max}-M_j \bigr| \left|\frac{\tilde{P}_j^{(i)}(x)-\phi_{i,j}^n}{M_j-\phi_{j}^n}\right|. 
 \end{align*} 
 In the same way as in the analysis of~\eqref{eq:ptilde}, we  notice that the 
  assumption that $\hat{\theta}$ given by~\eqref{thetahatass} satisfies $\hat{\theta} <1$  produces $\phi_{\max}-\phi_{j}^n<M_j-\phi_{j}^n$ which implies $M_j>\phi_{\max}$.}  
    {So, if $x_{*}\in I_j$ is such that $\smash{\sum_{\ell=1}^N\tilde{P}_j^{(\ell)}(x_{*}) = M_j}$, then we have 
   \begin{align*} 
    \phi(x_{*})\leq \phi_{\max}<M_j=\sum_{\ell=1}^N \tilde{P}_j^{(\ell)}(x_{*}),  
    \end{align*}  
    which implies the bound
\begin{align*}
	\bigl|\phi_{\max}-M_j\bigr|&\leq \left|\sum_{\ell=1}^N\tilde{P}_j^{(\ell)}(x_{*})- \phi(x_{*})\right|=\left| \sum_{\ell=1}^N \Big(\tilde{P}_j^{(\ell)}(x_{*})- \phi_{\ell}(x_{*})\Big)\right| =\left| \sum_{\ell=1}^N \Big(\tilde{P}_j^{(\ell)}(x_{*})-{P}_j^{(\ell)}(x_{*})+P_j^{(\ell)}(x_{*})- \phi_{\ell}(x_{*})\Big)\right| \\
	&\leq \sum_{\ell=1}^N \big|\tilde{P}_j^{(\ell)}(x_{*})-{P}_j^{(\ell)}(x_{*})\big|+\sum_{\ell=1}^N \big|P_j^{(\ell)}(x_{*})- \phi_{\ell}(x_{*})\big|\lesssim \bigl\|\boldsymbol{\tilde{P}}_j(x_{*})-	\boldsymbol{P}_j(x_{*})\bigr\| + \bigl\|\boldsymbol{\tilde{P}}_j(x_{*})-\Phi(x_{*})\bigr\| =  \mathcal{O}(\Delta x^{r+1}).
\end{align*}
Thus, we get $\smash{M_j -\phi_{\max}= \mathcal{O}(\Delta x^{r+1})}$. In addition, we notice that 
\begin{align*} 
M_j-\phi_{j}^n\geq \sum_{\ell=1}^{N} (M_j^{(\ell)}-\phi_{\ell,j}^n)\geq M_j^{(i)}-\phi_{i,j}^n, 
\end{align*} 
where $\smash{M_j^{(i)} = \max_{x\in I_j} \tilde{P}_j^{(i)}(x)}$  for each $i=1,\dots,N$. Consequently, 
 we arrive at the bound  
\begin{align*}
	\left|\frac{\tilde{P}_j^{(i)}(x)-\phi_{i,j}^n}{M_j-\phi_{j}^n}\right|\leq \left|\frac{\tilde{P}_j^{(i)}(x)-\phi_{i,j}^n}{M_j^{(i)}-\phi_{i,j}^n}\right|, \quad i=1,\dots, N.
\end{align*}
Therefore, we only have to establish upper bounds for the scalar quantities 
\begin{align*} 
\left|\frac{\tilde{P}_j^{(i)}(x)-\phi_{i,j}^n}{M_j^{(i)}-\phi_{i,j}^n}\right|, \quad i=1,\dots, N, 
\end{align*} 
which can be obtained by following again  \cite[Appendix C]{zhang2017positivity}. In this way, we have proved that 
\begin{align*} 
\bigl\|\boldsymbol{\bar{P}}_j(x)-	\boldsymbol{\tilde{P}}_j(x) \bigr\| = \mathcal{O}(\Delta x^{r+1})
 \quad \text{for any $x\in I_j$.}
\end{align*}} 
  {Finally, we obtain 
\begin{align*} 
\bigl\|\boldsymbol{\bar{P}}_j(x)-\Phi(x) \bigr\|\leq \bigl\|\boldsymbol{\bar{P}}_j(x)-\boldsymbol{\tilde{P}}_j(x) \bigr\|+
 \bigl\|\boldsymbol{\tilde{P}}_j(x)-\boldsymbol{P}_j(x) \bigr\|+ \bigl\|\boldsymbol{P}_j(x)-\Phi_j(x)\bigr\| 
  = \mathcal{O}(\Delta x^{r+1}) \quad \text{for any $x\in I_j$}. \end{align*} 
}

  {
	\begin{remark}\label{rmk:ei}
		Following  \cite[Section 5]{zhang2011maximum}, we point out that the quantities $\smash{m_j^{(i)}}$, $i=1,\dots,N$ and $M_j$ can be computed in a less expensive  way. For instance, if $i\in\{1,\dots,N\}$ is fixed, one can use the Mean Value Theorem to get that 
		\begin{align*} 
		P_j^{(i)} \bigl(\xi^{(i)} \bigr)=\sum_{\alpha=2}^{G-1}\frac{\hat{w}_{\alpha}}{1-2\hat{w}_1} P_j^{(i)}(\hat{x}_j^{\alpha}),
		 \quad \text{for some $\xi^{(i)}\in I_j$}. 
		 \end{align*} 
		 The convex combination $\smash{\phi_{i,j}^n = \sum_{\alpha=1}^{G}   \hat{w}_{\alpha} P_j^{(i)}(\hat{x}_j^{\alpha})}$ implies $\smash{\phi_{i,j}^n = \hat{w}_1  P_j^{(i)}(x_{j-1/2})+ \hat{w}_GP_j^{(i)}(x_{j+1/2})+ (1-2\hat{w}_1)P_j^{(i)}(\xi^{(i)})}$, from which we deduce  that 
		 \begin{align*} 
		 P_j^{(i)}(\xi^{(i)})=\frac{1}{1-2\hat{w}_1}\Big(\phi_{i,j}^n-\hat{w}_1  P_j^{(i)}(x_{j-1/2})- \hat{w}_GP_j^{(i)}(x_{j+1/2})\Big), 
		\end{align*}  
		 i.e. we can compute the value of $\smash{P_j^{(i)}(x)}$ at $x=\xi^{(i)}$ even though the value of $\xi^{(i)}$ is unknown. Thus the minimum value of $\smash{P_j^{(i)}(x)}$ over $I_j$ is computed as 
		\begin{equation}\label{eq:new-min}
			\widehat{m}_{j}^{(i)}=\min\Big\{P_j^{(i)}(\xi^{(i)}),P_j^{(i)}(x_{j-1/2}),P_j^{(i)}(x_{j+1/2})\Big\}, \quad i=1,\dots,N.
		\end{equation}
		Analogously, for the case of the maximum $M_j$, we apply the Mean Value Theorem to the polynomial 
		\begin{align*} 
		Q_j(x) \coloneqq \sum_{i=1}^{N} \tilde{P}_j^{(i)}(x),
		\end{align*} 
		 to get a value $\xi\in I_j$ such that 
		 \begin{align*} 
		 Q_j(\xi) = \sum_{\alpha=2}^{G-1}\frac{\hat{w}_{\alpha}}{1-2\hat{w}_1} Q_j(\hat{x}_j^{\alpha}). 
		 \end{align*} 
		 One can also notice that 
		\[\phi_j^n = \sum_{i=1}^N \phi_{i,j}^n = \sum_{i=1}^N \dfrac{1}{\Delta x} \int_{I_j}\tilde{P}_j^{(i)}(x)
		 \, \mathrm{d}x = \sum_{i=1}^N \sum_{\alpha=1}^{G}   \hat{w}_{\alpha} \tilde{P}_j^{(i)} \bigl(\hat{x}_j^{\alpha} \bigr) =  \sum_{\alpha=1}^{G}   \hat{w}_{\alpha} Q_j \bigl(\hat{x}_j^{\alpha} \bigr),\]
		hence  we can compute $Q_j(\xi)$ by the formula 
		\begin{align*} 
		Q_j(\xi)=\frac{1}{1-2\hat{w}_1}\Big(\phi_{j}^n-\hat{w}_1  Q_j(x_{j-1/2})- \hat{w}_G Q_j (x_{j+1/2})\Big). 
		\end{align*} 
		Consequently, the maximum value of $Q_j (x)$ over $I_j$ is computed as 
		\begin{equation}\label{eq:new-max}
			\widehat{M}_{j}=\max\Big\{Q_j(\xi),Q_j(x_{j-1/2}),Q_j(x_{j+1/2})\Big\}. 
		\end{equation}
		So, in the numerical computations we actually use the values:
		\begin{equation}
			\theta_i\coloneqq \min\Biggl\{\dfrac{\phi_{i,j}^n-\delta}{\phi_{i,j}^n-\widehat{m}_j^{(i)}},1\Biggr\},\, i=1,\dots, N, \quad \text{and}\quad \hat{\theta} \coloneqq 
			\min\Biggl\{\Biggl|\dfrac{\phi_{\mathrm{max}}-\phi_{j}^n}{\widehat{M}_j-\phi_{j}^n}\Biggr|,1\Biggr\}.
		\end{equation}
		First we observe that $$\smash{\phi_{i,j}^n = \hat{w}_1  P_j^{(i)}(x_{j-1/2})+ \hat{w}_GP_j^{(i)}(x_{j+1/2})+ (1-2\hat{w}_1)P_j^{(i)}(\xi^{(i)})}\geq \hat{w}_1	\widehat{m}_{j}^{(i)}+  \hat{w}_G\widehat{m}_{j}^{(i)}+ (1-2\hat{w}_1)\widehat{m}_{j}^{(i)} = \widehat{m}_{j}^{(i)},$$
		thus, $\theta_i\in [0,1],$ $i=1,\dots,N$. In addition, we have that \[P_j^{(i)} \bigl(\xi^{(i)} \bigr)=\sum_{\alpha=2}^{G-1}\frac{\hat{w}_{\alpha}}{1-2\hat{w}_1} P_j^{(i)}(\hat{x}_j^{\alpha})\geq  \dfrac{m_j^{(i)}}{1-2\hat{w}_1}\sum_{\alpha=2}^{G-1}\hat{w}_{\alpha}=\dfrac{m_j^{(i)}}{1-2\hat{w}_1}(1-2\hat{w}_1)=m_j^{(i)},\] hence $m_j^{(i)}\leq \widehat{m}_j^{(i)}$. In this way we get that \[    \frac{\phi_{i,j}^n-\delta}{\phi_{i,j}^n-m_j^{(i)}}\leq \frac{\phi_{i,j}^n-\delta}{\phi_{i,j}^n-\widehat{m}_j^{(i)}}, \quad i=1,\dots, N. \] Analogously, we notice that $\widehat{M}_j\leq M_j$ and therefore \[	\left|\frac{\phi_{\max}-\phi_{j}^n}{M_j-\phi_{j}^n}\right|\leq 	\left|\frac{\phi_{\max}-\phi_{j}^n}{\widehat{M}_j-\phi_{j}^n}\right|.\]
		So, the limiters with the values \eqref{eq:new-min} and \eqref{eq:new-max} are more relaxed than \eqref{eq:ptilde}
		 and \eqref{eq:phat}, respectively. Therefore the analysis of accuracy performed above is still valid.
\end{remark}}

\subsection{Time discretization}\label{subsec:rk}
High-order time integration is required  to ensure the appropriate  order of convergence of the method. In this work, a strong-stability preserving (SSP) third-order TVD Runge-Kutta time discretization approach is used \cite{liu1994weighted}. This time stepping method could be described as follows. Assume that $\Phi^n$ is the vector of approximate solutions  of $\Phi'=\mathcal{L}(\Phi)$ at $t=t_n$. Then the approximate values $\Phi^{n+1}$ associated with $t_{n+1}=t_n+\Delta t$ are calculated by
\begin{equation}\label{eq:TVDRK3}
	\begin{aligned}
		\Phi^{(1)}&=\Phi^{n}+\Delta t \mathcal{L}(\Phi^n),\\
		\Phi^{(2)}&=\dfrac{3}{4}\Phi^{n}+\dfrac{1}{4}\Phi^{(1)}+\dfrac{1}{4}\Delta t \mathcal{L}(\Phi^{(1)}),\\
		\Phi^{n+1}&=\dfrac{1}{3}\Phi^{n}+\dfrac{2}{3}\Phi^{(2)}+\dfrac{2}{3}\Delta t \mathcal{L}(\Phi^{(2)}).
	\end{aligned}
\end{equation}
We refer to \eqref{eq:TVDRK3} as SSPRK (3,3). The SSP time discretization methods are widely used for hyperbolic PDE because they preserve  nonlinear stability properties  necessary for problems with non-smooth solutions. On the other hand, due to convexity, the intermediate stages of the SSPRK methods have SSP properties (i.e., $\| \Phi^{n}\|\leq \| \Phi^{n-1}\|$ for the internal stages).  Consequently,  the present finite volume WENO scheme with this time discretization will still satisfy the maximum principle. Since  it is necessary to evaluate three times the operator $\mathcal{L}(\cdot)$   to advance one time step,  the effective SSP coefficient for SSPRK (3,3) the method 
 (which is defined as in \cite{gottlieb2009high};  namely    the SSP coefficient of the method divided by the number of stages) equals~$1/3$.

To satisfy the CFL condition \eqref{eq:CFLMMP} the time step  $\Delta t$ is computed adaptively for each time   step~$n$. More specifically, the solution $\Phi^{n +1}$ at
$t_{n + 1} = t_n + \Delta t$ is calculated from $\Phi^{n}$ by using   $\Delta t =  \mu (\hat{w}_1  \Delta x /\alpha)$, where $\mu$ is given by \eqref{eq:mu}. 
\section{Numerical examples}\label{sec:numexa}

\begin{figure}[t] 
	\centering
			\includegraphics{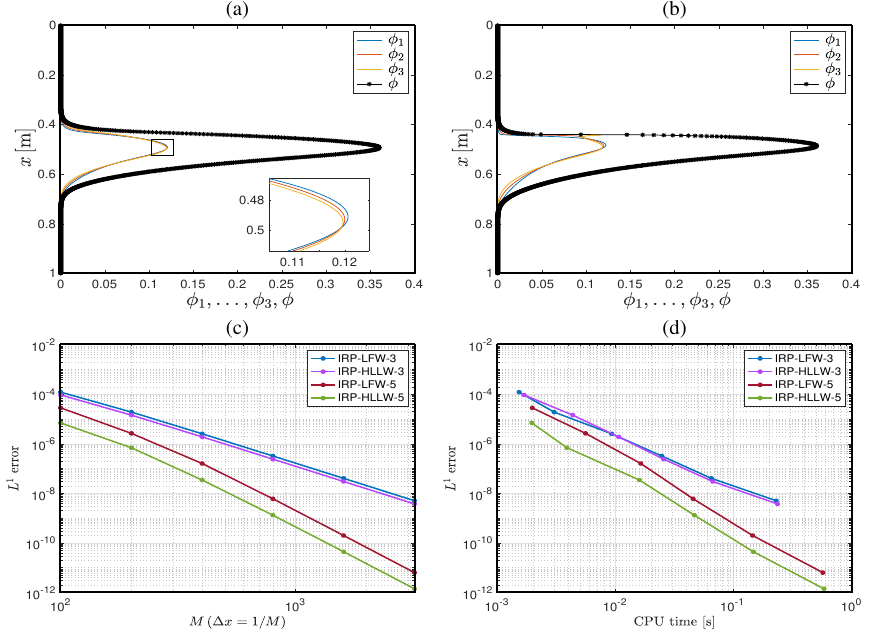} 
			\caption{Example 1 (MLB Model, $N = 3$): numerical results obtained by scheme IRP-HLLW-5  with $M = 1600$ at simulated time 
			 (a) $T = 5 \, \mathrm{s}$, (b) $T = 10 \, \mathrm{s}$,   (c) approximate $L^1$ errors for all schemes tested as function of~$M$,  (d) efficiency plot obtained for discretization levels $\Delta x = 1/M$ with $M = 100$, $200$, $400$,
		$800$, and~$1600$.}
	\label{fig:ex1-1}
\end{figure}

\begin{table}[t] \caption{Example 1 (MLB model, $N = 3$): $L^1$ errors and numerical order for IRP-LFW-3,  IRP-HLLW-3, IRP-LFW-5,  and IRP-HLLW-5 schemes applied
		to smooth initial conditions for $T = 5 \,\mathrm{s}$ (before shock formation) and $T = 10\, \mathrm{s}$ (after shock
		formation)   {for values of $\varepsilon = \Delta x^2$ and $\varepsilon = 10^{-6}$}. The reference solution is computed by IRP-HLLW-5 with $M_{\mathrm{ref}} = 12800.$}
	\begin{center} 
		{\footnotesize   \addtolength{\tabcolsep}{-2pt}
	  {		  
	\begin{tabular}{ccccccccccccc}\toprule
		\multicolumn{13}{c}{$\varepsilon = \Delta x^2$}\\
		\toprule
	$M$&$e_M^{\mathrm{tot}}$& $\theta_M$ & cpu [s] & $e_M^{\mathrm{tot}}$ & $\theta_M$& cpu [s]&$e_M^{\mathrm{tot}}$&$\theta_M$ &cpu [s] &$e_M^{\mathrm{tot}}$& $\theta_M$ &cpu [s]
	 \\  \cmidrule(lr){2-4}\cmidrule(lr){5-7}\cmidrule(lr){8-10}\cmidrule(lr){11-13}
	& \multicolumn{3}{c}{IRP-LFW-3, $T=5 \,\mathrm{s}$} & \multicolumn{3}{c}{IRP-LFW-3, $T=10 \, \mathrm{s}$}
	 & \multicolumn{3}{c}{IRP-HLLW-3, $T=5 \,\mathrm{s}$} & \multicolumn{3}{c}{IRP-HLLW-3, $T=10 \, \mathrm{s}$}	\\ \midrule
		100   & 1.21e-04  &  --  &  1.52e-03 & 7.90e-04 &  ---& 3.18e-03 &9.32e-05  &   --- &  1.67e-03  & 6.73e-04 & --- & 2.61e-03 \\ 
		200   & 1.90e-05 & 2.6 &  3.01e-03 & 3.51e-04 &  1.2 & 7.06e-03& 1.43e-05 &   2.7 & 4.33e-03  & 2.95e-04 &  1.2 & 6.21e-03 \\ 
		400   & 2.53e-06 & 2.9&  9.25e-03 & 1.58e-04 &  1.1 & 1.64e-02 & 1.89e-06 &  2.9 &  1.06e-02  & 1.30e-04 &  1.2  & 1.56e-02 \\ 
		800   &   3.2e-07 & 3.0 &  2.47e-02& 6.42e-05 &  1.3 & 4.08e-02 & 2.39e-07 & 3.0 &  2.55e-02 & 5.08e-05 &  1.3 &  4.10e-02 \\ 
		1600 & 4.01e-08 &  3.0 & 6.50e-02 & 2.39e-05 &  1.4 & 0.12 & 2.99e-08 & 3.0 & 6.67-02    & 1.81e-05  &  1.5  &  0.12 \\ 
		3200 & 5.01e-09 &  3.0 & 0.23         & 1.57e-05 & 0.6 &  0.47        & 3.73e-09 &  3.0 &   0.24        & 1.28e-05  & 0.5  &   0.48\\ \midrule 
		& \multicolumn{3}{c}{IRP-LFW-5, $T=5 \,\mathrm{s}$} & \multicolumn{3}{c}{IRP-LFW-5, $T=10 \, \mathrm{s}$}
	 & \multicolumn{3}{c}{IRP-HLLW-5, $T=5 \,\mathrm{s}$} & \multicolumn{3}{c}{IRP-HLLW-5, $T=10 \, \mathrm{s}$}  
	 \\ \midrule
		100   & 2.81e-05  &  --- &  1.96e-03 &5.95e-04 & ---   & 6.16e-03 & 6.92e-06  &  --- & 1.95e-03  &2.38e-04 &  ---  & 4.77e-03\\ 
		200   & 2.62e-06 & 3.4 &  5.57e-03 &2.92e-04 & 1.0   & 1.31e-02  & 6.9e-06  & 3.3 & 3.87e-03  &1.09e-04 &  1.1   & 1.32e-02 \\ 
		400   & 1.59e-07  & 4.0 &  1.63e-02 &1.20e-04 &  1.3   & 2.99e-02 & 3.43e-08  &4.3 &  1.59e-02  &4.80e-05 & 1.2   & 3.02e-02\\ 
		800   & 6.0e-09   & 4.7  & 4.55e-02 &4.21e-05 &  1.5   & 7.85e-02 & 1.3e-09    & 4.7 &  4.66e-02 &1.74e-05  & 1.5   & 8.07e-02  \\ 
		1600 & 1.97e-10  & 4.9  &   0.14        &1.18e-05  &  1.8   &  0.29        & 4.37e-11   & 4.9 &   0.15         &5.27e-06  & 1.7   & 0.29        \\ 
		3200 & 6.29e-12 & 5.0  &   0.57       & 9.94e-06 & 0.2   &    1.14       & 1.4e-12  & 5.0 &  0.58         & 5.15e-06  & 0.0  & 1.15          \\ \toprule
		\multicolumn{13}{c}{$\varepsilon = 10^{-6}$}\\
		\toprule
			$M$&$e_M^{\mathrm{tot}}$& $\theta_M$ & cpu [s] & $e_M^{\mathrm{tot}}$ & $\theta_M$& cpu [s]&$e_M^{\mathrm{tot}}$&$\theta_M$ &cpu [s] &$e_M^{\mathrm{tot}}$& $\theta_M$ &cpu [s]
		\\  \cmidrule(lr){2-4}\cmidrule(lr){5-7}\cmidrule(lr){8-10}\cmidrule(lr){11-13}
		& \multicolumn{3}{c}{IRP-LFW-3, $T=5 \,\mathrm{s}$} & \multicolumn{3}{c}{IRP-LFW-3, $T=10 \, \mathrm{s}$}
		& \multicolumn{3}{c}{IRP-HLLW-3, $T=5 \,\mathrm{s}$} & \multicolumn{3}{c}{IRP-HLLW-3, $T=10 \, \mathrm{s}$}	\\ \midrule
		100    &2.96e-04 & ---&  1.08e-03  &1.04e-03  & ---   & 3.17e-03 &2.2e-04  &  --- & 1.37e-03&8.63e-04 & ---  & 3.09e-03\\
		200   &5.03e-05 & 2.7&  3.89e-03  &4.29e-04 &  1.3  & 7.23e-03 &3.72e-05& 2.6 & 3.39e-03&3.54e-04 &  1.3 & 7.23e-03\\
		400   & 4.33e-06& 3.5&  1.03e-02  &1.68e-04  &  1.4  & 1.65e-02 &3.12e-06& 3.6 & 9.27e-03&1.37e-04  &  1.4 & 1.64e-02\\
		800   & 3.57e-07 & 3.6&  2.55e-02 &6.48e-05  & 1.4   &  4.1e-02  &2.64e-07& 3.5 & 2.48e-02&5.12e-05 &  1.4 &  4.15e-02\\
		1600 &3.36e-08 & 3.4&  6.55e-02  &2.33e-05  & 1.5   &  0.12        &2.55e-08& 3.4 & 6.64e-02&1.77e-05 & 1.5  &   0.12\\
		3200 &3.75e-09 & 3.2&   0.23         &1.51e-05   & 0.6 &  0.47         &2.89e-09& 3.1  &   0.24      &1.23e-05 & 0.5  &    0.48\\ \midrule 
		& \multicolumn{3}{c}{IRP-LFW-5, $T=5 \,\mathrm{s}$} & \multicolumn{3}{c}{IRP-LFW-5, $T=10 \, \mathrm{s}$}
		& \multicolumn{3}{c}{IRP-HLLW-5, $T=5 \,\mathrm{s}$} & \multicolumn{3}{c}{IRP-HLLW-5, $T=10 \, \mathrm{s}$}  
		\\ \midrule
	    100    &2.81e-05 &  --- & 2.07e-03& 6.58e-04&  --- & 4e-03       &5.75e-06 &  --- &  2.8e-03&2.68e-04  &   --- & 4.58e-03\\
		200   &2.62e-06 & 3.4 & 8.18e-03 &3.1e-04   & 1.1   & 8.67e-03  &6.86e-07 & 3.1  &  7.13e-03&1.15e-04   &  1.2  & 1.02e-02 \\
		400   & 1.64e-07 & 4.0 & 1.08e-02 &1.24e-04 &1.3   & 2.06e-02 &3.55e-08 & 4.3 & 1.83e-02 &4.89e-05  & 1.2  &  2.17e-02\\
		800   & 6.06e-09& 4.8 &4.37e-02 &4.24e-05 &1.5   & 9.56e-02 &1.32e-09 & 4.7  & 4.68e-02&1.75e-05   & 1.5  &  9.01e-02\\
		1600 &1.91e-10   &5.0  & 0.17        & 1.18e-05  & 1.8  &  0.38         &4.28e-11   & 5.0  &   0.14&5.23e-06  &  1.7  &   0.34 \\
		3200 &5.76e-12 & 5.1  &   0.63     &9.84e-06  & 0.3&    1.45       &1.33e-12  & 5.0   &   0.58       &5.1e-06      & 0.0  &    1.44\\ \bottomrule 
	\end{tabular}}} 
	\end{center} 
	\label{tab:tab1}
\end{table}

\begin{table}[t] \caption{Example 1 (MLB model, $N = 3$): minimum of the solutions $\phi_{i,j}^n$, $i=1,\dots,3$, and maximum of the solution $\phi_{j}^n$ obtained by schemes LFW-3, LFW-5 (without limiters) and  IRP-LFW-3, IRP-LFW-5 (with limiters) until specified times $T$.}
	\begin{center} 
		 {
		{\footnotesize  \addtolength{\tabcolsep}{-2pt} 
			\begin{tabular}{ccccccccc}\toprule
				$M$& $\underline{\phi}_{,1}$  & $\underline{\phi}_{,2}$ & $\underline{\phi}_{,3}$  &  $\overline{\phi}$  & $\underline{\phi}_{,1}$   & $\underline{\phi}_{,2}$ & $\underline{\phi}_{,3}$  &  $\overline{\phi}$
				\\  \cmidrule(lr){2-5}\cmidrule(lr){6-9} 
				& \multicolumn{4}{c}{LFW-3, $T=10 \,\mathrm{s}$} & \multicolumn{4}{c}{IRP-LFW-3, $T=10\, \mathrm{s}$}
				\\ \midrule
				100 &  -1.846e-04& -8.838e-10& -4.313e-10&0.3586& 7.32e-23&  7.32e-23&  7.32e-23& 0.3586\\ 
				200 &   1.659e-23 &   1.651e-23&   1.595e-23 &   0.3598& 3.969e-23& 3.969e-23& 3.969e-23& 0.3598\\ 
				400  & 2.203e-23&  2.226e-23&  2.258e-23& 0.3603 &3.002e-23& 3.002e-23& 3.002e-23& 0.3603\\ 
				800  &   2.293e-23&  2.307e-23&  2.327e-23& 0.3603&2.629e-23& 2.629e-23& 2.629e-23& 0.3603\\ 
				1600  & 2.307e-23&  2.314e-23&  2.324e-23& 0.3604 &2.465e-23& 2.465e-23& 2.465e-23& 0.3604 
				\\ \midrule 
				& \multicolumn{4}{c}{LFW-5, $T=10 \,\mathrm{s}$} & \multicolumn{4}{c}{IRP-LFW-5, $T=10 \, \mathrm{s}$}
				\\ \midrule
				100 & -1.515e-07& 5.313e-23& 5.313e-23& 0.3597 & 5.352e-23&  7.32e-23&  7.32e-23& 0.3596\\ 
				200 & 2.785e-23& 2.789e-23& 2.788e-23& 0.3601 & 2.271e-23& 3.013e-23& 3.969e-23& 0.3601\\
				400  &2.456e-23& 2.462e-23&2.469e-23& 0.3603 & 2.585e-23& 2.703e-23& 2.852e-23& 0.3603\\
				800  &2.374e-23& 2.378e-23& 2.382e-23& 0.3603 & 2.475e-23&   2.5e-23& 2.547e-23& 0.3603\\
				1600  & 2.341e-23& 2.343e-23& 2.345e-23& 0.3604 & 2.394e-23& 2.405e-23& 2.424e-23& 0.3604\\ \bottomrule 
		\end{tabular}}} 
	\end{center} 
	\label{tab:tab11}
\end{table}

\begin{figure}[t] 
	\centering
			\includegraphics[scale=0.88]{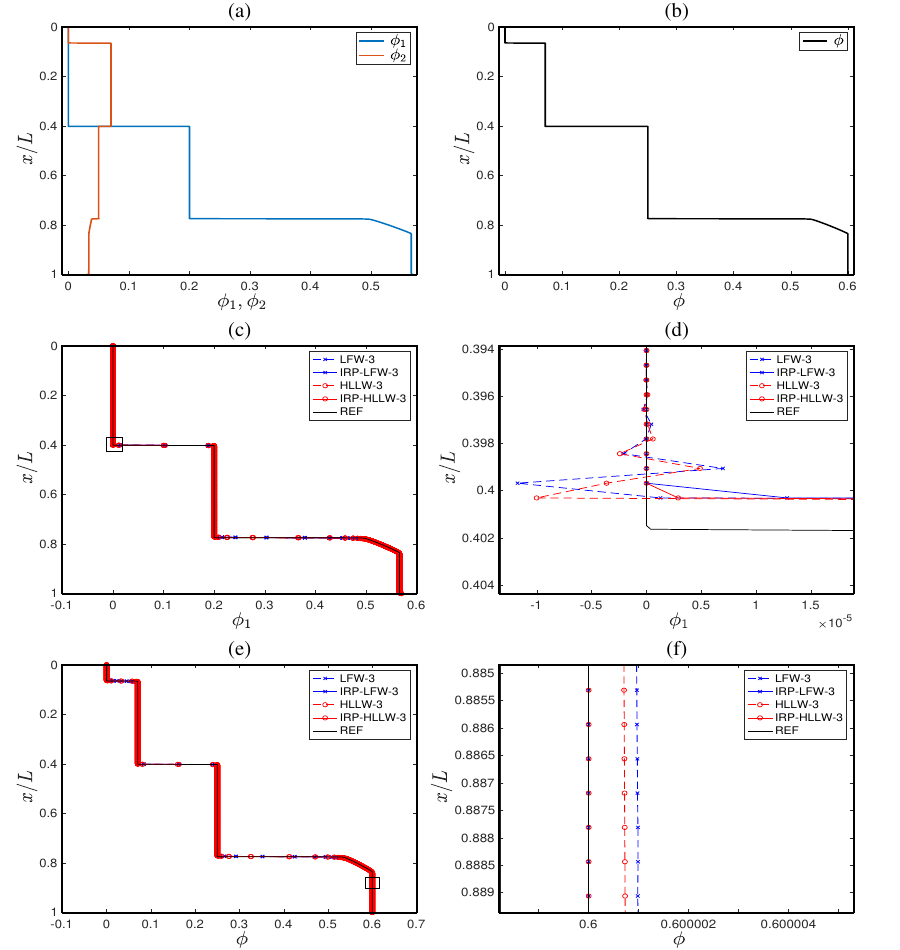}  
				\caption{Example 2 (MLB model, $N = 2$): reference solution for (a) $\phi_1$, $\phi_2$ and (b) 
				 $\phi$   at $T = 50 \, \mathrm{s}$  computed by scheme IRP-HLLW-3 with $M_{\mathrm{ref}} = 6400$, 
				  and   comparison of schemes  for (c)  $\phi_1$, (d) enlarged view of~(c), (e) $\phi$, 
				   and (f) enlarged view of~(e)  at $T = 50 \, \mathrm{s}$  with $M = 1600$.}
	\label{fig:fig2}
\end{figure}

\begin{figure}[t] 
	\centering
	\includegraphics{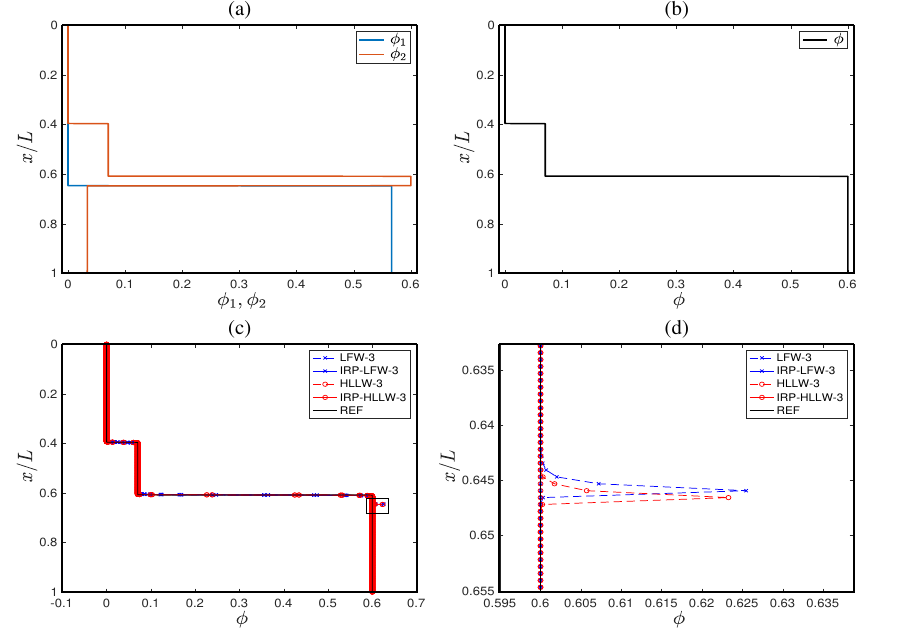}  
		\caption{Example 2 (MLB model, $N = 2$): reference solution for (a) $\phi_1$, $\phi_2$ and (b) 
				 $\phi$   at $T = 300 \, \mathrm{s}$  computed by scheme IRP-HLLW-3 with $M_{\mathrm{ref}} = 6400$, 
				  and   (c) comparison of schemes  for   $\phi_1$, (d) enlarged view of~(c), 
				   at $T = 300 \, \mathrm{s}$  with $M = 1600$.}
	\label{fig:fig3}
\end{figure}

\begin{figure}[t] 
	\centering
	\includegraphics{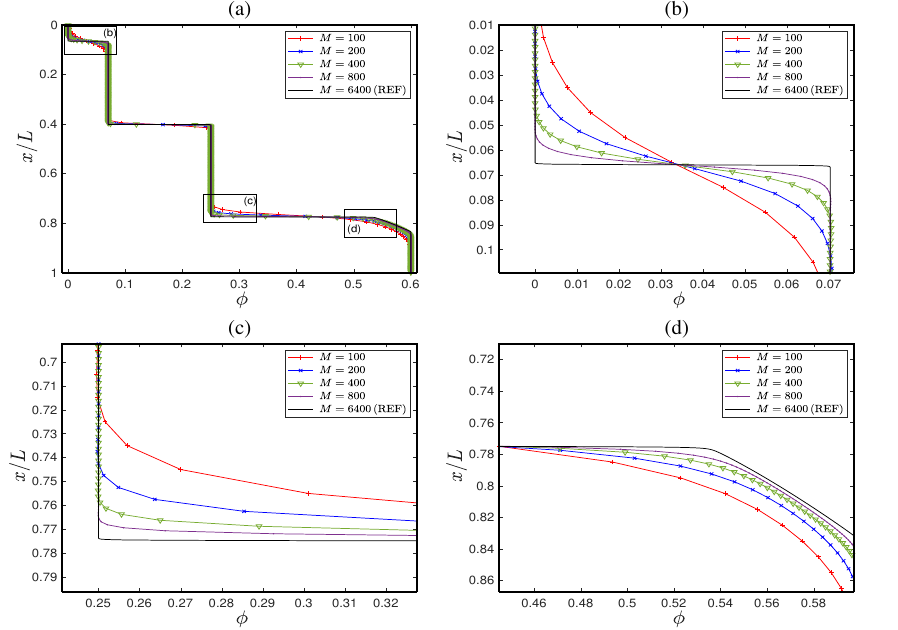}  
		\caption{Example 2 (MLB model, $N = 2$): (a) comparison of the numerical solutions  for  $\phi$  at $T = 50\, \mathrm{s}$ 
		  computed by IRP-HLLW-3 with various step sizes, (b)--(d)    enlarged views of~(a).}
	\label{fig:fig9}
\end{figure}

\begin{table}[t] \caption{Example 2 (MLB model, $N = 2$): $L^1$ errors, numerical order, and CPU time (seconds) for IRP-LFW-3 and IRP-HLLW-3  schemes at $T = 50\,\mathrm{s}$ and $T = 300\, \mathrm{s}$. The reference solution is computed by IRP-HLLW-3 with $M_{\mathrm{ref}} = 6400.$}
	\begin{center} 
		{\footnotesize \addtolength{\tabcolsep}{-2pt} 
	\begin{tabular}{ccccccccccccc}\toprule
	$M$&$e_M^{\mathrm{tot}}$& $\theta_M$ & cpu [s] & $e_M^{\mathrm{tot}}$ & $\theta_M$& cpu [s]&$e_M^{\mathrm{tot}}$&$\theta_M$ &cpu [s] &$e_M^{\mathrm{tot}}$& $\theta_M$ &cpu [s]
	 \\  \cmidrule(lr){2-4}\cmidrule(lr){5-7}\cmidrule(lr){8-10}\cmidrule(lr){11-13}
	& \multicolumn{3}{c}{IRP-LFW-3, $T=50 \,\mathrm{s}$} & \multicolumn{3}{c}{IRP-LFW-3, $T=300 \, \mathrm{s}$}
	 & \multicolumn{3}{c}{IRP-HLLW-3, $T=50 \,\mathrm{s}$} & \multicolumn{3}{c}{IRP-HLLW-3, $T=300 \, \mathrm{s}$}	\\ \midrule
		100 & 6.56e-03 &  ---   & 3.89e-02& 9.50e-03&   --- & 1.23e-01 & 4.78e-03&   ---& 5.99e-02 &4.93e-03& ---& 2.08e-01 \\
		200 & 3.83e-03 &   0.7 & 6.98e-02  &4.79e-03& 0.9& 3.54e-01&2.83e-03& 0.7&1.26e-01 &2.10e-03&  1.2& 6.67e-01 \\
		400 & 2.05e-03 &   0.9 &  2.33e-01 & 2.35e-03&  1.0&  1.36 &1.51e-03& 0.9&  5.04e-01  & 9.31e-04&  1.1&  2.59 \\
		800 & 1.02e-03 &   1.0 & 0.95e-01& 9.69e-04&  1.2&  5.56 &7.34e-04&  1.0&    2.04&4.59e-04& 1.0&  10.6 \\
		1600 & 4.67e-04 &   1.1 & 3.83 & 5.39e-04&  0.8&  21.8& 3.15e-04&  1.2&   8.08 &2.44e-04& 0.9&  41.7 \\
				 \bottomrule
	\end{tabular}} 
	\end{center} 
	\label{tab:tab3}
\end{table}

\begin{table}[t] \caption{Example 2 (MLB model, $N = 2$): minimum of the solutions $\phi_{i,j}^n$, $i=1,2$, and maximum of the solution $\phi_{j}^n$ obtained by 
 schemes LFW-3 (without limiters),  IRP-LFW-3 (with limiters),  HLLW-3 (without limiters), and IRP-HLLW-3 (with limiters) until specified times $T$.}
	\begin{center} 
		{\footnotesize  \addtolength{\tabcolsep}{-2pt} 
	\begin{tabular}{ccccccc}\toprule
	$M$& $\underline{\phi}_{,1}$  & $\underline{\phi}_{,2}$  &  $\overline{\phi}$  & $\underline{\phi}_{,1}$   & $\underline{\phi}_{,2}$  &  $\overline{\phi}$
	 \\  \cmidrule(lr){2-4}\cmidrule(lr){5-7} 
	& \multicolumn{3}{c}{LFW-3, $T=50 \,\mathrm{s}$} & \multicolumn{3}{c}{IRP-LFW-3, $T=50 \, \mathrm{s}$}
		\\ \midrule
		100 & -4.122e-004&    9.755e-004 & 0.600045& 2.056e-025& 9.763e-004 & 0.600000\\ 
		200 &  -2.043e-004&  -2.236e-004& 0.600026& 1.908e-039&7.263e-013& 0.600000\\
		400  &-1.013e-004&  -2.503e-004 &0.600016& 1.442e-067&  3.434e-013& 0.600000\\
		800  & -5.030e-005 &  -1.365e-004 & 0.600010 & 7.517e-124 & 1.163e-013& 0.600000\\
		1600  & -2.513e-005& -6.839e-005 & 0.600010& 2.073e-236& 2.049e-014& 0.600000
		 \\ \midrule 
		 & \multicolumn{3}{c}{LFW-3, $T=300 \,\mathrm{s}$} & \multicolumn{3}{c}{IRP-LFW-3, $T=300 \, \mathrm{s}$}
		\\ \midrule 
		100 &-4.130e-004&   -1.076e-003 & 0.618048&1.061e-046& 1.648e-013& 0.600000\\
		200 &  -2.047e-004&    -5.425e-004 & 0.654617 & 3.884e-083& 4.135e-014& 0.600000\\
		400  & -1.013e-004&  -2.724e-004 & 0.637268&1.829e-155& 3.924e-015& 0.600000\\
		800  & -5.035e-005&  -1.365e-004 & 0.650888&4.012e-300& 5.758e-017& 0.600000\\
		1600  & -2.513e-005& -6.839e-005 & 0.625597& 1.729e-322& 2.228e-020& 0.600000 \\ \midrule 
		 & \multicolumn{3}{c}{HLLW-3, $T=50 \,\mathrm{s}$} & \multicolumn{3}{c}{IRP-HLLW-3, $T=50 \, \mathrm{s}$}
		\\ \midrule  
		100 & -4.032e-004 & -2.881e-004 & 0.600047 & 1.767e-025&  2.279e-09& 0.600000\\
		200 & -1.994e-004  & -2.598e-004 & 0.600029 & 1.623e-039& 3.968e-014& 0.600000\\
		400  &  -9.873e-005   &  -1.308e-004 &  0.600020 &  1.222e-067& 4.873e-017& 0.600000\\
		800  & -4.895e-005  & -6.560e-005& 0.600024 &  6.608e-124& 3.263e-022& 0.600000\\
		1600  & -2.442e-005& -3.288e-005 & 0.600034 &  1.771e-236& 4.069e-032& 0.600000\\\midrule 
		& \multicolumn{3}{c}{HLLW-3, $T=300 \,\mathrm{s}$} & \multicolumn{3}{c}{IRP-HLLW-3, $T=300 \, \mathrm{s}$} \\ \midrule 
		100 &-4.032e-004 &  -5.196e-004& 0.656865&3.090e-047&   7.971e-020& 0.600000\\
		200 & -1.991e-004 & -2.609e-004 & 0.666836 &3.593e-083&  2.388e-027& 0.600000\\
		400  & -9.873e-005 & -1.308e-004 & 0.642500&1.635e-155& 5.389e-042& 0.600000\\
		800  &  -4.895e-005 & -6.560e-005 & 0.646103&3.455e-300&  6.799e-071& 0.600000\\
		1600  &-2.442e-005& -3.288e-005 & 0.623364& 2.734e-322& 2.678e-128& 0.600000
		\\ \bottomrule 
	\end{tabular}} 
	\end{center} 
	\label{tab:tab4}
\end{table}

\begin{figure}[t] 
	\centering
			\includegraphics[scale=0.88]{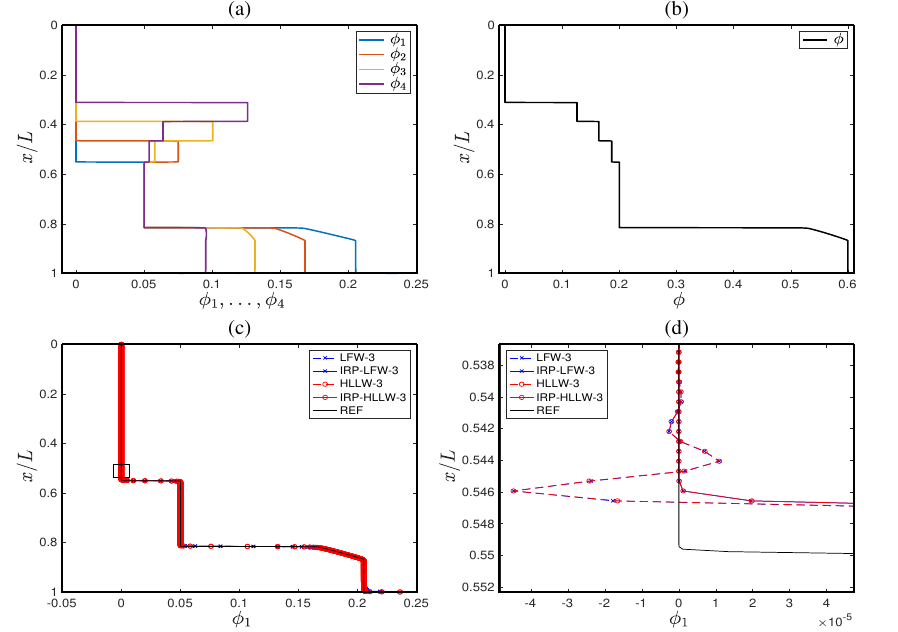}   
		\caption{Example 3 (MLB model, $N = 4$): reference solution for (a) $\phi_1, \dots, \phi_4$ and (b) 
				 $\phi$   at $T = 50 \, \mathrm{s}$  computed by scheme IRP-HLLW-3 with $M_{\mathrm{ref}} = 6400$, 
				  and   comparison of schemes  for (c)  $\phi_1$, (d) enlarged view of~(c) at $T = 50 \, \mathrm{s}$  with $M = 1600$.}
	\label{fig:fig4}
\end{figure}

\begin{figure}[t] 
	\centering
	\includegraphics{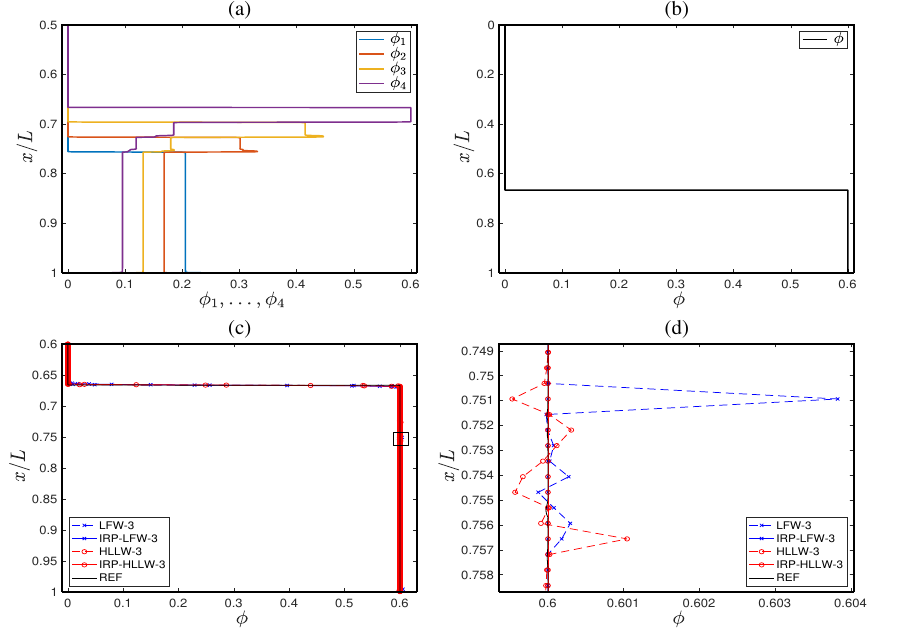}  
		\caption{Example 3 (MLB model, $N = 4$): reference solution for (a) $\phi_1, \dots, \phi_4$ and (b) 
				 $\phi$   at $T = 300 \, \mathrm{s}$  computed by scheme IRP-HLLW-3 with $M_{\mathrm{ref}} = 6400$, 
				  and  (c)  comparison of schemes  for   $\phi$, (d) enlarged view of~(c)
				    at $T = 300 \, \mathrm{s}$  with $M = 1600$.}
	\label{fig:fig5}
\end{figure}

\begin{table}[t] \caption{Example 3 (MLB model, $N = 4$): $L^1$ errors, numerical order, and CPU time (seconds) for IRP-LFW-3 and IRP-HLLW-3  schemes at $T = 50\,\mathrm{s}$ and $T = 300\, \mathrm{s}$. The reference solution is computed by IRP-HLLW-3 with $M_{\mathrm{ref}} = 6400.$}
	\begin{center} 
		{\footnotesize   \addtolength{\tabcolsep}{-2pt} 
	\begin{tabular}{ccccccccccccc}\toprule
	$M$&$e_M^{\mathrm{tot}}$& $\theta_M$ & cpu [s] & $e_M^{\mathrm{tot}}$ & $\theta_M$& cpu [s]&$e_M^{\mathrm{tot}}$&$\theta_M$ &cpu [s] &$e_M^{\mathrm{tot}}$& $\theta_M$ &cpu [s]
	 \\  \cmidrule(lr){2-4}\cmidrule(lr){5-7}\cmidrule(lr){8-10}\cmidrule(lr){11-13}
	& \multicolumn{3}{c}{IRP-LFW-3, $T=50 \,\mathrm{s}$} & \multicolumn{3}{c}{IRP-LFW-3, $T=300 \, \mathrm{s}$}
	 & \multicolumn{3}{c}{IRP-HLLW-3, $T=50 \,\mathrm{s}$} & \multicolumn{3}{c}{IRP-HLLW-3, $T=300 \, \mathrm{s}$}	\\ \midrule
		100 & 6.89e-03 &  ---   & 5.24e-02 & 1.83e-03 &  ---  & 1.75e-01 & 5.60e-03&  ---  & 8.17e-02 &1.75e-03&   ---  & 3.09e-01 \\
		200 & 3.54e-03 &   0.9 & 1.04e-01  &1.15e-03 & 0.7 & 5.37e-01  &2.86e-03&1.0 & 1.94e-01 &6.69e-04&1.4& 1.03  \\
		400 & 1.93e-03 &   0.8 & 4.11e-01  & 6.07e-04 & 0.9& 2.15&1.54e-03&0.9& 7.88e-01 &3.36e-04&0.9&  4.29\\
		800 & 1.05e-03 &   0.8 & 1.61& 2.82e-04 & 1.1& 8.59&8.48e-04&0.8& 3.21&1.64e-04&1.1 & 16.6\\
		1600 & 4.87e-04 &   1.1 & 6.31& 1.29e-04 & 1.1& 33.6& 3.78e-04&1.1& 12.8&9.00e-05&0.8 &67.6\\ 
				 \bottomrule
	\end{tabular}} 
	\end{center} 
	\label{tab:tab6}
\end{table}

\begin{table}[t] \caption{Example 3 (MLB model, $N = 4$): minimum of the solutions $\phi_{i,j}^n$, $i=1,\dots,4$, and maximum of the solution $\phi_{j}^n$ obtained by schemes LFW-3 (without limiters) and IRP-LFW-3 (with limiters) until specified times $T$.}
	\begin{center} 
		{\footnotesize   \addtolength{\tabcolsep}{-2pt} 
	\begin{tabular}{ccccccccccc}\toprule
	$M$& $\underline{\phi}_{,1}$  & $\underline{\phi}_{,2}$  &  $\underline{\phi}_{,3}$  & $\underline{\phi}_{,4}$  &   $\overline{\phi}$
	& $\underline{\phi}_{,1}$  & $\underline{\phi}_{,2}$  &  $\underline{\phi}_{,3}$  & $\underline{\phi}_{,4}$  &   $\overline{\phi}$
	 \\  \cmidrule(lr){2-6} \cmidrule(lr){7-11}
	& \multicolumn{5}{c}{LFW-3, $T=50 \,\mathrm{s}$} &  \multicolumn{5}{c}{IRP- LFW-3, $T=50 \,\mathrm{s}$}  \\ \midrule
		100 & -6.109e-004 &  -5.443e-004&    -4.791e-004&-4.305e-004  & 0.599930 & 8.428e-058&  1.248e-024&  2.774e-014& 6.634e-08& 0.598892 \\ 
		200 &  -3.687e-004 &  -3.199e-004&    -2.806e-004 & -2.157e-004 & 0.604964 &  2.437e-118&  1.606e-047&  1.825e-026& 2.213e-013& 0.599973 \\
		400  &-1.845e-004 &   -1.598e-004 &  -1.401e-004 & -1.080e-004 & 0.607519 &1.817e-239&  3.064e-093&  7.809e-051& 3.408e-024&  0.600000 \\
		800  & -9.229e-005 & -7.978e-005 &  -6.993e-005 & -5.411e-005 &0.609175 & 1.828e-322&  1.303e-184&  1.687e-099& 1.232e-045& 0.600000 \\
		1600  & -4.610e-005 & -3.983e-005 &  -3.491e-005 & -2.709e-005 &0.605247 & 1.976e-322& 3.112e-322& 1.033e-196& 2.660e-088& 0.600000 
				 \\ \midrule 
		& \multicolumn{5}{c}{LFW-3, $T=300 \,\mathrm{s}$}  & \multicolumn{5}{c}{IRP- LFW-3, $T=300 \,\mathrm{s}$}   \\ \midrule
		100 &-8.523e-004 & -6.713e-004 &  -5.662e-004 &  -4.305e-004 & 0.666358  & 1.172e-203&  5.930e-064&  5.930e-041&  1.990e-025& 0.600000 \\ 
		200 & -3.908e-004 &  -3.403e-004 & -2.806e-004 & -2.157e-004 & 0.609710 &  2.173e-322&  8.387e-128&  3.934e-081&  7.523e-050&      0.600000\\
		400  &-2.070e-004 & -1.701e-004 & -1.401e-004 & -1.080e-004  & 0.628049 &2.025e-322&  9.918e-256&  1.585e-161&  9.329e-099& 0.600000 \\
		800  &-1.376e-004 &  -1.329e-004 & -1.314e-004 & -5.411e-005 & 0.649020 & 1.729e-322& 3.112e-322& 6.175e-322& 1.184e-196& 0.600000 \\
		1600  & -5.020e-004 &  -5.744e-005&  -5.204e-005& -2.709e-005 & 0.605247 & 1.877e-322& 3.013e-322& 4.199e-322& 6.472e-322&0.600000 \\  \bottomrule 
	\end{tabular}} 
	\end{center} 
	\label{tab:tab7}
\end{table}

\begin{figure}[t] 
	\centering
	\includegraphics[scale=0.88]{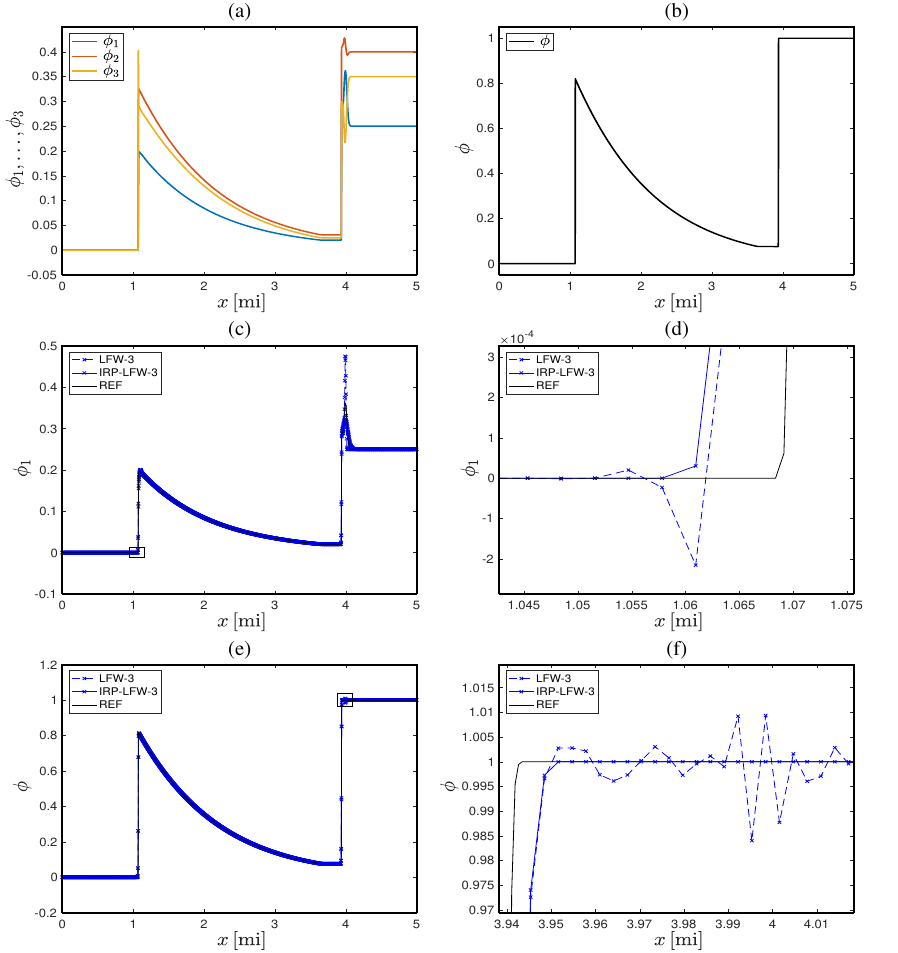}   
		\caption{Example 4 (MCLWR model, $N = 3$): reference solution for (a) $\phi_1, \phi_2, \phi_3$ and (b) 
				 $\phi$   at $T = 0.05 \, \mathrm{h}$  computed by scheme IRP-HLLW-3 with $M_{\mathrm{ref}} = 6400$, 
				  and   comparison of schemes  for (c)  $\phi_3$, (d) enlarged view of~(c), (e) $\phi$, 
				   and (f) enlarged view of~(e)  at $T = 0.05 \, \mathrm{h}$  with $M = 1600$.}
	\label{fig:fig6}
\end{figure}

\begin{figure}[t] 
	\centering
	\includegraphics[scale=0.88]{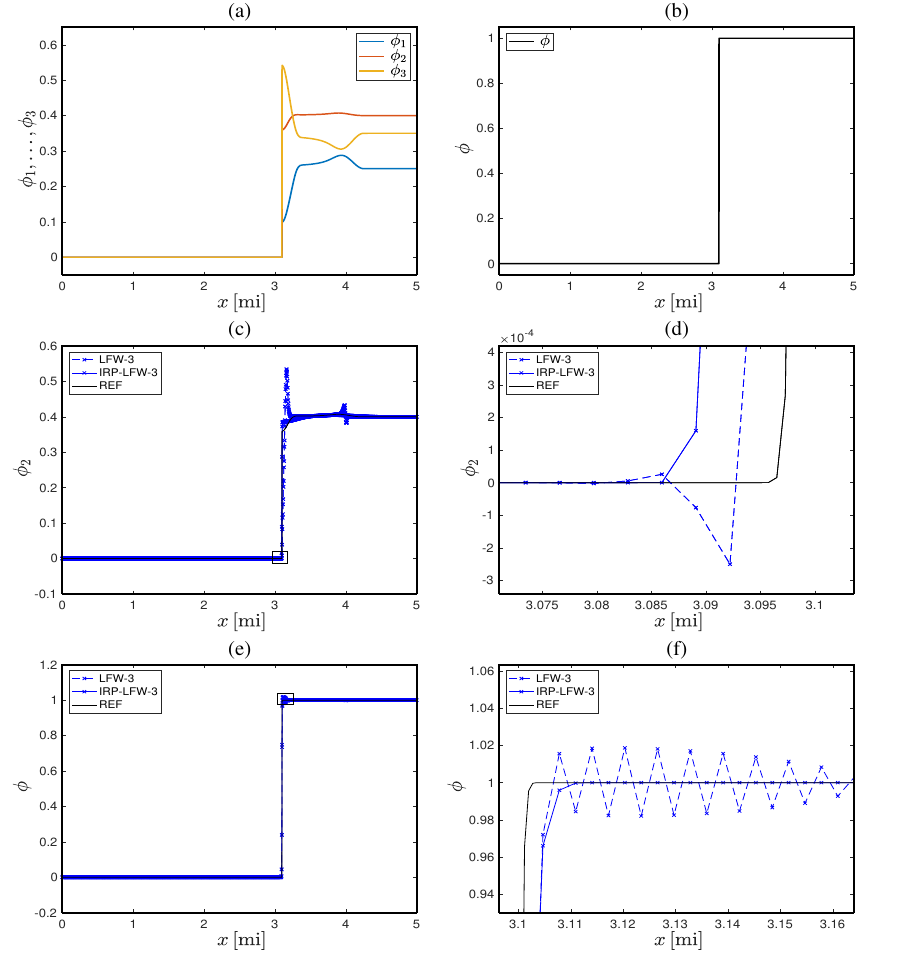}   
		\caption{Example 4 (MCLWR model, $N = 3$): reference solution for (a) $\phi_1, \phi_2, \phi_3$ and (b) 
				 $\phi$   at $T = 0.5 \, \mathrm{h}$  computed by scheme IRP-HLLW-3 with $M_{\mathrm{ref}} = 6400$, 
				  and   comparison of schemes  for (c)  $\phi_2$, (d) enlarged view of~(c),  (e) $\phi$,   
				   and (f) enlarged view of~(e)  at $T = 0.5 \, \mathrm{h}$  with $M = 1600$.}
	\label{fig:fig7}
\end{figure}

\begin{table}[t] \caption{Example 4 (MCLWR model, $N = 3$): $L^1$ errors, numerical order, and CPU time (seconds) for scheme IRP-LFW-3 
  at $T = 0.05 
 \,\mathrm{h}$ and $T = 0.5 \, \mathrm{h}$. The reference solution is computed by IRP-LFW-3 with $M_{\mathrm{ref}} = 6400.$}
	\begin{center} 
		{\footnotesize   \addtolength{\tabcolsep}{-2pt} 
	\begin{tabular}{ccccccc}\toprule
	$M$&$e_M^{\mathrm{tot}}$& $\theta_M$ & cpu [s] & $e_M^{\mathrm{tot}}$ & $\theta_M$& cpu [s] 
	 \\  \cmidrule(lr){2-4}\cmidrule(lr){5-7} 
	& \multicolumn{3}{c}{IRP-LFW-3, $T=0.05 \,\mathrm{h}$} & \multicolumn{3}{c}{IRP-LFW-3, $T=0.5 \, \mathrm{h}$} 	\\ \midrule
		100 & 7.45e-03&  --- & 9.67e-03 & 4.99e-03&  --- & 9.77e-02\\
		200 & 3.51e-03& 1.1& 2.96e-02 &2.67e-03& 0.9 & 2.26e-01 \\
		400 & 1.87e-03& 0.9 &  1.02e-01 & 1.32e-03& 1.0 &  1.05\\
		800 & 8.03e-04&  1.2&   4.23e-01& 6.26e-04&  1.0&  4.07\\
		1600 & 4.04e-04&  1.0&   1.6 & 2.70e-04&   1.2 &  16.2\\
				 \bottomrule
	\end{tabular}} 
	\end{center} 
	\label{tab:tab9}
\end{table}

\begin{table}[t] \caption{Example 4 (MCLWR model, $N = 3$): minimum of the solutions $\phi_{i,j}^n$, $i=1,\dots,3$, and maximum of the solution $\phi_{j}^n$ obtained by schemes LFW-3 (without limiters) and  IRP-LFW-3 (with limiters) until specified times $T$.}
	\begin{center} 
		{\footnotesize  \addtolength{\tabcolsep}{-2pt} 
	\begin{tabular}{ccccccccc}\toprule
	$M$& $\underline{\phi}_{,1}$  & $\underline{\phi}_{,2}$ & $\underline{\phi}_{,3}$  &  $\overline{\phi}$  & $\underline{\phi}_{,1}$   & $\underline{\phi}_{,2}$ & $\underline{\phi}_{,3}$  &  $\overline{\phi}$
	 \\  \cmidrule(lr){2-5}\cmidrule(lr){6-9} 
	& \multicolumn{4}{c}{LFW-3, $T=0.05 \,\mathrm{h}$} & \multicolumn{4}{c}{IRP-LFW-3, $T=0.05 \, \mathrm{h}$}
		\\ \midrule
		100 &   -4.580e-03 &  -4.976e-03  &   -4.914e-03 &  1.01102 & 0.00000& 0.00000& 0.00000& 1.00000\\ 
		200 &    -2.744e-03 &  -3.031e-03 &  -2.636e-03 & 1.01537 &   0.00000& 0.00000&  0.00000& 1.00000\\ 
		400  & -1.759e-03 &  -1.934e-03 & -1.891e-03 & 1.04909  & 0.00000& 0.00000& 0.00000& 1.00000\\ 
		800  &  -1.117e-03& -1.531e-03 & -1.344e-03 & 1.05402 &0.00000& 0.00000& 0.00000& 1.00000\\ 
		1600  & -9.064e-04 & -1.058e-03 & -9.975e-04 &  1.06708 & 0.00000& 0.00000& 0.00000& 1.00000\
				\\ \midrule 
				 & \multicolumn{4}{c}{LFW-3, $T=0.5 \,\mathrm{h}$} & \multicolumn{4}{c}{IRP-LFW-3, $T=0.5 \, \mathrm{h}$}
		\\ \midrule
  100 & -5.041e-03 &  -5.887e-03 &  -7.050e-03 & 1.04260 & 0.00000& 0.00000& 0.00000& 1.00000\\ 
  200 & -3.446e-03 &  -3.031e-03 & -3.518e-03 & 1.07077 & 0.00000& 0.00000& 0.00000& 1.00000\\
  400  &-2.447e-03 & -2.711e-03 & -1.891e-03 &   1.06805 & 0.00000& 0.00000& 0.00000& 1.00000\\
  800  &-2.137e-03 & -2.499e-03 & -1.344e-03 & 1.09395 & 0.00000& 0.00000& 0.00000& 1.00000\\
  1600  & -1.779e-03 & -1.225e-03 & -9.975e-04 & 1.09029 & 0.00000& 0.00000&  0.00000& 1.00000\\ \bottomrule 
	\end{tabular}} 
	\end{center} 
	\label{tab:tab10}
\end{table}

\begin{figure}
	\centering
	\includegraphics[scale=0.88]{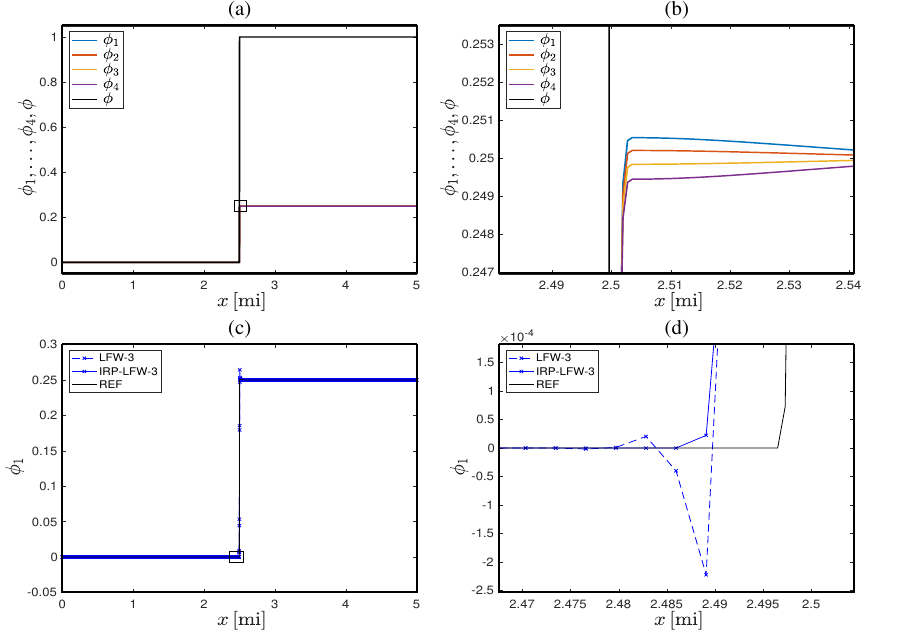} 
		\caption{Example 5 (MCLWR model, Daganzo's test, $N = 4$): 
		reference solution for (a) $\phi_1, \dots,  \phi_4$ and $\phi$, (b) enlarged view of~(a)
				   at $T = 0.5 \, \mathrm{h}$  computed by scheme IRP-HLLW-3 with $M_{\mathrm{ref}} = 6400$, 
				  and   comparison of schemes  for (c)  $\phi_1$, (d) enlarged view of~(c),  (e) $\phi$,   
				   and (f) enlarged view of~(e)  at $T = 0.5 \, \mathrm{h}$  with $M = 1600$.}
	\label{fig:fig8}
\end{figure}

\begin{table}[t] \caption{Example 5 (MCLWR model, Daganzo's test, $N = 4$): minimum of the solutions $\phi_{i,j}^n$, $i=1,\dots,4$, and maximum of the solution $\phi_{j}^n$ obtained by schemes LFW-3 (without limiters) and  IRP-LFW-3 (with limiters) until a specified time $T$.}
	\begin{center} 
		{\footnotesize  \addtolength{\tabcolsep}{-2pt} 
	\begin{tabular}{ccccccccccc}\toprule
	$M$& $\underline{\phi}_{,1}$  & $\underline{\phi}_{,2}$  &  $\underline{\phi}_{,3}$  & $\underline{\phi}_{,4}$  &   $\overline{\phi}$ &  $\underline{\phi}_{,1}$  & $\underline{\phi}_{,2}$  &  $\underline{\phi}_{,3}$  & $\underline{\phi}_{,4}$  &   $\overline{\phi}$
	 \\  \cmidrule(lr){2-6}  \cmidrule(lr){7-11} 
	& \multicolumn{5}{c}{LFW-3, $T=0.5 \,\mathrm{h}$}  & \multicolumn{5}{c}{IRP-LFW-3, $T=0.5 \,\mathrm{h}$}  \\ \midrule
		100 &-2.186e-003&  -2.673e-003 &  -3.177e-003 &  -3.691e-003 & 1.00403 & 0.00000&  0.00000& 0.00000&  0.00000& 1.00000 \\ 
		200 & -1.832e-003 &  -1.938e-003 &  -1.987e-003 &  -1.953e-003 & 1.01107 &  0.00000&  0.00000& 0.00000&  0.00000& 1.00000 \\
		400  &-9.020e-004 &  -6.778e-004 &   -6.999e-004 & -8.780e-004 & 1.00369 &0.00000&  0.00000& 0.00000&  0.00000& 1.00000 \\
		800  &-3.076e-004 & -3.848e-004 & -4.523e-004 & -4.882e-004 & 1.00107 & 0.00000&  0.00000& 0.00000&  0.00000& 1.00000 \\
		1600  & -2.266e-004 &  -2.356e-004 & -2.084e-004 & -1.960e-004 & 1.00122  & 0.00000& 0.00000& 0.00000& 0.00000& 1.00000	 
		 				 \\ \bottomrule \end{tabular}} 
	\end{center} 
	\label{tab:tab12}
\end{table}

\subsection{Preliminaries}  We discretize the domain $[0,L]\times [0,T]$ 
as outlined in Section~\ref{subsec3.1}. For the MLB model, we employ the first-order LLF  flux given by \eqref{eq:llf}  and the 
  HLL flux, defined by~\eqref{eq:hll}.  For the  MCLWR model,  we only use the  LLF method.   {For computing the nonlinear weights \eqref{eq:nonl-weights} in the WENO reconstruction procedure, we have set $p=2$ in all numerical test and we have used $\varepsilon = \Delta x^2$ and $\varepsilon = 10^{-6}$ in the first example (accuracy test) to show the robustness of the WENO scheme for both choices of the parameters $\varepsilon$ while in the rest of examples we use $\varepsilon = \Delta x^2$. The remaining parameters of the reconstructions are taken as usual (see \cite{levy2000compact,cravero2018cweno})}. In the first example, we compare numerical results obtained by the invariant region preserving WENO schemes of orders~$3$ and~$5$. For the rest of examples,  we 
   limit ourselves to third-order  WENO reconstructions. We denote by ``IRP-LFW-3'' and ``IRP-HLLW-3'' 
     the  LLF and HLL methods, respectively,  with third-order IRP  WENO reconstructions. Analogously,   ``IRP-LFW-5'' and ``IRP-HLLW-5''  denote the respective fifth-order versions. 

For comparison purposes, we compute reference solutions for numerical tests by the IRP-WENO-5 scheme with $M_{\mathrm{ref}} = 12800$ in Example 1, by IRP-HLLW-3  with $M_{\mathrm{ref}} = 6400$ in Examples 2 and 3, and by IRP-LFW-3 with $M_{\mathrm{ref}} = 6400$ in Examples 4 and 5. As in \cite{burger2011implementation,burger2016polynomial}, we compute approximate $L^1$ errors at different times for each scheme as follows. We denote by $(\phi_{i,j}^{n}(t))_{j=1}^{M}$ and $(\phi_{i,j}^{\mathrm{ref}}(t))_{j=1}^{M_{\mathrm{ref}}}$ the numerical solution for the $i-$th component at time $t$ calculated with $M$ and $M_{\mathrm{ref}}$ cells, respectively. We compute $\smash{\tilde{\phi}_{i,j}^{\mathrm{ref}}(t
)}$, for $j=1,\dots,M$,  by
\begin{align*}
\tilde{\phi}_{i,j}^{\mathrm{ref}}(t) = \dfrac{1}{R} \sum_{k=1}^{R} \phi_{i,R(j-1)+k}^{\mathrm{ref}} (t),\quad R=M_{\mathrm{ref}}/M.
\end{align*}

The total approximate $L^1$ error of the numerical solution on the $M$-cell grid at time $T$ is then given by
\begin{equation}\label{eq:err-tot}
	e_M^{\mathrm{tot}}(t) \coloneqq  \dfrac{1}{M}\sum_{i=1}^{N}\sum_{j=1}^{M} \big|\tilde{\phi}_{i,j}^{\mathrm{ref}}(t)-\phi_{i,j}^{M}(t)\big|.
\end{equation}
We may then calculate a numerical order of convergence from~$e_M^{\mathrm{tot}}(t)$ and~$e_{2M}^{\mathrm{tot}}(t)$ by
\begin{equation}\label{eq:conv-rate}
	\theta_M(t)\coloneqq  \log_2 \big(e_M^{\mathrm{tot}}(t)/e_{2M}^{\mathrm{tot}}(t)\big).
\end{equation}

 In order to study the effect of the limiters \eqref{eq:ptilde} and \eqref{eq:phat}, we compute the minimum values of the numerical solution in space and time for each component and the maximum values of the numerical solution in space and time for the total concentration with and without the limiters, i.e. we compute: 
 \begin{equation*}
 	 \underline{\phi}_{,i}\coloneqq\min\limits_{(j,n)\in \mathbb{Z}_M\times \mathbb{Z}_T}\{\phi_{i,j}^{n}\}, \quad i=1,\dots, N,\quad \text{ and }\quad \overline{\phi}\coloneqq\max\limits_{(j,n)\in \mathbb{Z}_M\times \mathbb{Z}_T}\{\phi_{j}^{n}\},
 \end{equation*}
 where $\mathbb{Z}_{T}\coloneqq\{0,\dots,N_T\}$.
 
\subsection{Example 1: MLB model, $N=3$, numerical order of accuracy}

In this example, we consider an experiment performed in \cite{boscarino2016linearly} to verify numerically  the convergence rate of the IRP-LFW-3,  IRP-HLLW-3, IRP-LFW-5, and IRP-HLLW-5 numerical schemes. We employ the MLB model introduced in Section \ref{sec:models} with $N=3$ species having normalized squared particle sizes $\boldsymbol{\delta} = (1, 0.8, 0.6)^{\mathrm{T}}$ 
  with density $\rho_{\mathrm{s}} = 2790\, \mathrm{kg}/\mathrm{m}^3$ and vessel height $L = 1\, \mathrm{m}$. The maximum total concentration is $\phi_{\mathrm{max}} = 0.66$. The hindered settling factor $V(\phi)$ is chosen according to \eqref{eq:hind-set-fun-2} with the exponent $n_{\mathrm{RZ}} = 4.7$. The remaining parameters are $g = 9.81 \, \mathrm{m}/\mathrm{s}^2$, $\mu_{\mathrm{f}} = 0.02416$ Pas and $\rho_{\mathrm{f}} = 1208\, \mathrm{kg}/\mathrm{m}^3$. We choose a smooth initial concentration profile given by 
    $\Phi_0(x) =  0.12 \exp(-200(x - 0.5)^2) (0.12,0.12,0.12)^{\mathrm{T}}$. We compute approximate solutions with $\Delta x = L/M$, and  $M =100 \cdot 2^{\ell}$,
$ \ell = 0\dots, 5,$ and a fixed time step $\Delta t = 50 \Delta x$, for the third-order method. Since the order of TVD Runge-Kutta time step solver is 3, in the case of fifth-order schemes,  we use a time step $\Delta t$ such that $(\Delta t)^3$ is maintained proportional
to $\Delta x^5$, i.e we just set $\Delta t \propto \Delta x^{5/3}$, in order to get the correct convergence rate. Figure \ref{fig:ex1-1} shows the numerical results for $M = 1600$ at $T = 5$ s (before shock formation, when the solution is still
smooth) and for $T = 10$ s (after shock formation).

The approximate $L^1$ errors $e_{M}^{\mathrm{tot}}(T)$ defined by \eqref{eq:err-tot} and their corresponding numerical orders $\theta_M (T)$ given by \eqref{eq:conv-rate} are displayed in Table~\ref{tab:tab1}  at times $T = 5$ s and
$T = 10$ s for both schemes. The reference solution is computed with $M_{\mathrm{ref}} = 12800$ cells by using the fifth-order scheme  IRP-HLLW-5. The
behavior of $\theta_M (5)$ for increasing values of $M$ confirms third-order convergence for smooth solutions for IRP-LFW-3,  IRP-HLLW-3, and fifth order for schemes IRP-LFW-5,  IRP-HLLW-5. The
results for $T = 10$ s indicate that accuracy is reduced to first order when shocks are present, as expected.

\subsection{Example 2: MLB model, $N=2$}\label{numex2}

This example corresponds to $N=2$ species  with density $\rho_{\mathrm{s}} = 2790\, \mathrm{kg}/\mathrm{m}^3$ and different diameters $D_1 = 4.96 \times 10^{-4}\, \mathrm{m}$ and $D_2 = 1.25 \times 10^{-4}\, \mathrm{m}$, such that  $\boldsymbol{\delta} = (1,0.063)^{\mathrm{T}}$. The (unnormalized) depth of the vessel in the original experiment is $L = 0.3\, \mathrm{m}$ \cite{schneider1985sediment}. The maximum total concentration is $\phi_{\mathrm{max}} = 0.6$, and the initial concentrations are $\Phi_0(x) = (0.2,0.05)^{\mathrm{T}}$. The remaining parameters are taken as in Example 1. 
The well-known solution of Example 2   has been used as a test case for a variety of methods  \cite{burger2007adaptive,burger2011implementation,burger2016polynomial}. For comparison purposes, we calculate numerical solutions for a sequence of spatial discretizations $\Delta x = L/M$ with IRP-LFW-3 and IRP-HLLW-3 schemes and compare the solutions with a reference solution  
with $M = M_{\mathrm{ref}} = 6400$ obtained by  the   IRP-HLLW-3 scheme. The reference solution is shown in Figures~\ref{fig:fig2} 
 and~\ref{fig:fig3}    for the simulated times $T = 50\,  \mathrm{s}$  and $T = 300\, \mathrm{s}$, respectively.  On the other hand, we plot the solutions of the total concentration $\phi$ with $M=2^{\ell}\cdot 100$, $\ell =0,1,\dots,3$ and the reference solution computed with $M=6400$ mesh points 
 in Figure~\ref{fig:fig9}.

In Table\ref{tab:tab3} we show approximate $L^1$ errors and CPU times for
both schemes at two selected simulated times, these approximate errors are computed by \eqref{eq:err-tot} and we observe the convergence of both methods.  The minimum values $\smash{\underline{\phi}_{,i}}$ for $i=1,2$ and the maximum value $\smash{\overline{\phi}}$, with and without the limiters are presented in Table~\ref{tab:tab4}. In the case of schemes without limiters we observe some negative values of $\smash{\underline{\phi}_{,i}}$ and some values of $\smash{\overline{\phi}}$ greater than $\phi_{\mathrm{max}}$ due to overshoots present in the numerical solution, see for example Figure \ref{fig:fig3} (d), while in the case of schemes with limiters the numerical solution belongs to $\mathcal{D}_{\phi_{\mathrm{max}}}$, as expected.

\subsection{Example 3: MLB model, $N=4$}

In this example we consider $N=4$  particle sizes with $D_1=4.96 \times 10^{-4}\, \mathrm{m}$, the rest of diameters $D_i$, $i=2,3,4$ are chosen such that $\boldsymbol{\delta}= (1,0.8,0.6,0.4)^{\mathrm{T}}$ and we set  $\Phi_0(x)=(0.05,0.05,0.05,0.05)^{\mathrm{T}}$. The other parameters are   as in Example~2. This example goes back to Greenspan and Ungarish \cite{greenspan1982hindered}, and was solved numerically in \cite{burger2000model} with the slightly different hindered settling factor $V(\phi) = (1-(5/3)\phi)^{2.7}$. The reference solution is shown in Figures~\ref{fig:fig4} 
 and~\ref{fig:fig5}   for the simulated times $T = 50\, \mathrm{s}$   and $T = 300\, \mathrm{s} $, respectively. In Table \ref{tab:tab6} we show approximate $L^1$ errors and CPU times for
both schemes at two selected simulated times and we observe the convergence of both methods. These approximate errors are computed by \eqref{eq:err-tot}.The minimum values $\smash{\underline{\phi}_{,i}}$ for $i=1,\dots,4$ and the maximum value $\smash{\overline{\phi}}$, with and without the limiters are presented in Table~\ref{tab:tab7}. In the case of schemes without limiters we observe some negative values of $\smash{\underline{\phi}_{,i}}$ and some values of $\smash{\overline{\phi}}$ greater than $\phi_{\mathrm{max}}$ due to overshoots present in the numerical solution, see for example Figure \ref{fig:fig5}(d), while in the case of schemes with limiters the numerical solution belongs to $\mathcal{D}_{\phi_{\mathrm{max}}}$, as expected.

\subsection{Example 4: MCLWR model, $N=3$}

We now study the  MCLWR model on  a road of length $L = 5$ mi (mi stands for miles) with $N = 3$ driver classes associated with $\beta_1 = 60\, \mathrm{mi}/\mathrm{h}$ , $\beta_2 = 55\, \mathrm{mi}/\mathrm{h}$ , and
$\beta_3 = 50\, \mathrm{mi}/\mathrm{h}$ . We employ the Dick–Greenberg model \eqref{eq:dgmodel} with $\phi_{\mathrm{max}}=1$ and choose (as in \cite{boscarino2016linearly}) 
$\phi_{\mathrm{c}}  = \exp(-7/\mathrm{e}) \approx 0.076142$. The velocity is then given by
\begin{align*}
	\begin{cases}
			v(\phi)  = 1, \, v'(\phi ) =0 &\text{for $0\leq \phi \leq \phi_{\mathrm{c}} = \exp(-1/C)\approx 0.076142$,} \\
		v(\phi)  = -C \ln (\phi), \, v'(\phi) = -C/\phi &\text{for $\phi_{\mathrm{c}}<\phi < 1$.} 
		\end{cases}
\end{align*}
 The initial density distribution is given by  $\Phi_0(x) =  p(x)(0.25, 0.4, 0.35)^{\mathrm{T}}$, where $p$ describes an isolated platoon on $1\leq x\leq 2$ followed by a constant maximum density function for $x\geq 4$, i.e.,    
\begin{align*}
	p(x) \coloneqq 
	\begin{cases}
		10(x-1)&\text{for $1\leq x\leq 1.1$,} \\
		1&\text{for $1.1\leq x\leq 1.9$,} \\
		-10(x-2)&\text{for $1.9\leq x\leq 2$,} \\
		1& \text{for $x\geq 4$,}\\
		0& \text{otherwise.}  
	\end{cases}
\end{align*}
 For this particular case we use zero flux boundary conditions in order to obtain a steady state solution. The reference solution is shown in Figure \ref{fig:fig6} (a)-(b) and Figure \ref{fig:fig7} (a)-(b), for the simulated times $T = 0.05$ h and $T = 0.5$ h, respectively. In Table \ref{tab:tab9} we show approximate $L^1$ errors and CPU times for the scheme at the two selected simulated times and we observe convergence of the method. These approximate errors are computed by \eqref{eq:err-tot}. In the case of the scheme without limiters we observe some negative values of 
  $\smash{\underline{\phi}_{,i}}$ and some values of $\smash{\overline{\phi}}$ greater than $\phi_{\mathrm{max}}$ due to overshoots present in the numerical solution, see for instance Figures~\ref{fig:fig7}(c) and~(d).

\subsection{Example 5: Daganzo’s test, $N=4$}

In this subsection,  we study a test that Daganzo suggested in \cite{daganzo1995requiem}. The parameters and boundary conditions are the same as in Example 4 with densities $\beta_1 = 60\, \mathrm{mi}/\mathrm{h}$ , $\beta_2 = 55\, \mathrm{mi}/\mathrm{h}$ , 
$\beta_3 = 50\, \mathrm{mi}/\mathrm{h}$ , and $\beta_4 = 45\, \mathrm{mi}/\mathrm{h}$ . In order to conduct a multiclass test that is appropriate for this scenario, we follow the ideas of Bürger et. al in \cite{burger2013antidiffusive} and we consider the initial condition 
\begin{align*}
	\Phi(x,0)=\begin{cases}
		\Phi_{\mathrm{L}} &\text{for $x<L/2$,} \\
		\Phi_{\mathrm{R}} &\text{for  $x\geq L/2$,} \\
	\end{cases}
\end{align*}
where $\Phi_{\mathrm{L}}=(0,0,0,0)^{\mathrm{T}}$ and $\Phi_{\mathrm{R}}=(0.25,0.25,0.25,0.25)^{\mathrm{T}}$. This density distribution
should be a stationary solution for the model. The reference solution is shown in Figures~\ref{fig:fig8}(a) and~(b) for the simulated time 
 $T = 0.5 \, \mathrm{h}$. The minimum values $\smash{\underline{\phi}_{,i}}$ for $i=1,\dots,4$ and the maximum value $\smash{\overline{\phi}}$, with and without the limiters are presented in Table \ref{tab:tab12}. See also Figures~\ref{fig:fig8}(c)--(f).

\section{Conclusions}  \label{sec:conc} 
In this work, we have designed high-order finite volume numerical schemes for multiclass kinematic flow problems, including polydisperse sedimentation and multiclass vehicular traffic models, whose numerical solution preserves  an  invariant region which corresponds to the space $\mathcal{D}_{\phi_{\max}}$ of physically relevant solutions of the model. The first contribution of the paper is the proof that first order schemes 
with the LLF  numerical flux given by \eqref{eq:llf} and the HLL numerical flux given by \eqref{eq:hll} are invariant region preserving (IRP) for the MLB polydisperse sedimentation model, under some appropriate CFL conditions \eqref{llfcfl} and \eqref{hllcfl}, respectively. The key part of the proof was to assume the conditions \eqref{mlbprop} and to consider the slightly lower bound for the eigenvalues of the model given by \eqref{sl34}. We also proof that these schemes are IRP for MCLWR model under the assumptions \eqref{mclwrprop}. One can also use  the  HW scheme (Scheme 4 in \cite{burger2008family}) which  satisfies  Theorem \ref{thm:lfinv} but happens to be too dissipative. Next, by considering a component-wise WENO reconstruction, we use a modification of Zhang and Shu's scaling limiter \cite{zhang2010maximum} which consists in two main steps: first, we define a linear scaling limiter for each species to get positive solutions in each component and then we consider a second linear scaling limiter in such a way that the total concentration is bounded by $\phi_{\mathrm{max}}$. Finally, we apply a strong-stability preserving (SSP) third-order TVD Runge-Kutta time discretization to obtain the fully-discrete scheme. With this one can show that, under a more restrictive CFL condition, the resulting finite volume scheme is IRP and the high order of accuracy is not destroyed.

The numerical results obtained support the theoretical findings. In all the numerical simulations the minimum values of the solution for each species and the maximum value of the total concentration are  computed and tabulated and the effect of the limiters is appreciated in both the tables and the plots of the solutions. As future perspectives, we wish to apply this strategy to multiclass vehicular dynamics with uneven space occupancy or with creeping, where some physical invariant regions different to $\mathcal{D}_{\phi_{\mathrm{max}}}$ are studied \cite{briani2021macroscopic,fan2015heterogeneous}. Moreover, we 
 intend to explore the extension of these modified scaling limiters to nonlocal two-dimensional multispecies models. 
   {On the other hand, we mention that   multispecies kinematic flow models  other than the MCLWR or MLB 
   models could  be handled by the present numerical method as well provided that  the velocity functions $v_i ( \Phi)$ have the properties 
    outlined in Section~\ref{subsec:scope}. For instance, kinematic models of polydisperse sedimentation of the form \eqref{eq:main-eq}
     but with velocity functions~$v_i ( \Phi)$ alternative to the MLB model are discussed in \cite{burger2010hyperbolicity,burger2011implementation}. We also mention  the model of gravity-driven separation of oil-water dispersions introduced in \cite{rosso2001gravity} whose 
     velocity functions are equivalent to those of the MCLWR model but which is equipped with zero-flux boundary conditions.} 

 {Finally, we comment that there are no fundamentally difficulties in appying similar WENO schemes  to two- or higher-dimensional systems of conservation laws that are multi-dimensional analogues of   \eqref{eq:main-eq}, for instance in two
 dimensions, 
\begin{align} \label{eq6.1} 
 \partial_t \Phi +  \partial_x \boldsymbol{f}^x ( \Phi) + \partial_y \boldsymbol{f}^y ( \Phi ) = \boldsymbol{0}, 
\end{align} 
with functions~$\boldsymbol{f}^x$  and~$\boldsymbol{f}^y$ that have properties similar to those of~$\boldsymbol{f}$ 
  in   \eqref{eq:main-eq}, and in particular are compactly supported on~$\mathcal{D}_{\phi_{\max}}$. However, 
   models of polydisperse sedimentation require in two or more space dimensions the computation of 
    a volume average flow field, which only in one space dimension is given by boundary conditions, 
     and for the case of batch settling considered herein vanishes at all (cf., e.g.,
      \cite{buerger2002model} for a complete formulation). That flow field  gives rise to an additional advection term 
      in the system of  conservation laws governing the evolution of the concentration vector, and needs to be computed 
        by a variant of the Stokes or Navier-Stokes equations. The flow is, moreover, driven by density fluctuations of the 
         mixture, that is, by fluctuations of concentrations. Thus, the two-dimensional analogue of the MLB model 
          is not just a multi-dimensional version of \eqref{eq:main-eq} but a coupled flow-transport problem, 
           which is currently being investigated by the authors. 
    }  

\section*{Acknowledgments}
 JBC, RB and LMV are supported by ANID (Chile) through Anillo project ANID/PIA/ACT210030 and   Centro de Modelamiento Matem\'{a}tico (CMM), project 
 FB210005 of BASAL funds for Centers of Excellence. RB is also supported by CRHIAM, projects ANID/Fondap/15130015 and ANID/Fondap/1523A0001 and Fondecyt project 1250676. JBC is supported by the National Agency for Research and Development, ANID-Chile through Scholarship Program, Beca Doctorado Nacional 2022, folio 21221387.  PM is supported by PID2020-117211GB-I00 and PID2023-146836NB-I00, granted by MCIN/ AEI /10.13039/ 501100011033, and 
CIAICO/2021/227, granted by GVA.


\appendix
\section{Interlacing property}
\label{app1}
Multispecies kinematic flow models are given by systems of conservation laws with fluxes as in \eqref{eq:main-eq} for sufficiently smooth functions $v_i$. However, in most applications,  these functions do not depend individually on each of the densities $\phi_1,\dots,\phi_N$ but rather on a small number $m\ll N$  of functions $p_1,\dots,p_m$ of $\phi_1,\dots,\phi_N$, as we presented in Section \ref{subsec:mlb} and \ref{subsec:mclwr}, i.e, we have
\begin{equation*}
	v_i = v_1(p_1,\dots,p_m),\quad p_{\ell} = p_{\ell}(\phi_1,\dots,\phi_N), \quad i=1,\dots,N,\quad \ell =1,\dots, m, \quad m\ll N.
\end{equation*}
Under this assumption we can write the Jacobian matrix $\boldsymbol{\mathcal{J}}_{\boldsymbol{f}} = \boldsymbol{\mathcal{J}}_{\boldsymbol{f}}(\Phi)$ in the form
\begin{align}
	\boldsymbol{\mathcal{J}}_{\boldsymbol{f}} &= \boldsymbol{D}+\boldsymbol{B}\boldsymbol{C}^{\mathrm{T}}, \quad \boldsymbol{D}\coloneqq  \diag(v_1,\dots,v_N),  \label{eq:jac}\\
	\boldsymbol{B}&\coloneqq  (b_{i\ell}) = \left(\phi_i \dfrac{\partial v_i}{\partial p_{\ell}}\right),\quad \boldsymbol{C}\coloneqq  (c_{j \ell}) = \left(\dfrac{\partial p_{\ell}}{\partial \phi_{j}}\right),\quad 1\leq i,j\leq N, \quad 1\leq \ell\leq m. \label{eq:matrices}
\end{align}
The following Theorem has been proved in \cite{anderson1996secular} and it is used to show hyperbolicity of selected multispecies kinematic flow models and for the construction of spectral numerical schemes.

\begin{theorem}[Secular equation \cite{donat2010secular,anderson1996secular}]\label{teo:seceq} Assume that $\boldsymbol{D}$ is a diagonal matrix as given by \eqref{eq:jac} with $v_i>v_j$ for $i<j$ and that $\boldsymbol{C}$ and $\boldsymbol{B}$ have the formats specified in \eqref{eq:matrices}. Let $\lambda\neq v_i$ for $i=1,\dots, N$. Then $\lambda$ is an eigenvalue of $\boldsymbol{D}+ \boldsymbol{BC}^{\mathrm{T}}$ if and only if
	\begin{equation*}
		R(\lambda)\coloneqq \det(M_{\lambda}) = 1+\sum_{i=1}^{N} \dfrac{\gamma_i}{v_i-\lambda} = 0.
	\end{equation*}
	The coefficients $\gamma_1,\dots, \gamma_N$ are given by the following expression, where $I\coloneqq \{i_1<\dots<i_k\}\in S_k^{N}$ and $J\coloneqq \{j_1<\dots<j_{\ell}\}\in S_{\ell}^{m}$ are index sets:
	\begin{align*}
		\gamma_i = \sum_{r=1}^{\min\{N,m\}} \sum_{i\in I\in S_r^{N}, J\in S_r^{m}} \dfrac{\det \boldsymbol{C}^{I,J}}{\prod_{\ell \in I, \ell \neq i} (v_{\ell}-v_i)}.
	\end{align*}	
\end{theorem}
\begin{corollary}[Interlacing property \cite{burger2010hyperbolicity}] \label{cor:inter-prop}  With the notation of the theorem above, assume that $\gamma_i \gamma_j>0$ for $i,j = 1,\dots,N$. Then $\boldsymbol{D}+ \boldsymbol{BC}^{\mathrm{T}}$ is diagonalizable with real eigenvalues $\lambda_1,\dots,\lambda_N$. If $\gamma_1,\dots,\gamma_N<0$, the interlacing property 
	\begin{equation}\label{eq:int-prop}
		M_1\coloneqq  v_N+\gamma_1+\cdots+\gamma_N<\lambda _N<v_N<\lambda_{N-1}<\cdots<\lambda_1<v_1 
	\end{equation}
	holds, while for $\gamma_1,\dots,\gamma_N>0$, the following analogous property holds:
	\begin{equation}
		v_N<\lambda_{N}<v_{N-1}<\lambda_{N-1}<\cdots<v_1<\lambda_1<M_2\coloneqq v_1+\gamma_1+\cdots+\gamma_N.
	\end{equation}
\end{corollary}
\begin{remark}
	Theorem \ref{teo:seceq} and Corollary \ref{cor:inter-prop} apply to the MCLWR and  MLB 
	   models (see \eqref{eq:bounds-MCLWR} and  \eqref{eq:bounds-MLB} for the respective specific bounds).
\end{remark}

\end{document}